\documentclass[twoside]{article}

\usepackage{epsfig}
\input{diagrams}
\diagramstyle[repositionpullbacks,midshaft,PostScript=dvips,nohug]

\newcommand{\mcm}[3]{\newcommand{#1}[#2]{{\ensuremath{#3}}}}

\usepackage{latexsym}
\usepackage{amssymb}

\mcm{\blank}{0}{(\emptybk)}
\mcm{\emptybk}{0}{\:\:}
\mcm{\hyph}{0}{\mbox{-}}

\mcm{\diagspace}{0}{\mbox{\hspace{2em}}}

\newcommand{\bref}[1]{(\ref{#1})}

\mcm{\cat}{1}{\mc{#1}}
\mcm{\fcat}{1}{\mb{#1}}
\mcm{\mc}{1}{\mathcal{#1}}
\mcm{\mr}{1}{\mathrm{#1}}
\mcm{\mi}{1}{\mathit{#1}}
\mcm{\mb}{1}{\mathbf{#1}}
\mcm{\scat}{1}{\Bbb{#1}}
\mcm{\twid}{1}{\widetilde{#1}}

\mcm{\elt}{0}{\in}
\mcm{\sub}{0}{\,\subseteq\,}
\mcm{\such}{0}{\:|\:}
\mcm{\without}{0}{\setminus}

\mcm{\iso}{0}{\,\cong\,}

\newcommand{\epsln}{\varepsilon}

\mcm{\blob}{0}{\scriptscriptstyle{\bullet}}

\mcm{\End}{0}{\mr{End}}
\mcm{\Hom}{0}{\mr{Hom}}
\mcm{\ob}{0}{\mr{ob}\,}
\mcm{\op}{0}{\mr{op}}

\mcm{\comp}{0}{\mi{comp}}
\mcm{\id}{0}{\mi{id}}
\mcm{\ids}{0}{\mi{ids}}
\mcm{\mult}{0}{\mi{mult}}
\mcm{\unit}{0}{\mi{unit}}

\mcm{\Cat}{0}{\fcat{Cat}}
\mcm{\Graph}{0}{\fcat{Graph}}
\mcm{\One}{0}{\fcat{1}}
\mcm{\Set}{0}{\fcat{Set}}

\mcm{\pr}{2}{\tuplebts{#1,#2}}

\mcm{\ftrcat}{2}{[#1,#2]}

\mcm{\go}{0}{\rTo}
\mcm{\goby}{1}{\rTo^{#1}}
\mcm{\goesto}{0}{\,\longmapsto\,}
\mcm{\goiso}{0}{\goby{\diso}}
\mcm{\og}{0}{\lTo}
\mcm{\ogby}{1}{\lTo^{#1}}

\mcm{\bktdvslob}{3}
	{\left(
	\begin{diagram}[height=1.5em]
	#1		\\
	\dTo>{\,#2}	\\
	#3		\\
	\end{diagram}
	\right)}

\newarrow{Edge}---->		
\newarrow{Equals}=====
\newarrow{Get}....>
\newarrow{Goesto}|---{->}
\newarrow{Monic}{vee}---{vee}

\newenvironment{tree}
	{\begin{diagram}[height=1em,width=.75em,abut,noPS,tight]}	
	{\end{diagram}}

\newcommand{\dn}{\dLine}
\mcm{\enode}{0}{\circ}
\newcommand{\lt}[1]{\ldLine(#1,2)}
\mcm{\nl}{1}{\stackrel{\textstyle #1}{\node}}
\mcm{\node}{0}{\bullet}
\newcommand{\rt}[1]{\rdLine(#1,2)}
\mcm{\utree}{0}{\node}

\mcm{\diso}{0}{\sim}
\newcommand{\pullshape}
	{\setlength{\unitlength}{1em}
	\begin{picture}(2,5)(-1,-5)
	\put(0,-5){\line(1,1){1}}
	\put(0,-5){\line(-1,1){1}}
	\end{picture}}
\newcommand{\Spbk}{\overprint{\raisebox{-2.5em}{\pullshape}}}
\mcm{\vdiso}{0}{\wr}

\newcommand{\piccy}[1]{\epsfig{file=#1}}

\mcm{\tuplebts}{1}{(#1)}
\mcm{\ovln}{1}{\overline{#1}}
\newarrow{Incl}C--->
\mcm{\tr}{0}{\fcat{tr}}
\mcm{\ladj}{0}{\,\dashv\,}

\mcm{\Top}{0}{\fcat{Top}}
\mcm{\act}{1}{\mi{act}_{#1}}
\mcm{\TCat}{1}{\fcat{Wk}\hyph#1\hyph\Cat}
\mcm{\littletree}{0}{\raisebox{-.3em}{\epsfig{file=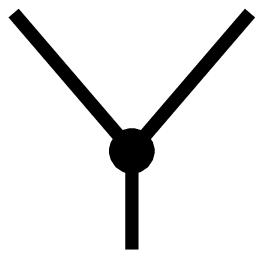,height=1.1em}}}
\mcm{\lefttree}{0}{\raisebox{-.3em}{\epsfig{file=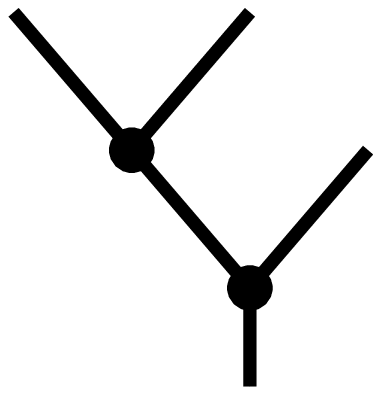,height=1.8em}}}
\mcm{\righttree}{0}{\raisebox{-.3em}{\epsfig{file=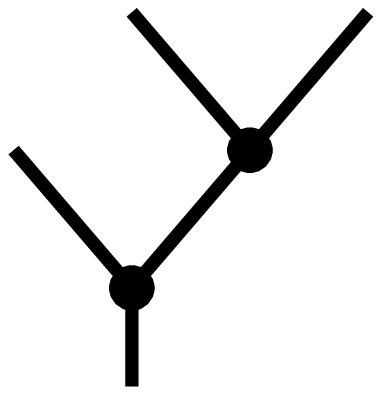,height=1.8em}}}

\mcm{\dashbk}{0}{-}

\newarrow{Labelling}----{triangle}
\mcm{\bof}{0}{\raisebox{0.08ex}{\ensuremath{\scriptstyle\bullet}}}
\mcm{\bofdim}{1}{\,\bof_{#1}\,}
\mcm{\vslob}{3}
	{\left(
	\begin{diagram}[height=1.5em]
	#1		\\
	\dTo>{\,#2}	\\
	#3		\\
	\end{diagram}
	\right)}
\mcm{\implies}{0}{\ \Rightarrow\ }

\mcm{\cod}{0}{\mi{cod}}
\mcm{\dom}{0}{\mi{dom}}
\mcm{\reason}{0}{\,\Rightarrow\,}
\mcm{\Reason}{0}{\,\Rrightarrow\,}
\newcommand{\Ddownarrow}{\epsfig{file=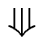}}

\mcm{\intvl}{1}{\langle #1 \rangle}
\mcm{\lnode}{1}{\node\makebox[0em]{$\,\,\scriptstyle{#1}$}} 
\mcm{\lanode}{1}{\node\raisebox{1.2ex}[0ex][0ex]%
{\makebox[0em]{$\!\scriptstyle{#1}$}}} 
\mcm{\lenode}{1}{\enode\makebox[0em]{$\,\,\scriptstyle{#1}$}} 
\mcm{\laenode}{1}{\enode\raisebox{1.2ex}[0ex][0ex]%
{\makebox[0em]{$\!\scriptstyle{#1}$}}}

\mcm{\Del}{0}{\Delta}
\mcm{\Deln}{1}{\Delta^{#1}}
\mcm{\Delop}{0}{\Delta^\op}
\mcm{\Delnop}{1}{(\Delta^{#1})^\op}

\mcm{\Trunc}{1}{#1\hyph\fcat{Trunc}}
\mcm{\inteqv}{0}{\,\sim\,}

\mcm{\ofdim}{1}{\,\of_{#1}\,}

\newcommand{\crossmark}{\textsf{\textit{x}}}

\newcounter{bean}
\newcommand{\bibentry}[2]{
	{\small
	\begin{list}{\refstepcounter{bean}\label{bib:#1}{[\thebean]}}{\setlength{\leftmargin}{\leftmargini}\setlength{\rightmargin}{\leftmargini}}
	\item #2
	\end{list}}\noindent}

\newcommand{\cyte}[1]{[\ref{bib:#1}]}
\newcommand{\slsh}{/\linebreak[0]}
\newcommand{\dblslsh}{//\linebreak[0]}
\newcommand{\dt}{.\linebreak[0]}

\newcommand{\jnl}[4]{\emph{#1}~\textbf{#2} (#3), #4} 

\newcommand{\contrib}[5]{in \emph{#1}, #2, #3, #4, pp.~#5}

\newcommand{\contb}[4]{in \emph{#1}, #2, #3, pp.~#4}

\newcommand{\contbau}[5]{in \emph{#1}, #2, #3, #4, pp.~#5}

\newcommand{\contribau}[6]{in \emph{#1}, #2, #3, #4, #5, pp.~#6}

\newcommand{\epr}[2]{e-print \url{#1}, #2}
 
\newcommand{\eprint}[3]{\epr{#1}{#2}, #3 pages}

\newcommand{\web}[1]{available via \url{#1}} 

\newcommand{\webprint}[3]{\web{#1}, #2, #3 pages} 

\newcommand{\concept}[1]{\subsection*{\textsl{#1}}}
\newcommand{\demph}[1]{\emph{#1}}
\newcommand{\ds}[1]{\textbf{#1}}  
\newcommand{\lp}{\ds{L$\mathbf{'}$}}  
\newcommand{\defnheading}[1]{\section*{Definition \ds{#1}}}	
\newcommand{\defnsheading}[1]{\section*{Definitions \ds{#1}}}	
\newcommand{\lowdimsheading}[1]{\section*{Definition \ds{#1} for $n\leq 2$}}
\newcommand{\lowsdimsheading}[1]{\section*{Definitions \ds{#1} for $n\leq 2$}}

\mcm{\parpair}{2}{\pile{\rTo^{\scriptstyle #1}\\ \rTo_{\scriptstyle #2}}}
\mcm{\parpairu}{0}{\pile{\rTo\\ \rTo}}
\newcommand{\url}[1]{{\tt #1}}
\mcm{\of}{0}{\raisebox{0.08ex}{\ensuremath{\scriptstyle\circ}}}
\mcm{\sof}{0}{\raisebox{0.08ex}{\ensuremath{\scriptscriptstyle\circ}}}

\newlength{\gwidth}	
\newlength{\gvert}	
\newlength{\gdrop}	
\newlength{\gbaredrop}	
\newlength{\goffset}	
\newlength{\gtemp}	

\newcommand{\present}[1]{%
\makebox[1\gwidth]{%
\rule[-1\gdrop]{0ex}{1\gvert}%
\raisebox{-1\gbaredrop}{#1}}}

\newcommand{\presentl}[1]{%
\makebox[1\gwidth][l]{%
\rule[-1\gdrop]{0ex}{1\gvert}%
\raisebox{-1\gbaredrop}{#1}}}

\newcommand{\presentr}[1]{%
\makebox[1\gwidth][r]{%
\rule[-1\gdrop]{0ex}{1\gvert}%
\raisebox{-1\gbaredrop}{#1}}}

\newcommand{\ginitdims}[2]{
\setlength{\unitlength}{1em}
\setlength{\goffset}{.25\unitlength}
\setlength{\gwidth}{#1\unitlength}
\setlength{\gvert}{#2\unitlength}
\setlength{\gdrop}{.5\gvert}
\addtolength{\gdrop}{-1\goffset}
\setlength{\gbaredrop}{1\gdrop}
\addtolength{\gvert}{.6\unitlength}
\addtolength{\gdrop}{.3\unitlength}}	

\newcommand{\cinitdims}[2]{
\setlength{\unitlength}{1em}
\setlength{\goffset}{.35\unitlength}
\setlength{\gwidth}{#1\unitlength}
\setlength{\gvert}{#2\unitlength}
\setlength{\gdrop}{.5\gvert}
\addtolength{\gdrop}{-1\goffset}
\setlength{\gbaredrop}{1\gdrop}
\addtolength{\gvert}{.6\unitlength}
\addtolength{\gdrop}{.3\unitlength}}	

\newcommand{\gsinitdims}[2]{
\setlength{\unitlength}{0.5em}
\setlength{\goffset}{.25\unitlength}
\setlength{\gwidth}{#1\unitlength}
\setlength{\gvert}{#2\unitlength}
\setlength{\gdrop}{.5\gvert}
\addtolength{\gdrop}{-1\goffset}
\setlength{\gbaredrop}{1\gdrop}
\addtolength{\gvert}{.6\unitlength}
\addtolength{\gdrop}{.3\unitlength}}	

\newcommand{\sidespic}[1]{%
\settowidth{\gtemp}{\ensuremath{#1}}%
\addtolength{\gwidth}{1\gtemp}}

\newcommand{\abovepic}[1]{%
\settoheight{\gtemp}{\ensuremath{#1}}%
\addtolength{\gvert}{1\gtemp}%
\settodepth{\gtemp}{\ensuremath{#1}}%
\addtolength{\gvert}{1\gtemp}}

\newcommand{\belowpic}[1]{%
\settoheight{\gtemp}{\ensuremath{#1}}%
\addtolength{\gvert}{1\gtemp}%
\addtolength{\gdrop}{1\gtemp}%
\settodepth{\gtemp}{\ensuremath{#1}}%
\addtolength{\gvert}{1\gtemp}%
\addtolength{\gdrop}{1\gtemp}}

\newcommand{\cell}[4]{\put(#1,#2){\makebox(0,0)[#3]{\ensuremath{#4}}}}
\mcm{\zmark}{0}{\scriptstyle{\bullet}}

\newcommand{\pregfst}[1]{%
\begin{picture}(0.5,0.2)(-0.5,-0.2)%
\cell{-0.1}{-0.2}{tr}{#1}%
\cell{0}{0}{c}{\zmark}%
\end{picture}}

\mcm{\gfst}{1}{%
\ginitdims{0.5}{0.4}%
\sidespic{#1}%
\belowpic{#1}%
\presentr{\pregfst{#1}}}

\newcommand{\preglst}[1]{%
\begin{picture}(0.5,0.2)(0,-0.2)%
\cell{0.1}{-0.2}{tl}{#1}%
\cell{0.05}{0}{c}{\zmark}%
\end{picture}}

\mcm{\glst}{1}{%
\ginitdims{.5}{.4}%
\sidespic{#1}%
\belowpic{#1}%
\presentl{\preglst{#1}}}

\newcommand{\preglft}[1]{%
\begin{picture}(0,0.2)(0,-0.2)%
\cell{-0.1}{-0.2}{tr}{#1}%
\cell{0.05}{0}{c}{\zmark}%
\end{picture}}

\mcm{\glft}{1}{%
\ginitdims{0}{.4}%
\belowpic{#1}%
\present{\preglft{#1}}}

\newcommand{\pregrgt}[1]{%
\begin{picture}(0,0.2)(0,-0.2)%
\cell{0.1}{-0.2}{tl}{#1}%
\cell{0.05}{0}{c}{\zmark}%
\end{picture}}

\mcm{\grgt}{1}{%
\ginitdims{0}{.4}%
\belowpic{#1}%
\present{\pregrgt{#1}}}

\newcommand{\pregblw}[1]{%
\begin{picture}(0,0.3)(0,-0.3)
\cell{0}{-0.3}{t}{#1}%
\cell{0.05}{0}{c}{\zmark}%
\end{picture}}

\mcm{\gblw}{1}{%
\ginitdims{0}{.6}%
\belowpic{#1}%
\present{\pregblw{#1}}}

\newcommand{\pregfbw}[1]{%
\begin{picture}(0,0.65)(0,-0.65)
\cell{0}{-0.65}{t}{#1}%
\cell{0.05}{0}{c}{\zmark}%
\end{picture}}

\mcm{\gfbw}{1}{%
\ginitdims{0}{1.3}%
\belowpic{#1}%
\present{\pregfbw{#1}}}

\newcommand{\pregzero}[1]{%
\begin{picture}(0.8,0.4)(-0.4,-0.4)
\cell{0}{-0.4}{t}{#1}%
\cell{0}{0}{c}{\zmark}%
\end{picture}}

\mcm{\gzero}{1}{%
\ginitdims{0.8}{.6}%
\belowpic{#1}%
\sidespic{#1}%
\present{\pregzero{#1}}}

\newcommand{\pregone}[1]{%
\begin{picture}(5,0.4)(0,-0.2)%
\cell{2.5}{0.2}{b}{#1}%
\put(0,0){\vector(1,0){4.9}}
\end{picture}}

\mcm{\gone}{1}{%
\ginitdims{5}{0.4}%
\abovepic{#1}%
\present{\pregone{#1}}}

\newcommand{\pregtwo}[3]{%
\begin{picture}(5,3.4)(0,-0.2)%
\cell{2.5}{3.2}{b}{#1}%
\cell{2.5}{-.2}{t}{#2}%
\cell{2.7}{1.5}{l}{#3}%
\qbezier(0,1.5)(2.5,4.5)(5,1.5)%
\qbezier(0,1.5)(2.5,-1.5)(5,1.5)%
\put(5,1.5){\vector(1,-1){0}}%
\put(5,1.5){\vector(1,1){0}}%
\put(2.5,2.5){\vector(0,-1){2}}%
\end{picture}}

\mcm{\gtwo}{3}{%
\ginitdims{5}{3.4}%
\abovepic{#1}%
\belowpic{#2}%
\present{\pregtwo{#1}{#2}{#3}}}

\newcommand{\pregthree}[5]{%
\begin{picture}(5,5.4)(0,-1.2)%
\cell{2.5}{4.2}{b}{#1}%
\cell{1.5}{1.7}{b}{#2}%
\cell{2.5}{-1.2}{t}{#3}%
\cell{2.7}{2.75}{l}{#4}%
\cell{2.7}{0.25}{l}{#5}%
\qbezier(0,1.5)(2.5,6.5)(5,1.5)%
\qbezier(0,1.5)(2.5,-3.5)(5,1.5)%
\put(0,1.5){\vector(1,0){5}}%
\put(2.5,3.5){\vector(0,-1){1.5}}%
\put(2.5,1){\vector(0,-1){1.5}}%
\put(5,1.5){\vector(1,-3){0}}%
\put(5,1.5){\vector(1,3){0}}%
\end{picture}}

\mcm{\gthree}{5}{%
\ginitdims{5}{5.4}%
\abovepic{#1}%
\belowpic{#3}%
\present{\pregthree{#1}{#2}{#3}{#4}{#5}}}

\newcommand{\pregfour}[7]{%
\begin{picture}(5,8.4)(0,-2.7)%
\cell{2.5}{5.7}{b}{#1}%
\cell{1.5}{2.8}{b}{#2}%
\cell{1.5}{0.2}{t}{#3}%
\cell{2.5}{-2.7}{t}{#4}%
\cell{2.7}{4.25}{l}{#5}%
\cell{2.7}{1.5}{l}{#6}%
\cell{2.7}{-1.25}{l}{#7}%
\qbezier(0,1.5)(2.5,9.5)(5,1.5)%
\qbezier(0,1.5)(2.5,4)(5,1.5)%
\qbezier(0,1.5)(2.5,-1)(5,1.5)%
\qbezier(0,1.5)(2.5,-6.5)(5,1.5)%
\put(2.5,5.25){\vector(0,-1){2}}%
\put(2.5,2.5){\vector(0,-1){2}}%
\put(2.5,-0.25){\vector(0,-1){2}}%
\put(5,1.5){\vector(1,-4){0}}%
\put(5,1.5){\vector(4,-3){0}}%
\put(5,1.5){\vector(4,3){0}}%
\put(5,1.5){\vector(1,4){0}}%
\end{picture}}

\mcm{\gfour}{7}{%
\ginitdims{5}{8.4}%
\abovepic{#1}%
\belowpic{#4}%
\present{\pregfour{#1}{#2}{#3}{#4}{#5}{#6}{#7}}}

\newcommand{\pregthreecell}[5]{%
\begin{picture}(8,5)(-4,-2.5)%
\cell{0}{2.5}{b}{#1}%
\cell{0}{-2.5}{t}{#2}%
\cell{-1.7}{0}{r}{#3}%
\cell{1.7}{0}{l}{#4}%
\cell{0}{0.2}{b}{#5}%
\qbezier(-4,0)(0,4.2)(4,0)%
\qbezier(-4,0)(0,-4.2)(4,0)%
\qbezier(-0.5,1.8)(-2.5,0)(-0.5,-1.8)%
\qbezier(0.5,1.8)(2.5,0)(0.5,-1.8)%
\put(-1,0){\vector(1,0){2}}%
\put(4,0){\vector(1,-1){0}}%
\put(4,0){\vector(1,1){0}}%
\put(-0.5,-1.8){\vector(1,-1){0}}%
\put(0.5,-1.8){\vector(-1,-1){0}}%
\end{picture}}

\mcm{\gthreecell}{5}{%
\ginitdims{8}{5}%
\abovepic{#1}%
\belowpic{#2}%
\present{\pregthreecell{#1}{#2}{#3}{#4}{#5}}}

\newcommand{\pregthreecellu}{%
\begin{picture}(5,3.4)(-0.5,-0.2)%
\qbezier(-.5,1.5)(2,4.5)(4.5,1.5)%
\qbezier(-.5,1.5)(2,-1.5)(4.5,1.5)%
\qbezier(1.5,2.7)(0.5,1.5)(1.5,0.3)%
\qbezier(2.5,2.7)(3.5,1.5)(2.5,0.3)%
\put(1.3,1.5){\vector(1,0){1.4}}%
\put(4.5,1.5){\vector(1,-1){0}}%
\put(4.5,1.5){\vector(1,1){0}}%
\put(1.5,0.3){\vector(2,-3){0}}%
\put(2.5,0.3){\vector(-2,-3){0}}%
\end{picture}}

\mcm{\gthreecellu}{0}{%
\ginitdims{5}{3.4}%
\present{\pregthreecellu}}

\newcommand{\pregtwocentre}[3]{%
\begin{picture}(5,3.4)(0,-0.2)%
\cell{2.5}{3.2}{b}{#1}%
\cell{2.5}{-.2}{t}{#2}%
\cell{2.5}{1.5}{c}{#3}%
\qbezier(0,1.5)(2.5,4.5)(5,1.5)%
\qbezier(0,1.5)(2.5,-1.5)(5,1.5)%
\put(5,1.5){\vector(1,-1){0}}%
\put(5,1.5){\vector(1,1){0}}%
\put(2.5,2.5){\vector(0,-1){2}}%
\end{picture}}

\mcm{\gtwocentre}{3}{%
\ginitdims{5}{3.4}%
\abovepic{#1}%
\belowpic{#2}%
\present{\pregtwocentre{#1}{#2}{#3}}}

\newcommand{\pregspecialone}[9]{%
\begin{picture}(8,8)(-4,-4)%
\cell{0}{3.9}{b}{#1}%
\cell{-2}{-0.2}{t}{#2}%
\cell{0}{-3.9}{t}{#3}%
\cell{-1.5}{1.1}{r}{#4}%
\cell{0.2}{1.5}{l}{#5}%
\cell{1.5}{1.1}{l}{#6}%
\cell{0.2}{-2}{l}{#7}%
\cell{-0.9}{2.3}{b}{#8}%
\cell{0.9}{2.3}{b}{#9}%
\qbezier(-4,0)(0,8)(4,0)%
\qbezier(-4,0)(0,-8)(4,0)%
\qbezier(-0.5,3.4)(-3.5,2)(-0.5,0.6)%
\qbezier(0.5,3.4)(3.5,2)(0.5,0.6)%
\put(-4,0){\vector(1,0){8}}%
\put(0,3.4){\vector(0,-1){2.8}}%
\put(0,-0.8){\vector(0,-1){2.4}}%
\put(-1.5,2.2){\vector(1,0){1.2}}%
\put(0.3,2.2){\vector(1,0){1.2}}%
\put(4,0){\vector(1,-2){0}}%
\put(4,0){\vector(1,2){0}}%
\put(-0.5,0.6){\vector(2,-1){0}}%
\put(0.5,0.6){\vector(-2,-1){0}}%
\end{picture}}

\mcm{\gspecialone}{9}{%
\ginitdims{8}{8}%
\abovepic{#1}%
\belowpic{#3}%
\present{\pregspecialone{#1}{#2}{#3}{#4}{#5}{#6}{#7}{#8}{#9}}}

\newcommand{\pregspecialtwo}{%
\begin{picture}(5,3.4)(0,-0.2)%
\qbezier(0,1.5)(2.5,4.5)(5,1.5)%
\qbezier(0,1.5)(2.5,-1.5)(5,1.5)%
\qbezier(1.7,2.5)(0,1.5)(1.7,0.5)%
\qbezier(3.3,2.5)(5,1.5)(3.3,0.5)%
\put(5,1.5){\vector(1,-1){0}}%
\put(5,1.5){\vector(1,1){0}}%
\put(1.7,0.5){\vector(3,-2){0}}%
\put(3.3,0.5){\vector(-3,-2){0}}%
\put(2.5,2.5){\vector(0,-1){2}}%
\put(1.2,1.5){\vector(1,0){1}}%
\put(2.8,1.5){\vector(1,0){1}}%
\end{picture}}

\mcm{\gspecialtwo}{0}{%
\ginitdims{5}{3.4}%
\present{\pregspecialtwo}}

\newcommand{\pregspecialthree}{%
\begin{picture}(5,5.4)(0,-1.2)%
\qbezier(0,1.5)(2.5,6.5)(5,1.5)%
\qbezier(0,1.5)(2.5,-3.5)(5,1.5)%
\qbezier(2,3.5)(1,2.75)(2,2)%
\qbezier(3,3.5)(4,2.75)(3,2)%
\qbezier(2,1)(1,0.25)(2,-0.5)%
\qbezier(3,1)(4,0.25)(3,-0.5)%
\put(0,1.5){\vector(1,0){5}}%
\put(1.5,2.75){\vector(1,0){2}}%
\put(1.5,0.25){\vector(1,0){2}}%
\put(5,1.5){\vector(1,-3){0}}%
\put(5,1.5){\vector(1,3){0}}%
\put(2,2){\vector(1,-1){0}}%
\put(3,2){\vector(-1,-1){0}}%
\put(2,-0.5){\vector(1,-1){0}}%
\put(3,-0.5){\vector(-1,-1){0}}%
\end{picture}}

\mcm{\gspecialthree}{0}{%
\ginitdims{5}{5.4}%
\present{\pregspecialthree}}

\newcommand{\pregonew}[1]{%
\begin{picture}(8,0.4)(0,-0.2)%
\cell{4}{0.2}{b}{#1}%
\put(0,0){\vector(1,0){8}}%
\end{picture}}

\mcm{\gonew}{1}{%
\ginitdims{8}{0.4}%
\abovepic{#1}%
\present{\pregonew{#1}}}

\mcm{\gzersu}{0}{%
\gsinitdims{0}{.6}%
\present{\pregblw{}}}

\mcm{\gonesu}{0}{%
\gsinitdims{5}{0.4}%
\present{\pregone{}}}

\mcm{\gtwosu}{0}{%
\gsinitdims{5}{3.4}%
\present{\pregtwo{}{}{}}}

\mcm{\gthreesu}{0}{%
\gsinitdims{5}{5.4}%
\present{\pregthree{}{}{}{}{}}}

\mcm{\gfoursu}{0}{%
\gsinitdims{5}{8.4}%
\present{\pregfour{}{}{}{}{}{}{}}}

\newcommand{\prectwo}[3]%
{\begin{picture}(4.2,3.4)(-0.1,-0.2)%
\cell{2}{3.2}{b}{#1}%
\cell{2}{-0.2}{t}{#2}%
\cell{2.2}{1.5}{l}{#3}%
\qbezier(0,2)(2,4)(4,2)%
\qbezier(0,1)(2,-1)(4,1)%
\put(4,2){\vector(1,-1){0}}%
\put(4,1){\vector(1,1){0}}%
\put(2,2.5){\vector(0,-1){2}}%
\end{picture}}

\mcm{\ctwo}{3}{%
\cinitdims{4.2}{3.4}%
\abovepic{#1}%
\belowpic{#2}%
\present{\prectwo{#1}{#2}{#3}}}

\newcommand{\precthree}[5]{%
\begin{picture}(4.2,5.4)(-0.1,-0.2)%
\cell{2}{5.2}{b}{#1}%
\cell{1}{2.7}{b}{#2}%
\cell{2}{-.2}{t}{#3}%
\cell{2.2}{3.75}{l}{#4}%
\cell{2.2}{1.25}{l}{#5}%
\qbezier(0,3)(2,7)(4,3)%
\qbezier(0,2)(2,-2)(4,2)%
\put(0,2.5){\vector(1,0){4}}%
\put(2,4.5){\vector(0,-1){1.5}}%
\put(2,2){\vector(0,-1){1.5}}%
\put(4,3){\vector(1,-3){0}}%
\put(4,2){\vector(1,3){0}}%
\end{picture}}

\mcm{\cthree}{5}{%
\cinitdims{4.2}{5.4}%
\abovepic{#1}%
\belowpic{#3}%
\present{\precthree{#1}{#2}{#3}{#4}{#5}}}

\newcommand{\prectwoop}[3]%
{\begin{picture}(4.2,3.4)(-0.1,-0.2)%
\cell{2}{3.2}{b}{#1}%
\cell{2}{-0.2}{t}{#2}%
\cell{2.2}{1.5}{l}{#3}%
\qbezier(0,2)(2,4)(4,2)%
\qbezier(0,1)(2,-1)(4,1)%
\put(0,2){\vector(-1,-1){0}}%
\put(0,1){\vector(-1,1){0}}%
\put(2,2.5){\vector(0,-1){2}}%
\end{picture}}

\mcm{\ctwoop}{3}{%
\cinitdims{4.2}{3.4}%
\abovepic{#1}%
\belowpic{#2}%
\present{\prectwoop{#1}{#2}{#3}}}

\newcommand{\prectwopar}[4]{%
\begin{picture}(4.2,3.4)(-0.1,-0.2)%
\cell{2}{3.2}{b}{#1}%
\cell{2}{-0.2}{t}{#2}%
\cell{1.6}{1.5}{r}{#3}%
\cell{2.4}{1.5}{l}{#4}%
\qbezier(0,2)(2,4)(4,2)%
\qbezier(0,1)(2,-1)(4,1)%
\put(4,2){\vector(1,-1){0}}%
\put(4,1){\vector(1,1){0}}%
\put(1.8,2.5){\vector(0,-1){2}}%
\put(2.2,2.5){\vector(0,-1){2}}%
\end{picture}}

\mcm{\ctwopar}{4}{%
\cinitdims{4.2}{3.4}%
\abovepic{#1}%
\belowpic{#2}%
\present{\prectwopar{#1}{#2}{#3}{#4}}}

\newcommand{\precthreein}[5]{%
\begin{picture}(4.2,5.4)(-0.1,-0.2)%
\cell{2}{5.2}{b}{#1}%
\cell{1}{2.7}{b}{#2}%
\cell{2}{-.2}{t}{#3}%
\cell{2.2}{3.75}{l}{#4}%
\cell{2.2}{1.25}{l}{#5}%
\qbezier(0,3)(2,7)(4,3)%
\qbezier(0,2)(2,-2)(4,2)%
\put(0,2.5){\vector(1,0){4}}%
\put(2,4.5){\vector(0,-1){1.5}}%
\put(2,0.5){\vector(0,1){1.5}}%
\put(4,3){\vector(1,-3){0}}%
\put(4,2){\vector(1,3){0}}%
\end{picture}}

\mcm{\cthreein}{5}{%
\cinitdims{4.2}{5.4}%
\abovepic{#1}%
\belowpic{#3}%
\present{\precthreein{#1}{#2}{#3}{#4}{#5}}}

\newcommand{\precthreecell}[5]{%
\begin{picture}(8.2,5)(-4.1,-2.5)%
\cell{0}{2.5}{b}{#1}%
\cell{0}{-2.5}{t}{#2}%
\cell{-1.7}{0}{r}{#3}%
\cell{1.7}{0}{l}{#4}%
\cell{0}{0.2}{b}{#5}%
\qbezier(-4,0.5)(0,4)(4,0.5)%
\qbezier(-4,-0.5)(0,-4)(4,-0.5)%
\qbezier(-0.5,2)(-2.5,0)(-0.5,-2)%
\qbezier(0.5,2)(2.5,0)(0.5,-2)%
\put(-1,0){\vector(1,0){2}}%
\put(4,0.5){\vector(1,-1){0}}%
\put(4,-0.5){\vector(1,1){0}}%
\put(-0.5,-2){\vector(1,-1){0}}%
\put(0.5,-2){\vector(-1,-1){0}}%
\end{picture}}

\mcm{\cthreecell}{5}{%
\cinitdims{8.2}{5}%
\abovepic{#1}%
\belowpic{#2}%
\present{\precthreecell{#1}{#2}{#3}{#4}{#5}}}

\newcommand{\precthreecellpar}[6]{%
\begin{picture}(8.2,5)(-4.1,-2.5)%
\cell{0}{2.5}{b}{#1}%
\cell{0}{-2.5}{t}{#2}%
\cell{-1.7}{0}{r}{#3}%
\cell{1.7}{0}{l}{#4}%
\cell{0}{0.4}{b}{#5}%
\cell{0}{-0.4}{t}{#6}%
\qbezier(-4,0.5)(0,4)(4,0.5)%
\qbezier(-4,-0.5)(0,-4)(4,-0.5)%
\qbezier(-0.5,2)(-2.5,0)(-0.5,-2)%
\qbezier(0.5,2)(2.5,0)(0.5,-2)%
\put(-1,0.2){\vector(1,0){2}}%
\put(-1,-0.2){\vector(1,0){2}}%
\put(4,0.5){\vector(1,-1){0}}%
\put(4,-0.5){\vector(1,1){0}}%
\put(-0.5,-2){\vector(1,-1){0}}%
\put(0.5,-2){\vector(-1,-1){0}}%
\end{picture}}

\mcm{\cthreecellpar}{6}{%
\cinitdims{8.2}{5}%
\abovepic{#1}%
\belowpic{#2}%
\present{\precthreecellpar{#1}{#2}{#3}{#4}{#5}{#6}}}

\newcommand{\prectwov}[5]{%
\begin{picture}(3.4,4.2)(0.8,0.9)%
\cell{2.5}{5.1}{b}{#1}%
\cell{2.5}{0.9}{t}{#2}%
\cell{0.8}{3}{r}{#3}%
\cell{4.2}{3}{l}{#4}%
\cell{2.5}{3.2}{b}{#5}%
\qbezier(2,5)(0,3)(2,1)%
\qbezier(3,5)(5,3)(3,1)%
\put(2,1){\vector(1,-1){0}}%
\put(3,1){\vector(-1,-1){0}}%
\put(1.5,3){\vector(1,0){2}}%
\end{picture}}

\mcm{\ctwov}{5}{%
\cinitdims{3.4}{4.2}%
\abovepic{#1}%
\belowpic{#2}%
\sidespic{#3}%
\sidespic{#4}%
\present{\prectwov{#1}{#2}{#3}{#4}{#5}}}

\newcommand{\precthreecellv}[7]{%
\begin{picture}(5,8.2)(0.5,-1.6)%
\cell{3}{6.6}{b}{#1}%
\cell{3}{-1.6}{t}{#2}%
\cell{0.5}{2.5}{r}{#3}%
\cell{5.5}{2.5}{l}{#4}%
\cell{3}{4.2}{b}{#5}%
\cell{3}{0.8}{t}{#6}%
\cell{3.2}{2.5}{l}{#7}%
\qbezier(3.5,6.5)(7,2.5)(3.5,-1.5)%
\qbezier(2.5,6.5)(-1,2.5)(2.5,-1.5)%
\put(2.5,-1.5){\vector(1,-1){0}}%
\put(3.5,-1.5){\vector(-1,-1){0}}%
\qbezier(1,3)(3,5)(5,3)%
\qbezier(1,2)(3,0)(5,2)%
\put(5,3){\vector(1,-1){0}}%
\put(5,2){\vector(1,1){0}}%
\put(3,3.5){\vector(0,-1){2}}%
\end{picture}}

\mcm{\cthreecellv}{7}{%
\cinitdims{5}{8.2}%
\abovepic{#1}%
\belowpic{#2}%
\sidespic{#3}%
\sidespic{#4}%
\present{\precthreecellv{#1}{#2}{#3}{#4}{#5}{#6}{#7}}}

\newcommand{\pretopez}[2]{%
\begin{picture}(2.6,2.3)(-1.3,-2.2)%
\cell{0}{-2.2}{t}{#1}%
\cell{0}{-1.2}{c}{#2}%
\qbezier(0,0)(-2,-2)(0,-2)%
\qbezier(0,0)(2,-2)(0,-2)%
\put(0,0){\vector(-1,1){0}}%
\end{picture}}

\mcm{\topez}{2}{%
\ginitdims{2.6}{2.3}%
\belowpic{#1}%
\present{\pretopez{#1}{#2}}}

\newcommand{\pretopea}[3]{%
\begin{picture}(4,1.9)(-2,-0,2)%
\cell{0}{1.7}{b}{#1}%
\cell{0}{-0.2}{t}{#2}%
\cell{0}{0.7}{c}{#3}%
\qbezier(-2,0)(0,3)(2,0)%
\put(-2,0){\vector(1,0){4}}%
\put(2,0){\vector(2,-3){0}}%
\end{picture}}

\mcm{\topea}{3}{%
\ginitdims{4}{1.9}%
\abovepic{#1}%
\belowpic{#2}%
\present{\pretopea{#1}{#2}{#3}}}

\newcommand{\pretopeb}[4]{%
\begin{picture}(4,2.2)(-2,-0.2)%
\cell{-1.1}{1}{br}{#1}%
\cell{1.1}{1}{bl}{#2}%
\cell{0}{-0.2}{t}{#3}%
\cell{0}{0.8}{c}{#4}%
\put(-2,0){\vector(1,1){2}}%
\put(0,2){\vector(1,-1){2}}%
\put(-2,0){\vector(1,0){4}}%
\end{picture}}

\mcm{\topeb}{4}{%
\ginitdims{4}{2.2}%
\belowpic{#3}%
\present{\pretopeb{#1}{#2}{#3}{#4}}}

\newcommand{\pretopec}[5]{%
\begin{picture}(4,2.2)(-2,-0.2)%
\cell{-1.8}{1}{br}{#1}%
\cell{0}{2.2}{b}{#2}%
\cell{1.8}{1}{bl}{#3}%
\cell{0}{-0.2}{t}{#4}%
\cell{0}{0.8}{c}{#5}%
\put(-2,0){\vector(1,2){1}}%
\put(-1,2){\vector(1,0){2}}%
\put(1,2){\vector(1,-2){1}}%
\put(-2,0){\vector(1,0){4}}%
\end{picture}}

\mcm{\topec}{5}{%
\ginitdims{4}{2.2}%
\sidespic{#1}%
\abovepic{#2}%
\sidespic{#3}%
\belowpic{#4}%
\present{\pretopec{#1}{#2}{#3}{#4}{#5}}}

\newcommand{\pretoped}[6]{%
\begin{picture}(4,2.5)(-2,-0.2)%
\cell{-2}{0.6}{br}{#1}%
\cell{-0.7}{2.2}{br}{#2}%
\cell{0.7}{2.2}{bl}{#3}%
\cell{2}{0.6}{bl}{#4}%
\cell{0}{-0.2}{t}{#5}%
\cell{0}{0.8}{c}{#6}%
\put(-2,0){\vector(1,3){0.5}}%
\put(-1.5,1.5){\vector(3,2){1.5}}%
\put(0,2.5){\vector(3,-2){1.5}}%
\put(1.5,1.5){\vector(1,-3){0.5}}%
\put(-2,0){\vector(1,0){4}}%
\end{picture}}

\mcm{\toped}{6}{%
\ginitdims{4}{2.5}%
\sidespic{#1}%
\abovepic{#2}%
\abovepic{#3}%
\sidespic{#4}%
\belowpic{#5}%
\present{\pretoped{#1}{#2}{#3}{#4}{#5}{#6}}}

\newcommand{\pretopeq}[5]{%
\begin{picture}(4,2.5)(-2,-0.2)%
\cell{-2}{0.6}{br}{#1}%
\cell{-1}{2.2}{br}{#2}%
\cell{2}{0.6}{bl}{#3}%
\cell{0}{-0.2}{t}{#4}%
\cell{0}{0.8}{c}{#5}%
\put(-2,0){\vector(1,3){0.5}}%
\put(-1.5,1.5){\vector(1,1){1}}%
\cell{0.9}{2.3}{c}{\ddots}
\put(1.5,1.5){\vector(1,-3){0.5}}%
\put(-2,0){\vector(1,0){4}}%
\end{picture}}

\mcm{\topeq}{5}{%
\ginitdims{4}{2.5}%
\sidespic{#1}%
\abovepic{#2}%
\sidespic{#3}%
\belowpic{#4}%
\present{\pretopeq{#1}{#2}{#3}{#4}{#5}}}

\newcommand{\pretopebase}[1]{%
\begin{picture}(4,0.4)(0,-0.2)%
\cell{2}{0.2}{b}{#1}%
\put(0,0){\vector(1,0){4}}%
\end{picture}}

\mcm{\topebase}{1}{%
\ginitdims{4}{0.4}%
\abovepic{#1}%
\present{\pretopebase{#1}}}

\newcommand{\pretopezs}[2]{%
\begin{picture}(2.6,2.3)(-1.3,-2.2)%
\cell{0}{-2.2}{t}{#1}%
\cell{0}{-1.2}{c}{#2}%
\qbezier(0,0)(-2,-2)(0,-2)%
\qbezier(0,0)(2,-2)(0,-2)%
\end{picture}}

\mcm{\topezs}{2}{%
\ginitdims{2.6}{2.3}%
\belowpic{#1}%
\present{\pretopezs{#1}{#2}}}

\newcommand{\pretopeas}[3]{%
\begin{picture}(4,1.9)(-2,-0,2)%
\cell{0}{1.7}{b}{#1}%
\cell{0}{-0.2}{t}{#2}%
\cell{0}{0.7}{c}{#3}%
\qbezier(-2,0)(0,3)(2,0)%
\put(-2,0){\line(1,0){4}}%
\end{picture}}

\mcm{\topeas}{3}{%
\ginitdims{4}{1.9}%
\abovepic{#1}%
\belowpic{#2}%
\present{\pretopeas{#1}{#2}{#3}}}

\newcommand{\pretopebs}[4]{%
\begin{picture}(4,2.2)(-2,-0.2)%
\cell{-1.1}{1}{br}{#1}%
\cell{1.1}{1}{bl}{#2}%
\cell{0}{-0.2}{t}{#3}%
\cell{0}{0.8}{c}{#4}%
\put(-2,0){\line(1,1){2}}%
\put(0,2){\line(1,-1){2}}%
\put(-2,0){\line(1,0){4}}%
\end{picture}}

\mcm{\topebs}{4}{%
\ginitdims{4}{2.2}%
\belowpic{#3}%
\present{\pretopebs{#1}{#2}{#3}{#4}}}

\newcommand{\pretopecs}[5]{%
\begin{picture}(4,2.2)(-2,-0.2)%
\cell{-1.8}{1}{br}{#1}%
\cell{0}{2.2}{b}{#2}%
\cell{1.8}{1}{bl}{#3}%
\cell{0}{-0.2}{t}{#4}%
\cell{0}{0.8}{c}{#5}%
\put(-2,0){\line(1,2){1}}%
\put(-1,2){\line(1,0){2}}%
\put(1,2){\line(1,-2){1}}%
\put(-2,0){\line(1,0){4}}%
\end{picture}}

\mcm{\topecs}{5}{%
\ginitdims{4}{2.2}%
\sidespic{#1}%
\abovepic{#2}%
\sidespic{#3}%
\belowpic{#4}%
\present{\pretopecs{#1}{#2}{#3}{#4}{#5}}}

\newcommand{\pretopeds}[6]{%
\begin{picture}(4,2.5)(-2,-0.2)%
\cell{-2}{0.6}{br}{#1}%
\cell{-0.7}{2.2}{br}{#2}%
\cell{0.7}{2.2}{bl}{#3}%
\cell{2}{0.6}{bl}{#4}%
\cell{0}{-0.2}{t}{#5}%
\cell{0}{0.8}{c}{#6}%
\put(-2,0){\line(1,3){0.5}}%
\put(-1.5,1.5){\line(3,2){1.5}}%
\put(0,2.5){\line(3,-2){1.5}}%
\put(1.5,1.5){\line(1,-3){0.5}}%
\put(-2,0){\line(1,0){4}}%
\end{picture}}

\mcm{\topeds}{6}{%
\ginitdims{4}{2.5}%
\sidespic{#1}%
\abovepic{#2}%
\abovepic{#3}%
\sidespic{#4}%
\belowpic{#5}%
\present{\pretopeds{#1}{#2}{#3}{#4}{#5}{#6}}}

\newcommand{\pretopeqs}[5]{%
\begin{picture}(4,2.5)(-2,-0.2)%
\cell{-2}{0.6}{br}{#1}%
\cell{-1}{2.2}{br}{#2}%
\cell{2}{0.6}{bl}{#3}%
\cell{0}{-0.2}{t}{#4}%
\cell{0}{0.8}{c}{#5}%
\put(-2,0){\line(1,3){0.5}}%
\put(-1.5,1.5){\line(1,1){1}}%
\cell{0.9}{2.3}{c}{\ddots}
\put(1.5,1.5){\line(1,-3){0.5}}%
\put(-2,0){\line(1,0){4}}%
\end{picture}}

\mcm{\topeqs}{5}{%
\ginitdims{4}{2.5}%
\sidespic{#1}%
\abovepic{#2}%
\sidespic{#3}%
\belowpic{#4}%
\present{\pretopeqs{#1}{#2}{#3}{#4}{#5}}}

\newcommand{\pretopebases}[1]{%
\begin{picture}(4,0.4)(0,-0.2)%
\cell{2}{0.2}{b}{#1}%
\put(0,0){\line(1,0){4}}%
\end{picture}}

\mcm{\topebases}{1}{%
\ginitdims{4}{0.4}%
\abovepic{#1}%
\present{\pretopebases{#1}}}

\newcommand{\pregdots}[6]{%
\begin{picture}(5,8.4)(0,-2.7)%
\cell{2.5}{5.7}{b}{#1}%
\cell{1.5}{2.8}{b}{#2}%
\cell{1.5}{0.2}{t}{#3}%
\cell{2.5}{-2.7}{t}{#4}%
\cell{2.7}{4.25}{l}{#5}%
\cell{2.7}{-1.25}{l}{#6}%
\qbezier(0,1.5)(2.5,9.5)(5,1.5)%
\qbezier(0,1.5)(2.5,4)(5,1.5)%
\qbezier(0,1.5)(2.5,-1)(5,1.5)%
\qbezier(0,1.5)(2.5,-6.5)(5,1.5)%
\put(2.5,5.25){\vector(0,-1){2}}%
\put(2.5,-0.25){\vector(0,-1){2}}%
\cell{2.5}{1.7}{c}{\vdots}%
\put(5,1.5){\vector(1,-4){0}}%
\put(5,1.5){\vector(4,-3){0}}%
\put(5,1.5){\vector(4,3){0}}%
\put(5,1.5){\vector(1,4){0}}%
\end{picture}}

\mcm{\gdots}{6}{%
\ginitdims{5}{8.4}%
\abovepic{#1}%
\belowpic{#4}%
\present{\pregdots{#1}{#2}{#3}{#4}{#5}{#6}}}

\newcommand{\presplitcoeqrhs}[2]%
{\begin{picture}(4.4,1.8)(-0.1,-0.2)%
\cell{2}{1.7}{b}{\scriptstyle{#1}}%
\cell{2}{-0.2}{t}{\scriptstyle{#2}}%
\qbezier(0,1)(2,-1)(4,1)%
\put(0,1.5){\vector(1,0){4}}%
\put(0,1){\vector(-1,1){0}}%
\end{picture}}

\mcm{\splitcoeqrhs}{2}{%
\cinitdims{4.4}{3.4}%
\abovepic{#1}%
\belowpic{#2}%
\present{\presplitcoeqrhs{#1}{#2}}}

\newcommand{\presplitcoeqlhs}[3]%
{\begin{picture}(4.4,1.8)(-0.1,-0.2)%
\cell{2}{1.8}{b}{\scriptstyle{#1}}%
\cell{2}{1.2}{t}{\scriptstyle{#2}}%
\cell{2}{-0.2}{t}{\scriptstyle{#3}}%
\qbezier(0,1)(2,-1)(4,1)%
\put(0,1.7){\vector(1,0){4}}%
\put(0,1.3){\vector(1,0){4}}%
\put(0,1){\vector(-1,1){0}}%
\end{picture}}

\mcm{\splitcoeqlhs}{3}{%
\cinitdims{4.4}{3.4}%
\abovepic{#1}%
\belowpic{#3}%
\present{\presplitcoeqlhs{#1}{#2}{#3}}}

\newcommand{\prectwocentre}[3]%
{\begin{picture}(4.2,3.4)(-0.1,-0.2)%
\cell{2}{3.2}{b}{#1}%
\cell{2}{-0.2}{t}{#2}%
\cell{2}{1.5}{c}{#3}%
\qbezier(0,2)(2,4)(4,2)%
\qbezier(0,1)(2,-1)(4,1)%
\put(4,2){\vector(1,-1){0}}%
\put(4,1){\vector(1,1){0}}%
\put(2,2.5){\vector(0,-1){2}}%
\end{picture}}

\mcm{\ctwocentre}{3}{%
\cinitdims{4.2}{3.4}%
\abovepic{#1}%
\belowpic{#2}%
\present{\prectwocentre{#1}{#2}{#3}}}

\newcommand{\prectwodotty}[3]%
{\begin{picture}(4.2,3.4)(-0.1,-0.2)%
\cell{2}{3.2}{b}{#1}%
\cell{2}{-0.2}{t}{#2}%
\cell{2.2}{1.5}{l}{#3}%
\qbezier(0,2)(2,4)(4,2)%
\qbezier(0,1)(2,-1)(4,1)%
\put(4,2){\vector(1,-1){0}}%
\put(4,1){\vector(1,1){0}}%
\multiput(2,2.5)(0,-0.25){7}{\makebox(0,0)[c]{$\cdot$}}
\put(2,0.5){\vector(0,-1){0}}%
\end{picture}}

\mcm{\ctwodotty}{3}{%
\cinitdims{4.2}{3.4}%
\abovepic{#1}%
\belowpic{#2}%
\present{\prectwodotty{#1}{#2}{#3}}}

\newcommand{\pregtwodotty}[3]{%
\begin{picture}(5,3.4)(0,-0.2)%
\cell{2.5}{3.2}{b}{#1}%
\cell{2.5}{-.2}{t}{#2}%
\cell{2.7}{1.5}{l}{#3}%
\qbezier(0,1.5)(2.5,4.5)(5,1.5)%
\qbezier(0,1.5)(2.5,-1.5)(5,1.5)%
\put(5,1.5){\vector(1,-1){0}}%
\put(5,1.5){\vector(1,1){0}}%
\multiput(2.5,2.5)(0,-0.25){7}{\makebox(0,0)[c]{$\cdot$}}
\put(2.5,0.5){\vector(0,-1){0}}%
\end{picture}}

\mcm{\gtwodotty}{3}{%
\ginitdims{5}{3.4}%
\abovepic{#1}%
\belowpic{#2}%
\present{\pregtwodotty{#1}{#2}{#3}}}

\newcommand{\prebundleint}[3]%
{\begin{picture}(4.4,1.8)(-0.1,-0.2)%
\cell{2}{1.6}{b}{\scriptstyle{#1}}%
\cell{2}{0.4}{b}{\scriptstyle{#2}}
\cell{2}{-0.35}{t}{\scriptstyle{#3}}%
\put(0,1.5){\vector(1,0){4}}%
\qbezier(0,1.25)(2,-0.75)(4,1.25)%
\qbezier(0,0.85)(2,-1.15)(4,0.85)%
\put(0,1.25){\vector(-1,1){0}}%
\put(0,0.85){\vector(-1,1){0}}%
\end{picture}}

\mcm{\bundleint}{3}{%
\cinitdims{4.4}{3.4}%
\abovepic{#1}%
\belowpic{#3}%
\present{\prebundleint{#1}{#2}{#3}}}

\newcommand{\precone}[1]{%
\begin{picture}(4.2,0.4)(-0.3,-0.2)%
\cell{1.8}{0.2}{b}{#1}%
\put(0,0){\vector(1,0){3.6}}%
\end{picture}}

\mcm{\cone}{1}{%
\cinitdims{4.2}{0.4}%
\abovepic{#1}%
\present{\precone{#1}}}

\newcommand{\preghappy}[1]{%
\begin{picture}(5,3.4)(0,-0.2)%
\cell{2.5}{-.2}{t}{#1}%
\qbezier(0,1.5)(2.5,-1.5)(5,1.5)%
\put(5,1.5){\vector(1,1){0}}%
\end{picture}}

\mcm{\ghappysu}{0}{%
\gsinitdims{5}{3.4}%
\present{\preghappy{}}}

\newcommand{\pregunhappy}[1]{%
\begin{picture}(5,3.4)(0,-0.2)%
\cell{2.5}{3.2}{b}{#1}%
\qbezier(0,1.5)(2.5,4.5)(5,1.5)%
\put(5,1.5){\vector(1,-1){0}}%
\end{picture}}

\mcm{\gunhappysu}{0}{%
\gsinitdims{5}{3.4}%
\present{\pregunhappy{}{}{}}}

\newcommand{\pregtwowide}[3]{%
\begin{picture}(10,3.4)(0,-0.2)%
\cell{5}{3.2}{b}{#1}%
\cell{5}{-.2}{t}{#2}%
\cell{5.2}{1.5}{l}{#3}%
\qbezier(0,1.5)(5,4.5)(10,1.5)%
\qbezier(0,1.5)(5,-1.5)(10,1.5)%
\put(10,1.5){\vector(2,-1){0}}%
\put(10,1.5){\vector(2,1){0}}%
\put(5,2.5){\vector(0,-1){2}}%
\end{picture}}

\mcm{\gtwowidesu}{0}{%
\gsinitdims{10}{3.4}%
\present{\pregtwowide{}{}{}}}

\newcommand{\preghole}{%
\begin{picture}(5,0)(0,0)%
\end{picture}}

\mcm{\ghole}{0}{%
\gsinitdims{5}{0}%
\present{\preghole}}

\mcm{\gzeros}{1}{%
\gsinitdims{0.8}{.6}%
\belowpic{\scriptstyle #1}%
\sidespic{\scriptstyle #1}%
\present{\pregzero{\scriptstyle #1}}}

\mcm{\gfsts}{1}{%
\gsinitdims{0.5}{0.4}%
\sidespic{\scriptstyle #1}%
\belowpic{\scriptstyle #1}%
\presentr{\pregfst{\scriptstyle #1}}}

\mcm{\glsts}{1}{%
\gsinitdims{.5}{.4}%
\sidespic{\scriptstyle #1}%
\belowpic{\scriptstyle #1}%
\presentl{\preglst{\scriptstyle #1}}}

\mcm{\gblws}{1}{%
\gsinitdims{0}{.6}%
\belowpic{#1}%
\present{\pregblw{\scriptstyle #1}}}

\mcm{\gfbws}{1}{%
\gsinitdims{0}{1.3}%
\belowpic{\scriptstyle #1}%
\present{\pregfbw{\scriptstyle #1}}}

\mcm{\glfts}{1}{%
\gsinitdims{0}{.4}%
\belowpic{\scriptstyle #1}%
\present{\preglft{\scriptstyle #1}}}

\mcm{\grgts}{1}{%
\gsinitdims{0}{.4}%
\belowpic{\scriptstyle #1}%
\present{\pregrgt{\scriptstyle #1}}}

\mcm{\gones}{1}{%
\gsinitdims{5}{0.4}%
\abovepic{\scriptstyle #1}%
\present{\pregone{\scriptstyle #1}}}

\mcm{\gtwos}{3}{%
\gsinitdims{5}{3.4}%
\abovepic{\scriptstyle #1}%
\belowpic{\scriptstyle #2}%
\present{\pregtwo{\scriptstyle #1}{\scriptstyle #2}{\scriptstyle #3}}}

\mcm{\gthrees}{5}{%
\gsinitdims{5}{5.4}%
\abovepic{\scriptstyle #1}%
\belowpic{\scriptstyle #3}%
\present{\pregthree{\scriptstyle #1}{\scriptstyle #2}{\scriptstyle #3}{\scriptstyle #4}{\scriptstyle #5}}}

\mcm{\gfours}{7}{%
\gsinitdims{5}{8.4}%
\abovepic{\scriptstyle #1}%
\belowpic{\scriptstyle #4}%
\present{\pregfour{\scriptstyle #1}{\scriptstyle #2}{\scriptstyle
#3}{\scriptstyle #4}{\scriptstyle #5}{\scriptstyle #6}{\scriptstyle #7}}}

\mcm{\gfstsu}{0}{%
\gsinitdims{0.5}{0.4}%
\presentr{\pregfst{}}}

\mcm{\glstsu}{0}{%
\gsinitdims{0.5}{0.4}%
\presentl{\preglst{}}}

\newcommand{\pretopeavar}[5]{%
\begin{picture}(4.4,1.9)(-2.2,-0.2)%
\cell{-2.2}{0}{br}{#1}
\cell{2.2}{0}{bl}{#2}
\cell{0}{1.7}{b}{#3}%
\cell{0}{-0.2}{t}{#4}%
\cell{0.2}{0.7}{l}{#5}%
\qbezier(-2,0)(0,3)(2,0)%
\put(-2,0){\vector(1,0){4}}%
\cell{-2}{0}{c}{\zmark}%
\cell{2}{0}{c}{\zmark}%
\put(2,0){\vector(2,-3){0}}%
\put(0,1.2){\vector(0,-1){1}}
\end{picture}}

\mcm{\topeavar}{5}{%
\ginitdims{4}{1.9}%
\sidespic{#1}%
\sidespic{#2}%
\abovepic{#3}%
\belowpic{#4}%
\sidespic{\zmark}%
\present{\pretopeavar{#1}{#2}{#3}{#4}{#5}}}

\renewcommand{\theenumi}{\roman{enumi}}
\begin{document}

\title{A Survey of Definitions of $n$-Category}

\author{Tom Leinster\\ \\
        \normalsize{Department of Pure Mathematics, University of
        Cambridge}\\ 
        \normalsize{Email: leinster@dpmms.cam.ac.uk}\\
        \normalsize{Web: http://www.dpmms.cam.ac.uk/$\sim$leinster}}

\date{}

\begin{titlepage} 

\maketitle
\thispagestyle{empty}

\vfill

\begin{center}	\bfseries
Abstract
\end{center}
\vspace*{-2ex}
\begin{quotation}
Many people have proposed definitions of `weak $n$-category'.  Ten of them
are presented here.  Each definition is given in two pages, with a further
two pages on what happens when $n\leq 2$.  The definitions can be read
independently.  Chatty bibliography follows.
\end{quotation}

\vfill

\begin{center}
\bfseries Contents
\end{center} 

\begin{center}
\begin{tabular}{lr}
Introduction		&\pageref{p:intro}	\\
Background		&\pageref{p:background}	\\
			&			\\
Definition \ds{Tr}	&\pageref{p:tr}		\\
Definition \ds{P}	&\pageref{p:p}		\\
Definitions \ds{B}	&\pageref{p:b}		\\
Definitions \ds{L}	&\pageref{p:l}		\\
Definition \ds{\lp}	&\pageref{p:lprime}	\\
Definition \ds{Si}	&\pageref{p:si}		\\
Definition \ds{Ta}	&\pageref{p:ta}		\\
Definition \ds{J}	&\pageref{p:j}		\\
Definition \ds{St}	&\pageref{p:st}		\\
Definition \ds{X}	&\pageref{p:x}		\\
			&			\\
Further Reading		&\pageref{p:biblio}
\end{tabular}
\end{center}

\end{titlepage}

\section*{Introduction}		\label{p:intro}

\begin{quote}
\textit{%
L\'evy \ldots\ once remarked to me that reading other mathematicians' research
gave him actual physical pain.}
\end{quote}
---J. L. Doob on the probabilist Paul L\'evy, \emph{Statistical Science}
   \textbf{1}, no.~1, 1986. 

\begin{quote}
\textit{%
Hell is other people.}
\end{quote}
---Jean-Paul Sartre, \emph{Huis Clos}.

\paragraph*{}

The last five years have seen a vast increase in the literature on
higher-dimensional categories.  Yet one question of central concern remains
resolutely unanswered: what exactly is a weak $n$-category?  There have,
notoriously, been many proposed definitions, but there seems to be a general
perception that most of these definitions are obscure, difficult and long.  
I hope that the present work will persuade the reader that this is not the
case, or at least does not \emph{need} to be: that while no existing
approach is without its mysteries, it is quite possible to state the
definitions in a concise and straightforward way.

\subsection*{What's in here, and what's not}

The sole purpose of this paper is to state several possible definitions of
weak $n$-category.  In particular, I have made no attempt to compare the
proposed definitions with one another (although certainly I hope that this
work will help with the task of comparison).  So the definitions of weak
$n$-category that follow may or may not be `equivalent'; I make no comment.
Moreover, I have not included any notions of weak functor or equivalence
between weak $n$-categories, which would almost certainly be required before
one could make any statement such as `Professor Yin's definition of weak
$n$-category is equivalent to Professor Yang's'.

I have also omitted any kind of motivational or introductory material.  The
`Further Reading' section lists various texts which attempt to explain the
relevance of $n$-categories and other higher categorical structures to
mathematics at large (and to physics and computer science).  I will just
mention two points here for those new to the area.  Firstly, it is easy to
define \emph{strict} $n$-categories (see `Preliminaries'), and it is true
that every weak $2$-category is equivalent to a strict $2$-category, but the
analogous statement fails for $n$-categories when $n>2$: so the difference
between the theories of weak and strict $n$-categories is nontrivial.
Secondly, the issue of comparing definitions of weak $n$-category is a
slippery one, as it is hard to say what it even \emph{means} for two such
definitions to be equivalent.  For instance, suppose you and I each have in
mind a definition of algebraic variety and of morphism of varieties; then we
might reasonably say that our definitions of variety are `equivalent' if your
category of varieties is equivalent to mine.  This makes sense because the
structure formed by varieties and their morphisms is a category.  It is
widely held that the structure formed by weak $n$-categories and the
functors, transformations, \ldots\ between them should be a weak
$(n+1)$-category; and if this is the case then the question is whether your
weak $(n+1)$-category of weak $n$-categories is equivalent to mine---but
whose definition of weak $(n+1)$-category are we using here\ldots?

This paper gives primary importance to $n$-categories, with other higher
categorical structures only mentioned where they have to be.  In writing it
this way I do not mean to imply that $n$-categories are the only interesting
structures in higher-dimensional category theory: on the contrary, I see the
subject as including a whole range of interesting structures, such as operads
and multicategories in their various forms, double and $n$-tuple categories,
computads and string diagrams, homotopy-algebras, $n$-vector spaces, and
structures appropriate for the study of braids, knots, graphs, cobordisms,
proof nets, flowcharts, circuit diagrams, \ldots.  Moreover, consideration of
$n$-categories seems inevitably to lead into consideration of some of these
other structures, as is borne out by the definitions below.  However,
$n$-categories are here allowed to upstage the other structures in what is
probably an unfair way.

Finally, I do not claim to have included \emph{all} the definitions of weak
$n$-category that have been proposed by people; in fact, I am aware that I
have omitted a few.  They are omitted purely because I am not familiar with
them.  More information can be found under `Further Reading'.

\subsection*{Layout}

The first section is `Background'.  This is mainly for reference, and it is
not really recommended that anyone starts reading here.  It begins with a
page on ordinary category theory, recalling those concepts that will be used
in the main text and fixing some terminology.  Everything here is completely
standard, and almost all of it can be found in any introductory book or
course on the subject; but only a small portion of it is used in each
definition of weak $n$-category.  There is then a page each on strict
$n$-categories and bicategories, again recalling widely-known material.

Next come the ten definitions of weak $n$-category.  They are absolutely
independent and self-contained, and can be read in any order.  No
significance should be attached to the order in which they are presented; I
tried to arrange them so that definitions with common themes were grouped
together in the sequence, but that is all.  (Some structures just don't fit
naturally into a single dimension.)

Each definition of weak $n$-category is given in two pages, so that if this
is printed double-sided then the whole definition will be visible on a
double-page spread.  This is followed, again in two pages, by an explanation
of the cases $n=0,1,2$.  We expect weak $0$-categories to be sets, weak
$1$-categories to be categories, and weak $2$-categories to be
bicategories---or at least, to resemble them to some reasonable degree---and
this is indeed the case for all of the definitions as long as we interpret
the word `reasonable' generously.  Each main definition is given in a formal,
minimal style, but the analysis of $n\leq 2$ is less formal and more
explanatory; partly the analysis of $n\leq 2$ is to show that the proposed
definition of $n$-category is a reasonable one, but partly it is for
illustrative purposes.  The reader who gets stuck on a definition might
therefore be helped by looking at $n\leq 2$.

Taking a definition of weak $n$-category and performing a rigorous comparison
between the case $n=2$ and bicategories is typically a long and tedious
process.  For this reason, I have not checked all the details in the $n\leq
2$ sections.  The extent to which I feel confident in my assertions can be
judged from the number of occurrences of phrases such as `probably' and `it
appears that', and by the presence or absence of references under `Further
Reading'.

There are a few exceptions to this overall scheme.  The section labelled
\ds{B} consists, in fact, of \emph{two} definitions of weak $n$-category, but
they are so similar in their presentation that it seemed wasteful to give
them two different sections.  (The reason for the name \ds{B} is explained
below.)  The same goes for definition \ds{L}, so we have definitions of weak
$n$-category called \ds{B1}, \ds{B2}, \ds{L1} and \ds{L2}.  A variant for
definition \ds{St} is also given (in the $n\leq 2$ section), but this goes
nameless.  However, definition \ds{X} is not strictly speaking a mathematical
definition at all: I was unable to find a way to present it in two pages, so
instead I have given an informal version, with one sub-definition (opetopic
set) done by example only.  The cases $n\leq 2$ are clear enough to be
analysed precisely.

Another complicating factor comes from those definitions which include a
notion of weak $\omega$-category ($=$ weak $\infty$-category).  There, the
pattern is very often to define weak $\omega$-category and then to define a
weak $n$-category as a weak $\omega$-category with only trivial cells in
dimensions $>n$.  This presents a problem when one comes to attempt a precise
analysis of $n\leq 2$, as even to determine what a weak $0$-category is
involves considering an infinite-dimensional structure.  For this reason it
is more convenient to redefine weak $n$-category in a way which never
mentions cells of dimension $>n$, by imitating the original definition of
weak $\omega$-category.  Of course, one then has to show that the two
different notions of weak $n$-category are equivalent, and again I have not
always done this with full rigour (and there is certainly not the space to
give proofs here).  So, this paper actually contains significantly more than
ten possible definitions of weak $n$-category.

`Further Reading' is the final section.  To keep the definitions of
$n$-category brief and self-contained, there are no citations at all in the
main text; so this section is a combination of reference list, historical
notes, and general comments, together with a few pointers to literature in
related areas.

\subsection*{Overview of the definitions}

Table~\ref{table:defns} shows some of the main features of the definitions of
weak $n$-category.  
\begin{table}
\centering
\begin{tabular}{lllll}
\emph{Definition}	&
\emph{Author(s)}	&
\emph{Shapes used}	&
\emph{A/the}		&
\emph{$\omega$?}	\\
\\
\ds{Tr}			&
Trimble			&
path parametrizations	&
the			&
\crossmark		\\
\ds{P}			&
Penon			&
globular		&
the			&
\checkmark		\\
\ds{B}			&
Batanin			&
globular		&
the (\ds{B1}), a (\ds{B2})&
\checkmark		\\
\ds{L}			&
Leinster		&
globular		&
the (\ds{L1}), a (\ds{L2})&
\checkmark		\\
\ds{\lp}		&
Leinster		&
globular		&
a			&
\checkmark		\\
\ds{Si}			&
Simpson			&
simplicial/globular	&
a			&
\crossmark		\\
\ds{Ta}			&
Tamsamani		&
simplicial/globular	&
a			&
\crossmark		\\
\ds{J}			&
Joyal			&
globular/simplicial	&
a			&
\checkmark		\\
\ds{St}			&
Street			&
simplicial		&
a			&
\checkmark		\\
\ds{X}			&
see text		&
opetopic		&
a			&
\crossmark		\\
\end{tabular}
\caption{Some features of the definitions}
\label{table:defns}
\end{table}
Each definition is given a name such as \ds{A} or \ds{Z}, according to the
name of the author from whom the definition is derived.  (Definition \ds{X}
is a combination of the work of many people, principally Baez, Dolan,
Hermida, Makkai and Power.)  The point of these abbreviations is to put some
distance between the definitions as proposed by those authors and the
definitions as stated below.  At the most basic level, I have in all cases
changed some notation and terminology.  Moreover, taking what is often a long
paper and turning it into a two-page definition has seldom been just a matter
of leaving out words; sometimes it has required a serious reshaping of the
concepts involved.  Whether the end result (the definition of weak
$n$-category) is mathematically the same as that of the original author is
not something I always know: on various occasions there have been passages in
the source paper that have been opaque to me, so I have guessed at the
author's intended meaning.  Finally, in several cases only a definition of
weak $\omega$-category was explicitly given, leaving me to supply the
definition of weak $n$-category for finite $n$.  In summary, then, I do
believe that I have given ten reasonable definitions of weak $n$-category,
but I do not guarantee that they are the same as those of the authors listed
in Table~\ref{table:defns}; ultimately, the responsibility for them is mine.

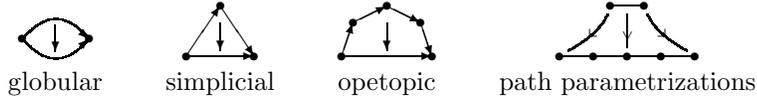
\begin{figure}
\centering
\begin{tabular}{ccccccc}%
\gfstsu\gtwosu\glstsu%
&&%
\raisebox{-3.4mm}{%
\setlength{\unitlength}{0.75mm}%
\begin{picture}(12,11)(-6,-2)
\cell{-6}{0}{c}{\zmark}
\cell{0}{9}{c}{\zmark}
\cell{6}{0}{c}{\zmark}
\put(-6,0){\vector(2,3){6}}
\put(0,9){\vector(2,-3){6}}
\put(-6,0){\vector(1,0){12}}
\put(0,6){\vector(0,-1){5}}
\end{picture}}%
&&%
\raisebox{-3.4mm}{%
\setlength{\unitlength}{0.75mm}%
\begin{picture}(16,11)(-8,-2)
\cell{-8}{0}{c}{\zmark}
\cell{-6}{6}{c}{\zmark}
\cell{0}{9}{c}{\zmark}
\cell{6}{6}{c}{\zmark}
\cell{8}{0}{c}{\zmark}
\put(-8,0){\vector(1,3){2}}
\put(-6,6){\vector(2,1){6}}
\put(0,9){\vector(2,-1){6}}
\put(6,6){\vector(1,-3){2}}
\put(-8,0){\vector(1,0){16}}
\put(0,6){\vector(0,-1){5}}
\end{picture}}%
&&%
\raisebox{-3.4mm}{%
\setlength{\unitlength}{0.75mm}%
\begin{picture}(24,11)(-12,-2)
\cell{-3}{9}{c}{\zmark}
\cell{3}{9}{c}{\zmark}
\cell{-12}{0}{c}{\zmark}
\cell{-6}{0}{c}{\zmark}
\cell{0}{0}{c}{\zmark}
\cell{6}{0}{c}{\zmark}
\cell{12}{0}{c}{\zmark}
\put(-3,9){\line(1,0){6}}
\put(-12,0){\line(1,0){24}}
\qbezier(-3.5,7.5)(-5,4)(-10.5,1)
\qbezier(3.5,7.5)(5,4)(10.5,1)
\put(0,7.5){\line(0,-1){6}}
\cell{-7}{3.5}{bl}{\scriptstyle\llcorner}
\cell{7.1}{3.5}{br}{\scriptstyle\lrcorner}
\cell{0.1}{3}{b}{\scriptscriptstyle\vee}
\end{picture}}%
\\
globular	&&
simplicial	&&
opetopic	&&
path parametrizations\\
\end{tabular}
\caption{Shapes used in the definitions}
\label{fig:shapes}
\end{figure}

The column headed `shapes used' refers to the different shapes of $m$-cell
(or `$m$-arrow', or `$m$-morphism') employed in the definitions.  These are
shown in Figure~\ref{fig:shapes}.

It has widely been observed that the various definitions of $n$-category fall
into two groups, according to the attitude one takes to the status of
composition.  This distinction can be explained by analogy with products.
Given two sets $A$ and $B$, one can define \emph{a product} of $A$ and $B$ to
be a triple $(P,p_1,p_2)$ where $P$ is a set and $p_1: P \go A$, $p_2: P \go
B$ are functions with the usual universal property.  This is of course the
standard thing to do in category theory, and in this context one can strictly
speaking never refer to \emph{the} product of $A$ and $B$.  On the other
hand, one could define \emph{the product} of $A$ and $B$ to be the set
$A\times B$ of ordered pairs $(a,b) = \{ \{a\}, \{a,b\} \}$ with $a\in A$ and
$b\in B$; this has the virtue of being definite and allowing one to speak of
\emph{the} product in the customary way, but involves a wholly artificial
construction.  Similarly, in some of the proposed definitions of weak
$n$-category, one can never speak of \emph{the} composite of morphisms $g$
and $f$, only of \emph{a} composite (of which there may be many, all equally
valid); but in some of the definitions one does have definite composites
$g\of f$, \emph{the} composite of $g$ and $f$.  (The use of the word `the' is
not meant to imply strictness, e.g.\ the three-fold composite $h\of (g\of f)$
will in general be different from the three-fold composite $(h\of g) \of f$.)
So this is the meaning of the column headed `a/the'; it might also have been
headed `indefinite/definite', `relational/functional', `universal/coherent',
or even `geometric/algebraic'.

All of the sections include a definition of weak $n$-category for
natural numbers $n$, but some also include a definition of weak
$\omega$-category (in which there are $m$-cells for all natural $m$).  This
is shown in the last column.

Finally, I warn the reader that the words `contractible' and `contraction'
occur in many of the definitions, but mean different things from
definition to definition.  This is simply to save having to invent new words
for concepts which are similar but not identical, and to draw attention to
the common idea.

\subsection*{Acknowledgements}

I first want to thank Eugenia Cheng and Martin Hyland.  Their involvement in
this project has both made it much more pleasurable for me and provided a
powerful motivating force.  Without them, I suspect it would still not be
done.

I prepared for writing this by giving a series of seminars (one definition
per week) in Cambridge in spring 2001, and am grateful to the participants:
the two just mentioned, Mario C\'accamo, Marcelo Fiore, and Joe Templeton.  I
would also like to thank those who have contributed over the years to the
many Cambridge Category Theory seminars on the subject of $n$-categories,
especially Jeff Egger (who introduced me to Tamsamani's definition), Peter
Johnstone, and Craig Snydal (with whom I have also had countless interesting
conversations on the subject).

Todd Trimble was generous enough to let me publish his definition for the
first time, and to cast his eye over a draft of what appears below as
definition \ds{Tr}---though all errors, naturally, are mine.

I am also grateful to the other people with whom I have had helpful
communications, including Michael Batanin (who told me about Penon's
definition), David Carlton, Jack Duskin, Anders Kock, Peter May, Carlos
Simpson, Ross Street, Bertrand Toen, Dominic Verity, Marek Zawadowski (who
told me about Joyal's definition), and surely others, whose names I apologize
for omitting.

Many of the diagrams were drawn using Paul Taylor's commutative diagrams
package.

It is a pleasure to thank St John's College, Cambridge, where I hold the
Laurence Goddard Fellowship, for their support. 

\clearpage

\section*{Background}		\label{p:background}

\concept{Category Theory}

Here is a summary of the categorical background and terminology needed in order
to read the entire paper.  The reader who isn't familiar with everything below
shouldn't be put off: each individual Definition only uses some of it.

I assume familiarity with \demph{categories}, \demph{functors},
\demph{natural transformations}, \demph{adjunctions}, \demph{limits}, and
\demph{monads} and their \demph{algebras}.  Limits include \demph{products},
\demph{pullbacks} (with the pullback of a diagram $X \go Z \og Y$ sometimes
written $X \times_Z Y$), and \demph{terminal objects} (written $1$,
especially for the terminal set $\{ * \}$); we also use
\demph{initial objects}.  A monad $(T,\eta,\mu)$ is often abbreviated to $T$.

I make no mention of the difference between sets and classes (`small
and large collections').  All the Definitions are really of \emph{small}
weak $n$-category.

Let \cat{C} be a category.  $X\in \cat{C}$ means that $X$ is an object of
$\cat{C}$, and $\cat{C}(X,Y)$ is the set of morphisms (or \demph{maps}, or
\demph{arrows}) from $X$ to $Y$ in \cat{C}.  If $f\in \cat{C}(X,Y)$ then $X$
is the \demph{domain} or \demph{source} of $f$, and $Y$ the \demph{codomain}
or \demph{target}.

\Set\ is the category (sets $+$ functions), and \Cat\ is (categories $+$
functors).  A set is just a \demph{discrete category} (one in which the only
maps are the identities).

$\cat{C}^\op$ is the \demph{opposite} or \demph{dual} of a category
\cat{C}.  $\ftrcat{\cat{C}}{\cat{D}}$ is the category of functors from
$\cat{C}$ to $\cat{D}$ and natural transformations between them.  Any object
$X$ of $\cat{C}$ induces a functor $\cat{C}(X, \dashbk): \cat{C} \go \Set$,
and a natural transformation from $\cat{C}(X, \dashbk)$ to
$F: \cat{C} \go \Set$ is the same thing as an element of $FX$ (the
\demph{Yoneda Lemma}); dually for $\cat{C}(\dashbk,X): \cat{C}^\op \go \Set$.

A functor $F: \cat{C} \go \cat{D}$ is an \demph{equivalence} if these
equivalent conditions hold: (i) $F$ is full, faithful and essentially
surjective on objects; (ii) there exist a functor $G: \cat{D} \go \cat{C}$ (a
\demph{pseudo-inverse} to $F$) and natural isomorphisms $\eta: 1 \go GF$,
$\epsln: FG \go 1$ ; (iii) as~(ii), but with $(F,G,\eta,\epsln)$ also being
an adjunction.

Any set $\cat{D}_0$ of objects of a category \cat{C} determines a \demph{full
subcategory} \cat{D} of \cat{C}, with object-set $\cat{D}_0$ and
$\cat{D}(X,Y) = \cat{C}(X,Y)$.  Every category \cat{C} has a
\demph{skeleton}: a subcategory whose inclusion into \cat{C} is an
equivalence and in which no two distinct objects are isomorphic.  If $F, G:
\cat{C} \go \Set$, $GX \sub FX$ for each $X \in \cat{C}$, and $F$ and $G$
agree on morphisms of \cat{C}, then $G$ is a
\demph{subfunctor} of $F$.

A \demph{total order} on a set $I$ is a reflexive transitive
relation $\leq$ such that if $i \neq j$ then exactly one of $i\leq j$ and
$j\leq i$ holds.  
$(I,\leq)$ can be seen as a category with object-set $I$
in which each hom-set has at most one element.
An \demph{order-preserving map} $(I,\leq) \go (I',\leq')$
is a function $f$ such that $i \leq j \implies f(i) \leq' f(j)$.

Let \Del\ be the category with objects $[k]=\{0,\ldots,k\}$ for $k\geq 0$,
and order-preserving functions as maps.  A \demph{simplicial set} is a
functor $\Delop \go \Set$.  Every category \cat{C} has a \demph{nerve} (the
simplicial set $N\cat{C}: [k] \goesto \Cat([k],\cat{C})$), giving a full and
faithful functor $N: \Cat \go \ftrcat{\Delop}{\Set}$.  So \Cat\ is equivalent
to the full subcategory of \ftrcat{\Delop}{\Set} with objects $\{ X \such X
\iso N\cat{C} \textrm{ for some } \cat{C} \}$; there are various
characterizations of such $X$, but we come to that in the main text.

Leftovers: a \demph{monoid} is a set (or more generally, an object of a
monoidal category) with an associative binary operation and a two-sided unit.
\Cat\ is monadic over the category of directed graphs.  The \demph{natural
numbers} start at $0$.

\clearpage

\concept{Strict $n$-Categories}

If \cat{V} is a category with finite products then there is a category
$\cat{V}\hyph\Cat$ of \cat{V}-enriched categories and \cat{V}-enriched
functors, and this itself has finite products.  (A \demph{\cat{V}-enriched
category} is just like an ordinary category, except that the `hom-sets' are
now objects of \cat{V}.)  Let $0\hyph\Cat = \Set$ and, for $n\geq 0$,
$(n+1)\hyph\Cat = (n\hyph\Cat)\hyph\Cat$; a \demph{strict $n$-category} is an
object of $n\hyph\Cat$.  Note that $1\hyph\Cat = \Cat$.

Any finite-product-preserving functor $U: \cat{V} \go \cat{W}$
induces a finite-product-preserving functor $U_*: \cat{V}\hyph\Cat \go
\cat{W}\hyph\Cat$, so we can define functors $U_n: (n+1)\hyph\Cat \go
n\hyph\Cat$ by taking $U_0$ to be the objects functor and $U_{n+1} =
(U_n)_*$.  The category $\omega\hyph\Cat$ of \demph{strict
$\omega$-categories} is the limit of the diagram
\[
\cdots 
\goby{U_{n+1}} (n+1)\hyph\Cat  \goby{U_n} 
\cdots
\goby{U_1} 1\hyph\Cat = \Cat
\goby{U_0} 0\hyph\Cat = \Set.
\]

Alternatively: a \demph{globular set} (or \demph{$\omega$-graph}) $A$
consists of sets and functions
\[
\cdots 
\parpair{s}{t} A_m  \parpair{s}{t} A_{m-1} \parpair{s}{t} 
\cdots 
\parpair{s}{t} A_0
\]
such that for $m\geq 2$ and $\alpha\in A_m$, $ss(\alpha) = st(\alpha)$ and
$ts(\alpha) = tt(\alpha)$.  An element of $A_m$ is called an
\demph{$m$-cell}, and we draw a $0$-cell $a$ as $\gzeros{a}$, a $1$-cell $f$
as $\gfsts{a}\gones{f}\glsts{b}$ (where $a=s(f), b=t(f)$), a 2-cell $\alpha$
as $\gfsts{a}\gtwos{f}{g}{\alpha}\glsts{b}$, etc.  For $m > p \geq 0$, write
$ A_m \times_{A_p} A_m = \{ (\alpha',\alpha) \in A_m \times A_m \such
s^{m-p}(\alpha') = t^{m-p}(\alpha) \}$.

A \demph{strict $\omega$-category} is a globular set $A$ together with a
function $\ofdim{p}: A_m \times_{A_p} A_m \go A_m$ (\demph{composition}) for
each $m > p \geq 0$ and a function $i: A_m \go A_{m+1}$ (\demph{identities},
usually written $i(\alpha) = 1_\alpha$) for each $m\geq 0$, such that
\begin{enumerate}
\item 	\label{part:strict-n:source-comp}
if $m > p \geq 0$ and
$(\alpha',\alpha) \in A_m \times_{A_p} A_m$ then
\[
\begin{array}{llll}
s(\alpha' \ofdim{p} \alpha) = 
s(\alpha) 			&	
\textrm{and}			&
t(\alpha' \ofdim{p} \alpha) = 
t(\alpha') 			&
\textrm{for }
m=p+1	\\
s(\alpha' \ofdim{p} \alpha) = 
s(\alpha') \ofdim{p} s(\alpha)	&
\textrm{and}			&
t(\alpha' \ofdim{p} \alpha) = 
t(\alpha') \ofdim{p} t(\alpha)	&
\textrm{for }
m\geq p+2	
\end{array}
\]
\item  	\label{part:strict-n:source-id}
if $m\geq 0$ and $\alpha\in A_m$ then $s(i(\alpha)) = \alpha =
t(i(\alpha))$ 
\item \label{part:strict-n:ass-and-id} 
if $m > p \geq 0$ and $\alpha \in A_m$ then $i^{m-p}(t^{m-p}(\alpha))
\ofdim{p} \alpha = \alpha = \alpha \ofdim{p}$\linebreak
$i^{m-p}(s^{m-p}(\alpha))$; if also $\alpha', \alpha''$ are such that
$(\alpha'', \alpha'), (\alpha', \alpha) \in A_m \times_{A_p} A_m$, then
$(\alpha'' \ofdim{p} \alpha') \ofdim{p} \alpha = \alpha'' \ofdim{p} (\alpha'
\ofdim{p} \alpha)$
\item  	\label{part:strict-n:int}
if $p>q\geq 0$ and $(f',f) \in A_p \times_{A_q} A_p$ then
$i(f') \ofdim{q} i(f) = i(f' \ofdim{q} f)$; 
if also $m>p$ and $\alpha,\alpha',\beta,\beta'$ are such that
$(\beta',\beta), (\alpha', \alpha) \in A_m \times_{A_p} A_m$, 
$(\beta',\alpha'), (\beta, \alpha) \in A_m \times_{A_q} A_m$,
then 
$(\beta' \ofdim{p} \beta) \ofdim{q} (\alpha' \ofdim{p} \alpha) 
= 
(\beta' \ofdim{q} \alpha') \ofdim{p} (\beta \ofdim{q} \alpha)$.
\end{enumerate}

The composition $\ofdim{p}$ is `composition of $m$-cells by gluing along
$p$-cells'.  Pictures for $(m,p) = (2,1), (1,0), (2,0)$ are in the
Bicategories section below. 

\demph{Strict $n$-categories} are defined similarly, but with the globular
set only going up to $A_n$.  \demph{Strict $n$- and $\omega$-functors} are
maps of globular sets preserving composition and identities; the categories
$n\hyph\Cat$ and $\omega\hyph\Cat$ thus defined are equivalent to the ones
defined above.  The comments below on the two alternative definitions of
bicategory give an impression of how this equivalence works. 

\clearpage

\concept{Bicategories}

Bicategories are the traditional and best-known formulation of `weak
2-category'.  

A \demph{bicategory} $B$ consists of
\begin{itemize}
\item a set $B_0$, whose elements $a$ are called \demph{0-cells} or
\demph{objects} of $B$ and drawn
$\gzeros{a}$
\item for each $a,b \in B_0$, a category $B(a,b)$, whose objects $f$ are
called \demph{1-cells} and drawn $\gfsts{a} \gones{f} \glsts{b}$, whose
arrows $\alpha: f \go g$ are called \demph{2-cells} and drawn $\gfsts{a}
\gtwos{f}{g}{\alpha} \glsts{b}$, and whose composition $\gfsts{a}
\gthrees{f}{g}{h}{\alpha}{\beta} \glsts{b} \goesto \gfsts{a}
\gtwos{f}{h}{\!\!\!\!\!\! \beta \sof \alpha} \glsts{b}$ is called
\demph{vertical composition} of 2-cells
\item for each $a \in B_0$, an object $1_a \in B(a,a)$ (the \demph{identity}
on $a$); and for each $a,b,c \in B_0$, a functor $B(b,c) \times B(a,b) \go
B(a,c)$, which on objects is called \emph{1-cell composition},
$\gfsts{a}\gones{f}\gblws{b}\gones{g}\glsts{c} \goesto 
\gfsts{a}\gones{g\sof f}\glsts{c}$, and on arrows is called \demph{horizontal
composition} of 2-cells, $\gfsts{a} \gtwos{f}{g}{\alpha} \gfbws{a'}
\gtwos{f'}{g'}{\alpha'} \glsts{a''} \goesto \gfsts{a} \gtwos{f' \sof f}{g' \sof
g}{\!\!\!\!\!\! \alpha' * \alpha} \glsts{a''}$ 
\item \demph{coherence 2-cells}: for each $f \in B(a,b), g \in B(b,c), h \in
B(c,d)$, an \demph{associativity isomorphism} $\xi_{h,g,f}: (h\of g)\of f
\go h\of (g\of f)$; and for each $f \in B(a,b)$, \demph{unit isomorphisms}
$\lambda_f: 1_b \of f \go f$ and $\rho_f: f \of 1_a \go f$
\end{itemize}
satisfying the following \demph{coherence axioms}:
\begin{itemize}
\item $\xi_{h,g,f}$ is natural in $h$, $g$ and $f$, and $\lambda_f$ and
$\rho_f$ are natural in $f$
\item if $f \in B(a,b), g \in B(b,c), h \in B(c,d), k
\in B(d,e)$, then 
$
\xi_{k,h,g\sof f} \,\of\, \xi_{k\sof h, g, f} = 
(1_k * \xi_{h,g,f}) \,\of\, \xi_{k,h\sof g,f} \,\of\, (\xi_{k,h,g} * 1_f)
$
(the \demph{pentagon axiom});
and if $f \in B(a,b), g \in B(b,c)$, then 
$
\rho_g * 1_f =
(1_g * \lambda_f) \,\of\, \xi_{g,1_b,f}
$
(the \demph{triangle axiom}).
\end{itemize}

An alternative definition is that a bicategory consists of sets and functions
$B_2 \parpair{s}{t} B_1 \parpair{s}{t} B_0$ satisfying $ss=st$ and $ts=tt$,
together with functions determining composition, identities and coherence
cells (in the style of the second definition of strict $\omega$-category
above).  The idea is that $B_m$ is the set of $m$-cells and that $s$ and $t$
give the source and target of a cell.  Strict 2-categories can be identified
with bicategories in which the coherence 2-cells are all identities.

A 1-cell $\gfsts{a}\gones{f}\glsts{b}$ in a bicategory $B$ is called an
\demph{equivalence} if there exists a 1-cell $\gfsts{b}\gones{g}\glsts{a}$
such that $g\of f \iso 1_a$ and $f\of g \iso 1_b$.  

A \demph{monoidal category} can be defined as a bicategory with only one
0-cell: for if the 0-cell is called $\star$ then the bicategory just consists
of a category $B(\star,\star)$ equipped with an object $I$, a functor
$\otimes: B(\star,\star)^2 \go B(\star,\star)$, and associativity and unit
isomorphisms satisfying coherence axioms.

We can consider \demph{strict functors} of bicategories, in which composition
etc is preserved strictly; more interesting are \demph{weak functors} $F$, in
which there are isomorphisms $Fg \of Ff \go F(g \of f)$, $1_{Fa} \go
F(1_a)$ satisfying coherence axioms.

\defnheading{Tr}	\label{p:tr}

\concept{Topological Background}

\paragraph{Spaces}
Let \Top\ be the category of topological spaces and continuous maps.  Recall
that compact spaces are exponentiable in \Top: that is, if $K$ is compact
then the set $Z^K$ of continuous maps from $K$ to a space
$Z$ can be given a topology (namely, the compact-open topology) in such a way
that there is an isomorphism $\Top(Y, Z^K) \iso \Top(K \times Y, Z)$
natural in $Y,Z \in \Top$.

\paragraph{Operads}
A \demph{(non-symmetric, topological) operad} $D$ is a sequence
$(D(k))_{k\geq 0}$ of spaces together with an element (the
\demph{identity}) of $D(1)$ and for each $k, r_1, \ldots, r_k \geq 0$ a map
\[
D(k) \times D(r_1) \times \cdots \times D(r_k) 	
\go						
D(r_1 + \cdots + r_k)
\]
(\demph{composition}), obeying unit and associativity laws.  (Example: fix an
object $M$ of a monoidal category \cat{M}, and define $D(k) =
\cat{M}(M^{\otimes k},M)$.)

\paragraph{The All-Important Operad}
There is an operad $E$ in which $E(k)$ is the space of continuous
endpoint-preserving maps from $[0,1]$ to $[0,k]$.  (`Endpoint-preserving'
means that $0$ maps to $0$ and $1$ to $k$.)  The identity element of $E(1)$
is the identity map, and composition in the operad is by substitution.

\paragraph{Path Spaces}
For any space $X$ and $x,x'\in X$, a \demph{path} from $x$ to $x'$ in $X$ is
a map $p: [0,1] \go X$ satisfying $p(0) = x$ and $p(1) = x'$.  There is a
space $X(x,x')$ of paths from $x$ to $x'$, a subspace of the exponential
$X^{[0,1]}$.

\paragraph{Operad Action on Path Spaces}
Fix a space $X$. For any $k\geq 0$ and $x_0, \ldots, x_k \in X$, there is a
canonical map
\[
\act{x_0,\ldots,x_k}: 
E(k) \times X(x_0,x_1) \times\cdots\times X(x_{k-1},x_k)
\go
X(x_0,x_k).
\]
These maps are compatible with the composition and identity of the operad
$E$,
and the construction is functorial in $X$.

\paragraph{Path-Components} Let $\Pi_0: \Top \go \Set$ be the functor
assigning to each space its set of path-components, and note that $\Pi_0$
preserves finite products.

\concept{The Definition}

We will define inductively, for each $n\geq 0$, a category \TCat{n} with
finite products and a functor $\Pi_n: \Top \go \TCat{n}$ preserving finite
products.  A \demph{weak $n$-category} is an object of \TCat{n}.
(Maps in \TCat{n} are to be thought of as \emph{strict}
$n$-functors.)

\paragraph{Base Case}

$\TCat{0} = \Set$, and $\Pi_0: \Top \go \Set$ is as above.

\paragraph{Objects of \TCat{(n+1)}} 

Inductively, a \demph{weak (n+1)-category} \pr{A}{\gamma} consists of
\begin{itemize}
\item a set $A_0$
\item a family $(A(a,a'))_{a,a'\in A_0}$ of weak $n$-categories
\item for each $k\geq 0$ and $a_0, \ldots, a_k \in A_0$, a map
\[
\gamma_{a_0, \ldots, a_k}:
\Pi_n(E(k)) \times A(a_0,a_1) \times\cdots\times A(a_{k-1},a_k)
\go
A(a_0,a_k)
\]
in \TCat{n},
\end{itemize}
such that the $\gamma_{a_0, \ldots, a_k}$'s satisfy compatibility axioms of
the same form as those satisfied by the $\act{x_0,\ldots,x_k}$'s.  (All this
makes sense because $\Pi_n$ preserves finite products and \TCat{n} has them.)

\paragraph{Maps in \TCat{(n+1)}} 

A \demph{map $\pr{A}{\gamma}\go\pr{B}{\delta}$} in \TCat{(n+1)} consists of
\begin{itemize}
\item a function $F_0: A_0 \go B_0$
\item for each $a,a'\in A_0$, a map $F_{a,a'}: A(a,a') \go B(F_0 a, F_0 a')$ of
weak $n$-categories,
\end{itemize}
satisfying the axiom
\[
F_{a_0,a_k} \,\of\, \gamma_{a_0, \ldots, a_k}
=
\delta_{F_0 a_0, \ldots, F_0 a_k} \,\of\, 
(1_{\Pi_n(E(k))} \times F_{a_0,a_1} \times\cdots\times F_{a_{k-1},a_k})
\]
for all $k\geq 0$ and $a_0, \ldots, a_k \in A_0$.

\paragraph{Composition and Identities in \TCat{(n+1)}}

Obvious.

\paragraph{$\Pi_{n+1}$ on Objects}

For a space $X$ we define $\Pi_{n+1}(X) = (A,\gamma)$, where 
\begin{itemize}
\item $A_0$ is the underlying set of $X$
\item $A(x,x') = \Pi_n(X(x,x'))$
\item for $x_0, \ldots, x_k \in X$, the map $\gamma_{x_0, \ldots, x_k}$ is the
composite 
\begin{eqnarray*}
\lefteqn{\Pi_n(E(k)) \times \Pi_n(X(x_0,x_1)) \times\cdots\times
\Pi_n(X(x_{k-1},x_k))} 							\\
	&\goiso	&\Pi_n(E(k) \times X(x_0,x_1) \times\cdots\times
		 X(x_{k-1},x_k))					\\
	&\goby{\Pi_n(\act{x_0,\ldots,x_k})}
		&\Pi_n(X(x_0,x_k)).
\end{eqnarray*}
\end{itemize}

\paragraph{$\Pi_{n+1}$ on Maps}

The functor $\Pi_{n+1}$ is defined on maps in the obvious way.

\paragraph{Finite Products Behave}

It is easy to show that \TCat{(n+1)} has finite products and that
$\Pi_{n+1}$ preserves finite products: so the inductive definition goes
through.

\clearpage

\lowdimsheading{Tr}

First observe that the space $E(k)$ is contractible for each $k$ (being, in a
suitable sense, convex).  In particular this tells us that $E(k)$ is
path-connected, and that the path space $E(k)(\theta,\theta')$ is
path-connected for every $\theta, \theta' \in E(k)$.

\concept{$n=0$}

By definition, $\TCat{0}=\Set$ and $\Pi_0: \Top\go\Set$ is the
path-components functor.

\concept{$n=1$}

\paragraph{The Category \TCat{1}}

A weak 1-category \pr{A}{\gamma} consists of
\begin{itemize}
\item a set $A_0$
\item a set $A(a,a')$ for each $a, a' \in A_0$
\item for each $k\geq 0$ and $a_0, \ldots, a_k \in A_0$, a function
\[
\gamma_{a_0, \ldots, a_k}:
\Pi_0(E(k)) \times A(a_0,a_1) \times\cdots\times A(a_{k-1},a_k)
\go
A(a_0,a_k)
\]
\end{itemize}
such that these functions satisfy certain axioms.  So a weak 1-category looks
something like a category: $A_0$ is the set of objects, $A(a,a')$ is
the set of maps from $a$ to $a'$, and $\gamma$ provides some kind of
composition.  Since $E(k)$ is path-connected, we may strike out $\Pi_0(E(k))$
from the product above; and then we may suggestively write
\[
(f_k \of\cdots\of f_1) = \gamma_{a_0, \ldots, a_k}(f_1, \ldots, f_k).
\]
The axioms on these `$k$-fold composition functions' mean that a weak
1-category is, in fact, exactly a category.  Maps in \TCat{1} are just
functors, and so $\TCat{1}$ is equivalent to $\Cat$.

\paragraph{The Functor $\Pi_1$} 

For a space $X$, the (weak 1-)category
$\Pi_1(X) = (A,\gamma)$ is given by
\begin{itemize}
\item $A_0$ is the underlying set of $X$
\item $A(x,x')$ is the set of path-components of the path-space $X(x,x')$:
that is, the set of homotopy classes of paths from $x$ to $x'$
\item Let $x_0 \goby{p_1} \ \cdots\ \goby{p_k} x_k$ be a sequence of paths in
$X$, and write $[p]$ for the homotopy class of a path $p$.  Then
\[
([p_k] \of \cdots \of [p_1]) = 
[\act{x_0,\ldots,x_k}(\theta, p_1, \ldots p_k)]
\]
where $\theta$ is any member of $E(k)$---it doesn't matter which.  In other
words, composition of paths is by laying them end to end.
\end{itemize}
Hence $\Pi_1(X)$ is the usual fundamental groupoid of $X$, and indeed $\Pi_1:
\Top \go$ $\Cat$ is the usual fundamental groupoid functor.

\concept{$n=2$} 

A weak 2-category \pr{A}{\gamma} consists of
\begin{itemize}
\item a set $A_0$
\item a category $A(a,a')$ for each $a, a' \in A_0$
\item for each $k\geq 0$ and $a_0, \ldots, a_k \in A_0$, a functor
\[
\gamma_{a_0, \ldots, a_k}:
\Pi_1(E(k)) \times A(a_0,a_1) \times\cdots\times A(a_{k-1},a_k)
\go
A(a_0,a_k)
\]
\end{itemize}
such that these functors satisfy axioms expressing compatibility with the
composition and identity of the operad $E$. 

By the description of $\Pi_1$ and the initial observations of this section,
the category $\Pi_1(E(k))$ is indiscrete (i.e.\ all hom-sets have one
element) and its objects are the elements of $E(k)$.  So $\gamma$ assigns to
each $\theta \in E(k)$ and $a_i\in A_0$ a functor
\[
\ovln{\theta}: 
A(a_0,a_1) \times\cdots\times A(a_{k-1},a_k)
\go
A(a_0,a_k),
\] 
and to each $\theta, \theta' \in E(k)$ and $a_i\in A_0$ a natural isomorphism
\[
\omega_{\theta,\theta'}: \ovln{\theta} \goiso \ovln{\theta'}.
\]
(Really we should add `$a_0, \ldots, a_k$' as a subscript to $\ovln{\theta}$
and to $\omega_{\theta,\theta'}$.)  Functoriality of $\gamma_{a_0, \ldots,
a_k}$ says that
\[
\omega_{\theta,\theta} = 1, \diagspace 
\omega_{\theta, \theta''} = 
\omega_{\theta',\theta''} \,\of\, \omega_{\theta,\theta'}.
\]
The `certain axioms' say firstly that 
\[
\ovln{\theta \,\of\, (\theta_1, \ldots, \theta_k)} =
\ovln{\theta} \,\of\, (\ovln{\theta_1} \times \cdots \times \ovln{\theta_k}),
\diagspace
\ovln{1} = 1
\]
for $\theta\in E(k)$ and $\theta_i\in E(r_i)$, where the left-hand sides of
the two equations refer respectively to composition and identity in the
operad $E$; and secondly that the natural isomorphisms
$\omega_{\theta,\theta'}$ fit together in a coherent way.

So a weak 2-category is probably not a structure with which we are already
familiar.  However, it nearly is.  For define $\tr(k)$ to be the set of
$k$-leafed rooted trees which are `unitrivalent' (each vertex has either 0 or
2 edges coming up out of it); and suppose we replaced $\Pi_1(E(k))$ by the
indiscrete category with object-set $\tr(k)$, so that the $\theta$'s above
would be trees.  A weak 2-category would then be exactly a bicategory: e.g.\
if $\theta=\littletree$ then $\ovln{\theta}$ is binary composition, and if
$(\theta,\theta') = (\lefttree,\righttree)$ then $\omega_{\theta,\theta'}$
is the associativity isomorphism.  And in some sense, a $k$-leafed tree might
be thought of as a discrete version of an endpoint-preserving map $[0,1] \go
[0,k]$.

With this in mind, any weak $2$-category $(A,\gamma)$ gives rise to a
bicategory $B$ (although the converse process seems less straightforward).
First pick at random an element $\theta_2$ of $E(2)$, and let $\theta_0$ be
the unique element of $E(0)$.  Then take $B_0 = A_0$, $B(a,a')=A(a,a')$,
binary composition to be $\ovln{\theta_2}$, identities to be
$\ovln{\theta_0}$, the associativity isomorphism to be $\omega_{\theta_2\sof
(1,\theta_2), \theta_2\sof (\theta_2,1)}$, and similarly units.  The coherence
axioms on $B$ follow from the coherence axioms on $\omega$: and so we have
a bicategory.

\clearpage

\defnheading{P}		\label{p:p}

\concept{Some Globular Structures}

\paragraph{Reflexive Globular Sets}

Let \scat{R} be the category whose objects are the natural numbers
$0,1,\ldots$, and whose arrows are generated by
\vspace{-4ex}
\[
\cdots \diagspace
m+1 \,
\splitcoeqlhs{\sigma_{m+1}}{\tau_{m+1}}{\iota_{m+1}} 
\, m \, 
\splitcoeqlhs{\sigma_m}{\tau_m}{\iota_m}
\diagspace \cdots \diagspace
\splitcoeqlhs{\sigma_1}{\tau_1}{\iota_1} 
0
\]
subject to the equations
\[
\sigma_m \of \sigma_{m+1} = \sigma_m \of \tau_{m+1},
\diagspace
\tau_m \of \sigma_{m+1} = \tau_m \of \tau_{m+1},
\diagspace
\sigma_m \of \iota_m = 1 = \tau_m \of \iota_m
\]
($m\geq 1$).  A functor $A: \scat{R} \go \Set$ is called a \demph{reflexive
globular set}.  I will write $s$ for $A(\sigma_m)$, and $t$ for $A(\tau_m)$,
and $1_a$ for $(A(\iota_m))(a)$ when $a \in A(m-1)$. 

\paragraph{Strict $\omega$-Categories, and $\omega$-Magmas}

A \demph{strict $\omega$-category} is a reflexive globular set $S$ together
with a function (\demph{composition}) $\ofdim{p}: S(m) \times_{S(p)} S(m) \go
S(m)$ for 
each $m > p \geq 0$, satisfying
\begin{itemize}
\item axioms determining the source and target of a composite
(part~\bref{part:strict-n:source-comp} in the Preliminary section `Strict
$n$-Categories')
\item 	
strict associativity, unit and interchange axioms (parts
~\bref{part:strict-n:ass-and-id} and~\bref{part:strict-n:int}). 
\end{itemize}

An \demph{$\omega$-magma} is like a strict $\omega$-category, but only
satisfying the first group of axioms~(\bref{part:strict-n:source-comp}) and
not necessarily the second
(\bref{part:strict-n:ass-and-id},~\bref{part:strict-n:int}).
A \demph{map of $\omega$-magmas} is a map of reflexive globular sets which
commutes with all the composition operations.  (A strict $\omega$-functor
between strict $\omega$-categories is, therefore, just a map of the
underlying $\omega$-magmas.)

\concept{Contractions}

Let $\phi: A \go B$ be a map of reflexive globular sets.  For $m\geq 1$,
define
\[
V_\phi(m) =
\{ (f_0,f_1) \in A(m) \times A(m) \such
s(f_0) = s(f_1), t(f_0) = t(f_1), \phi(f_0) = \phi(f_1) \},
\]
and define
\[
V_\phi(0) = 
\{ (f_0, f_1) \in A(0) \times A(0) \such \phi(f_0) = \phi(f_1) \}.
\]
A \demph{contraction} $\gamma$ on $\phi$ is a family of functions
\[
(\gamma_m : V_\phi(m) \go A(m+1))_{m\geq 0}
\]
such that for all $m\geq 0$ and $(f_0,f_1) \in V_\phi (m)$,
\[
s(\gamma_m(f_0,f_1)) = f_0, 
\ \ \ 
t(\gamma_m(f_0,f_1)) = f_1,
\ \ \ 
\phi(\gamma_m(f_0,f_1)) = 1_{\phi(f_0)} (= 1_{\phi(f_1)}),
\]
and for all $m\geq 0$ and $f \in A(m)$,
\[
\gamma_m(f,f) = 1_f.
\]

\concept{The Mysterious Category \cat{Q}}

\paragraph{Objects}

\begin{figure}
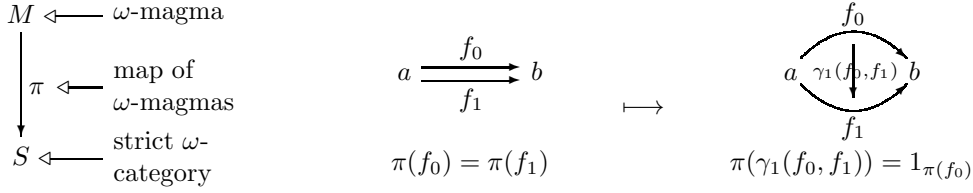

\[
\begin{diagram}[width=3em,height=2.7em,tight]
M		&\lLabelling	&\parbox{5em}{\raggedright
				$\omega$-magma}			\\
\dTo>{\pi}	&\lLabelling	&\parbox{5em}{\raggedright 
				map of $\omega$-magmas}		\\
S		&\lLabelling	&\parbox{5em}{\raggedright
				strict $\omega$-\\category}	\\
\end{diagram}
\mbox{\hspace{5em}}
\begin{array}{ccc}
\begin{diagram}[width=2.5em]
a&\pile{\rTo^{f_0}\\ \rTo_{f_1}}&b
\end{diagram} 				
&
\mbox{\hspace{1.6em}} \raisebox{-3ex}{\goesto} \mbox{\hspace{1.3em}}
&
a \ctwocentre{f_0}{f_1}{\,\scriptstyle \gamma_1(f_0,f_1)} b	\\
\pi(f_0) = \pi(f_1)			
&
&
\pi(\gamma_1(f_0,f_1)) = 1_{\pi(f_0)} 		
\end{array}
\]
\caption{An object of $\cat{Q}$, with $\gamma$ shown for $m=1$} 
\label{fig:object}
\end{figure}

An object of \cat{Q} (see Fig.~\ref{fig:object}) is a quadruple
$(M,S,\pi,\gamma)$ in which
\begin{itemize}
\item $M$ is an $\omega$-magma
\item $S$ is a strict $\omega$-category
\item $\pi$ is a map of $\omega$-magmas from $M$ to (the underlying
$\omega$-magma of) $S$
\item $\gamma$ is a contraction on $\pi$.
\end{itemize}

\paragraph{Maps}

A map $(M,S,\pi,\gamma) \go (M',S',\pi',\gamma')$ in \cat{Q} is a pair $(M
\goby{\chi} M',$ $S \goby{\zeta} S')$ commuting with everything in sight.
That is, $\chi$ is a map of $\omega$-magmas, $\zeta$ is a strict
$\omega$-functor, $\pi' \of \chi = \zeta \of \pi$, and $\gamma'_m(\chi(f_0),
\chi(f_1)) = \chi (\gamma_m(f_0, f_1))$ for all $(f_0, f_1) \in V_M(m)$.

\paragraph{Composition and Identities}

These are defined in the obvious way.

\concept{The Definition}

\paragraph{An Adjunction}

Let $U: \cat{Q} \go \ftrcat{\scat{R}}{\Set}$ be the functor sending
$(M,S,\pi,\gamma)$ to the underlying reflexive globular set of the
$\omega$-magma $M$.  It can be shown that $U$ has a left adjoint: so
there is an induced monad $T$ on \ftrcat{\scat{R}}{\Set}. 

\paragraph{Weak $\omega$-Categories}

A \demph{weak $\omega$-category} is a $T$-algebra.

\paragraph{Weak $n$-Categories}

Let $n\geq 0$.  A reflexive globular set $A$ is \demph{$n$-dimensional} if
for all $m \geq n$, the map $A(\iota_{m+1}): A(m) \go A(m+1)$ is an
isomorphism (and so $s = t = (A(\iota_{m+1}))^{-1}$).  A
\demph{weak $n$-category} is a weak $\omega$-category whose underlying
reflexive globular set is $n$-dimensional.

\clearpage

\lowdimsheading{P}

\concept{Direct Interpretation}

\paragraph{The Left Adjoint in Low Dimensions}

Here is a description of what the left adjoint $F$ to $U$ does in dimensions
$\leq 2$.  It is perhaps not obvious that $F$ as described does form the left
adjoint; we come to that later.

For a reflexive globular set $A$, write
\[
F(A) = 
\left(
\begin{diagram}[height=1.5em]
A^\#		\\
\dTo>{\pi_A}	\\
A^*		\\
\end{diagram},\ 
\gamma_A
\right).
\]
$A^*$ is, in fact, relatively easy to describe: it is the free strict
$\omega$-category on $A$, in which an $m$-cell is a formal pasting-together
of cells of $A$ of dimension $\leq m$.

\begin{description}
\item[Dimension $0$] We have $A^\#(0) = A^*(0) = A(0)$  and $(\pi_A)_0
= \id$.
\item[Dimension $1$] Next, $A^*(1)$ is the set of formal paths of 1-cells in
$A$, where we identify each identity cell $1_a$ with the identity path on
$a$.  The set $A^\#(1)$ and the functions $s,t: A^\#(1) \go A(0)$ are generated
by the following recursive clauses:
\begin{itemize}
\item if $a_0 \goby{f} a_1$ is a 1-cell in $A$ then $A^\#(1)$ contains an
element celled $f$, with $s(f)=a_0$ and $t(f)=a_1$
\item if $w,w' \in A^\#(1)$ with $t(w)=s(w')$ then $A^\#(1)$ contains an
element called $(w' \bofdim{0} w)$, with $s(w' \bofdim{0} w)=s(w)$ and $t(w'
\bofdim{0} w)=t(w')$. 
\end{itemize}
The identity map $A(0) \go A^\#(1)$ sends $a$ to $1_a \in A(1) \sub A^\#(1)$,
the map $\pi_A$ removes parentheses and sends $\bofdim{0}$ to $\ofdim{0}$,
and the contraction $\gamma_A$ is given by $\gamma_A(a,a) = 1_a$ (for $a\in
A(0)$).
\item[Dimension $2$] $A^*(1)$ is the set of formal pastings of
2-cells in $A$, again respecting the identities.  $A^\#(2)$ and $s, t:
A^\#(2) \go A^\#(1)$ are generated by:
\begin{itemize}
\item if $\alpha$ is a 2-cell in $A$ then $A^\#(2)$ has an element
called $\alpha$, with the evident source and target
\item if $a \parpair{w_0}{w_1} b$ in $A^\#(1)$ with $\pi_A(w_0)=\pi_A(w_1)$
then $A^\#(2)$ has an element called $\gamma_A(w_0,w_1)$, with source $w_0$ and
target $w_1$
\item if $x,x' \in A^\#(2)$ with $t(x)=s(x')$ then $A^\#(2)$ has an
element called $(x' \bofdim{1} x)$, with source $s(x)$ and target $t(x')$
\item if $x,x' \in A^\#(2)$ with $tt(x)=ss(x')$ then $A^\#(2)$ has an element
called $(x' \bofdim{0} x)$, with source $s(x')\bofdim{0} s(x)$ and target
$t(x')\bofdim{0} t(x)$;
\end{itemize}
furthermore, if $f \in A(1)$ then $1_f$ (from the first clause) is to be
identified with $\gamma_A(f,f)$ (from the second).  The identity map $A^\#(1)
\go A^\#(2)$ sends $w$ to $\gamma_A(w,w)$.  The map $\pi_A$ sends cells of
the form $\gamma_A(w_0,w_1)$ to identity cells, and otherwise acts as in
dimension 1.  The contraction $\gamma_A$ is defined in the way suggested by
the notation.

\end{description}

\paragraph{Adjointness}

We now have to see that this $F$ is indeed left adjoint to $U$.  First
observe that there is a natural embedding of $A(m)$ into $A^\#(m)$ (for
$m\leq 2$); this gives the unit of the adjunction.  Adjointness then says:
given $(M,S,\pi,\gamma) \in \cat{Q}$ and a map $A \goby{\phi} M$
of reflexive globular sets, there's a unique map
\[
(\chi,\zeta): \vslob{A^\#}{\pi_A}{A^*} \go \vslob{M}{\pi}{S}
\]
in \cat{Q} such that $\chi$ extends $\phi$.  This can be seen from the
description above.

\paragraph{Weak $2$-Categories}

A weak $2$-category consists of a $2$-dimensional reflexive globular set $A$
together with:
\begin{itemize}
\item (a map $A^\#(0) \go A(0)$ obeying axioms---which force it to be the
identity)
\item a map $A^\#(1) \go A(1)$ obeying axioms, which amounts to a binary
composition on the 1-cells of $A$ (\emph{not} obeying any axioms)
\item similarly, vertical and horizontal binary compositions of 2-cells, not
obeying any axioms `yet'
\item for each string $\cdot \goby{f_1} \cdots \goby{f_k} \cdot$ of 1-cells,
and each pair $\tau,\tau'$ of $k$-leafed binary trees, a 2-cell
$
\omega_{\tau,\tau'}: 
\of_{\tau}(f_1, \ldots, f_k) \go \of_{\tau'}(f_1, \ldots, f_k), 
$
where $\of_{\tau}$ indicates the iterated composition dictated by the shape
of $\tau$
\item amongst other things in dimension $3$: whenever we have some $2$-cells
$(\alpha_i)$, and two different ways of composing all the $\alpha_i$'s and
some $\omega_{\tau,\tau'}$'s to obtain new 2-cells $\beta$ and $\beta'$
respectively, and these satisfy $s(\beta)=s(\beta')$ and
$t(\beta)=t(\beta')$, then there is assigned a 3-cell $\beta\go\beta'$.
\end{itemize}
Since `the only $3$-cells of $A$ are equalities', we get $\beta=\beta'$ in
the last item.  Analysing this precisely, we find that the category of weak
2-categories is equivalent to the category of bicategories and strict
functors.  And more easily, a weak 1-category is just a category and a weak
0-category is just a set.

\concept{Indirect Interpretation}

An alternative way of handling weak $n$-categories is to work only with
$n$-dimensional (not infinite-dimensional) structures throughout: e.g.\
reflexive globular sets $A$ in which $A(m)$ is only defined for $m\leq n$.
We then only speak of contractions on a map $\phi$ if $(f_0,f_1) \in
V_\phi(n) \implies f_0=f_1$ (and in particular, the map $\pi$ must satisfy
this condition in order for $(M,S,\pi,\gamma)$ to qualify as an object of
\cat{Q}).  Our new category of weak $n$-categories appears to be equivalent
to the old one, taking algebra maps as the morphisms in both cases.

The analysis of $n=2$ is easier now: we can write down the left adjoint $F$
explicitly, and so get an explicit description of the monad $T$ on the
category of `reflexive 2-globular sets'.  This monad is presumably the free
bicategory monad.

\clearpage

\defnsheading{B}		\label{p:b}

\concept{Globular Operads and their Algebras}

\paragraph{Globular Sets}
 
Let \scat{G} be the category whose objects are the natural numbers
$0,1,\ldots$, and whose arrows are generated by
$
\sigma_m, \tau_m: m \go m-1
$
for each $m\geq 1$, subject to equations
\[
\sigma_{m-1} \of \sigma_m = \sigma_{m-1} \of \tau_m,
\diagspace
\tau_{m-1} \of \sigma_m = \tau_{m-1} \of \tau_m
\]
($m\geq 2$).  A functor $A: \scat{G} \go \Set$ is called a \demph{globular
set}; I will write $s$ for $A(\sigma_m)$, and $t$ for $A(\tau_m)$. 

\paragraph{The Free Strict $\omega$-Category Monad}

Any (small) strict $\omega$-category has an underlying globular set $A$, in
which $A(m)$ is the set of $m$-cells and $s$ and $t$ are the source and
target maps.  We thus obtain a forgetful functor $U$ from the category of
strict $\omega$-categories and strict $\omega$-functors to the category
\ftrcat{\scat{G}}{\Set} of globular sets.  $U$ has a left adjoint, so there
is an induced monad $(T, \id \goby{\eta} T, T^2 \goby{\mu} T)$ on
\ftrcat{\scat{G}}{\Set}.

\paragraph{Collections}

We define a monoidal category \fcat{Coll} of collections.  Let $1$ be the
terminal globular set.  A \demph{(globular) collection} is a map $C \goby{d}
T1$ into $T1$ in \ftrcat{\scat{G}}{\Set}; a \demph{map of collections} is a
commutative triangle.  The \demph{tensor product} of collections $C \goby{d}
T1$, $C' \goby{d'} T1$ is the composite along the top row of
\begin{diagram}[size=2em]
\SEpbk C \otimes C'&\rTo&TC'	&\rTo^{Td'}	&T^2 1	&\rTo^{\mu_1}&T1\\
\dTo		&	&\dTo>{T!}&		&	&	&	\\
C		&\rTo^d	&T1,	&		&	&	&	\\
\end{diagram}
where the right-angle symbol means that the square containing it is a
pullback, and $!$ denotes the unique map to $1$.  The \demph{unit} for the
tensor is $1 \goby{\eta_1} T1$.

\paragraph{Globular Operads}

A \demph{(globular) operad} is a monoid in the monoidal category \fcat{Coll};
a \demph{map of operads} is a map of monoids.  

\paragraph{Algebras}

Any operad $C$ induces a monad $C\cdot\dashbk$ on \ftrcat{\scat{G}}{\Set}.
For an object $A$ of \ftrcat{\scat{G}}{\Set}, this is defined by pullback: 
\begin{diagram}[size=2em]
\SEpbk C\cdot A &\rTo		&TA		\\
\dTo		&		&\dTo>{T!}	\\
C		&\rTo^d		&T1.		\\
\end{diagram}
The multiplication and unit of the monad come from the multiplication and
unit of the operad.  A $C$-\emph{algebra} is an algebra for the monad
$C\cdot\dashbk$.  Note that every $C$-algebra has an underlying globular
set.

\concept{Contractions and Systems of Composition}

\paragraph{Contractions}

Let $C \goby{d} T1$
be a collection.  For $m\geq 0$ and $\nu \in (T1)(m)$, write
$C(\nu) = \{ \theta \in C(m) \such d(\theta) = \nu \}$.  For $m\geq 1$ and
$\nu \in (T1)(m)$, define
\[
Q_C (\nu) =	\{ (\theta_0, \theta_1) \in C(\nu) \times C(\nu) \such
	s(\theta_0) = s(\theta_1) \mbox{ and } t(\theta_0) = t(\theta_1)\},
\]
and for $\nu \in (T1)(0)$, define $Q_C(\nu) = C(\nu) \times C(\nu)$.  Part of
the strict $\omega$-category structure on $T1$ is that each element $\nu \in
(T1)(m)$ gives rise to an element $1_\nu \in (T1)(m+1)$.  A
\demph{contraction} on $C$ is a family of functions
\[
(\gamma_\nu: Q_C (\nu) \go C(1_\nu))_{m\geq 0, \nu\in (T1)(m)} 
\]
satisfying
\[
s(\gamma_\nu \pr{\theta_0}{\theta_1} ) = \theta_0,
\diagspace
t(\gamma_\nu \pr{\theta_0}{\theta_1} ) = \theta_1
\]
for every $m\geq 0$, $\nu\in (T1)(m)$ and $\pr{\theta_0}{\theta_1} \in
Q_C (\nu)$. 

\paragraph{Systems of Compositions}

The map $\eta_1: 1 \go T1$ picks out an element $\eta_{1,m}$ of $(T1)(m)$
for each $m\geq 0$.  The strict $\omega$-category structure on $T1$ then
gives an element
\[
\beta_p^m = \eta_{1,m} \ofdim{p} \eta_{1,m} \in (T1)(m)
\]
for each $m > p \geq 0$; also put $\beta_m^m = \eta_{1,m}$.  Defining
$B(m) = \{ \beta_p^m \such m \geq p \geq 0 \} \sub (T1)(m)$, we obtain a
collection $B \rIncl T1$.

Also, the elements $\beta_m^m = \eta_{1,m} \in (T1)(m)$ determine a map $1
\go B$.

A \demph{system of compositions} in an operad $C$ is a map $B \go C$ of
collections such that the composite $1 \go B \go C$ is the unit of the
operad $C$.

\paragraph{Initial Object}  

Let \fcat{OCS} be the category in which an object is an operad equipped with
both a contraction and a system of compositions, and in which a map is a map
of operads preserving both the specified contraction and the specified system
of compositions.  Then \fcat{OCS} can be shown to have an initial object,
whose underlying operad will be written $K$.

\concept{The Definitions}

\paragraph{Definition \ds{B1}} 

A \demph{weak $\omega$-category} is a $K$-algebra.

\paragraph{Definition \ds{B2}} 

A \demph{weak $\omega$-category} is a pair $(C,A)$, where $C$ is an operad
satisfying $C(0) \iso 1$ and on which there exist a contraction and a system
of compositions, and $A$ is a $C$-algebra.

\paragraph{Weak $n$-Categories}

Let $n\geq 0$.  A globular set $A$ is \emph{$n$-dimensional} if for all
$m\geq n$,
\[
s=t: A(m+1) \go A(m)
\]
and this map is an isomorphism.  A \demph{weak $n$-category} is a weak
$\omega$-category whose underlying globular set is $n$-dimensional.  This can
be interpreted according to either \ds{B1} or \ds{B2}.

\clearpage

\lowsdimsheading{B}

\concept{Definition \ds{B1}}

An alternative way of handling weak $n$-categories is to work with only $n$-
(not infinite-) dimensional structures throughout.  So we replace $\scat{G}$
by its full subcategory $\scat{G}_n$ with objects $0, \ldots, n$, and $T$ by
the free strict $n$-category monad $T_n$, to obtain definitions of
\demph{$n$-collection}, \demph{$n$-operads}, and their \demph{algebras}.
\demph{Contractions} are defined as before, except that we only speak of
contractions on $C$ if
\begin{equation}	\label{eq:n-contr-b}
\forall \nu \in (T_n 1)(n), (\theta_0, \theta_1) \in Q_C(\nu)
\implies 
\theta_0 = \theta_1.
\end{equation}
There is an initial $n$-operad $K_n$ equipped with a contraction and a system
of compositions, and the category of weak $n$-categories turns out to be
equivalent to the category of $K_n$-algebras.  The latter is easier to
analyse.

\paragraph{$n=0$} 

We have $\ftrcat{\scat{G}_0}{\Set} \iso \Set$, $T_0 = \id$, and
$0\hyph\fcat{Coll} \iso \Set$; a 0-operad $C$ is a monoid, and a $C$-algebra
is a set with a $C$-action.  By~\bref{eq:n-contr-b}, the only 0-operad with a
contraction is the one-element monoid, so a weak $0$-category is just a set.

\paragraph{$n=1$}

$\ftrcat{\scat{G}_1}{\Set}$ is the category of directed graphs and $T_1$ is
the free category monad.  $K_1$ is the terminal 1-operad, by arguments
similar to those under `$n=2$' below.  It follows that the induced monad $K_1
\cdot \dashbk$ is just $T_1$, and so a weak $1$-category is just a
$T_1$-algebra, that is, a category.

\paragraph{$n=2$} 

A functor $A: \scat{G}_2 \go \Set$ consists of a set of $0$-cells (drawn
$\gzeros{a}$), a set of $1$-cells ($\gfsts{a}\gones{f}\glsts{b}$), and a set
of $2$-cells ($\gfsts{a}\gtwos{f}{g}{\alpha}\glsts{b}$).  A 2-collection $C$
consists of a set $C(0)$, a set $C(\nu_k)$ for each $k\geq 0$ (where $\nu_k$
indicates the `1-pasting diagram' $\gfstsu\gonesu\ \ldots\ \gonesu\glstsu$
with $k$ arrows), and a set $C(\pi)$ for each `2-pasting diagram' $\pi$ such
as the $\pi_i$ in Fig.~\ref{fig:op-comp-l}, together with source and target
functions.

\begin{figure}
\piccy{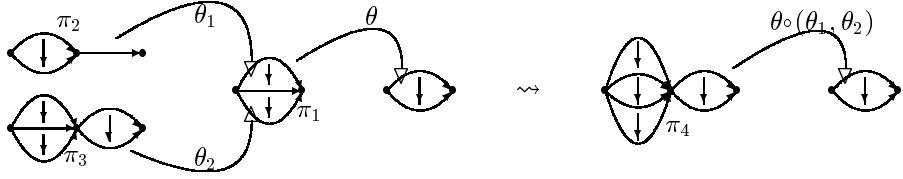}
\caption{Composition of operations in a globular operad}
\label{fig:op-comp-b}
\end{figure}

A 2-operad is a 2-collection $C$ together with `composition' functions such as 
\[
\begin{array}{ccc}
C(\nu_3) \times [C(\nu_2) \times_{C(0)} C(\nu_1) \times_{C(0)} C(\nu_2)] 
&\go &
C(\nu_5), \\
C(\pi_1) \times [C(\pi_2) \times_{C(\nu_2)} C(\pi_3)] &\go & C(\pi_4).
\end{array}
\]
In the first, the point is that there are $3$ terms $2,1,2$ and their sum is
$5$.  This makes sense if an element of $C(\nu_k)$ is regarded as an
operation which takes a string of $k$ 1-cells and turns it into a single
$1$-cell.  (The $\times_{C(0)}$'s denote pullbacks.)  Similarly for the
second; see Fig.~\ref{fig:op-comp-b}.  There are also identities for the
compositions.
A $C$-algebra is a functor $A: \scat{G}_2 \go \Set$ together with functions 
\[
\begin{array}{l}
\ovln{\psi}:	A(0) \go A(0) 
\textrm{ for each } \psi \in C(0),
\\
\ovln{\phi}: 	\{ \textrm{diagrams } \gfsts{a_0}\gones{f_1}\ \cdots \
\gones{f_k} \glsts{a_k} \textrm{ in } A \} \go A(1)			
\textrm{ for each } \phi \in C(\nu_k),
\\
\ovln{\theta}:	\{ \textrm{diagrams }
\gfsts{a}%
\gthrees{f}{g}{h}{\alpha}{\beta}%
\gfbws{\ b}\gtwos{l}{m}{\!\gamma}\glsts{c} 
\textrm{ in } A \} \go A(2)
\textrm{ for each } \theta \in C(\gfstsu\gthreesu\gzersu\gtwosu\glstsu)
\end{array}
\]
(etc), all compatible with the source, target, composition and identities in
$C$.

$K_2$ is generated from the empty collection by adding in the minimal amount
to obtain a 2-operad with contraction and system of compositions.  We have
the identity $1\in K_2(0)$.  Then, contraction gives an element (`1-cell
identities') of $K_2(\nu_0)$, the system of compositions gives an element
(`1-cell composition') of $K_2(\nu_2)$, and composition in $K_2$ gives one
element of $K_2(\nu_k)$ for each $k$-leafed tree in which every vertex has 0
or 2 edges coming up out of it.  Contraction at the next level gives
associativity and unit isomorphisms and identity 2-cells; the system of
compositions gives vertical and horizontal 2-cell composition.
Condition~\bref{eq:n-contr-b} gives coherence axioms.  Thus a weak 2-category
is exactly a bicategory.

\concept{Definition \ds{B2}}

Definition \ds{B2} of weak $n$-category refers to
infinite-dimensional globular operads.  So in order to do a concrete analysis
of $n\leq 2$, we redefine a \demph{weak $n$-category} as a pair $(C,A)$ where
$C$ is a $n$-operad admitting a contraction and a system of compositions and
with $C(0) \iso 1$, and $A$ is a $C$-algebra.  (Temporarily, call such a $C$
\demph{good}.)  I do not know to what extent this is equivalent to \ds{B2},
but the spirit, at least, is the same.

\paragraph{$n=0$} 

From `$n=0$' above we see that a weak $0$-category is just a set.

\paragraph{$n=1$}

The only good $1$-operad is the terminal $1$-operad, so by `$n=1$' above, a
weak $1$-category is just a category.

\paragraph{$n=2$}

A \demph{(non-symmetric) classical operad} $D$ is a sequence $(D(k))_{k\geq
0}$ of sets together with an element (the \demph{identity}) of $D(1)$ and for
each $k, r_1, \ldots, r_k \geq 0$ a map $ D(k) \times D(r_1) \times \cdots
\times D(r_k) \go D(r_1 + \cdots + r_k) $ (\demph{composition}), obeying unit
and associativity laws.  It turns out that good $2$-operads $C$ correspond
one-to-one with classical operads $D$ such that $D(k) \neq \emptyset$ for
each $k$, via $D(k) = C(\nu_k)$.  A $C$-algebra is then something like a
2-category or bicategory, with one way of composing a string of $k$ $1$-cells
for each element of $D(k)$, and all the appropriate coherence $2$-cells.
E.g.\ if $D=1$ then a $C$-algebra is a 2-category; if $D(k)$ is the set of
$k$-leafed trees in which each vertex has either $0$ or $2$ edges coming up
out of it then a $C$-algebra is a bicategory.  $C$ can, therefore, be
regarded as a theory of (more or less weak) 2-categories, and $A$ as a model
for such a theory.

\defnsheading{L}		\label{p:l}

\concept{Globular Operads and their Algebras}

\paragraph{Globular Sets}
 
Let \scat{G} be the category whose objects are the natural numbers
$0,1,\ldots$, and whose arrows are generated by
$
\sigma_m, \tau_m: m \go m-1
$
for each $m\geq 1$, subject to equations
\[
\sigma_{m-1} \of \sigma_m = \sigma_{m-1} \of \tau_m,
\diagspace
\tau_{m-1} \of \sigma_m = \tau_{m-1} \of \tau_m
\]
($m\geq 2$).  A functor $A: \scat{G} \go \Set$ is called a \demph{globular
set}; I will write $s$ for $A(\sigma_m)$, and $t$ for $A(\tau_m)$. 

\paragraph{The Free Strict $\omega$-Category Monad}

Any 
(small) 
strict $\omega$-category has an underlying globular set $A$, in
which $A(m)$ is the set of $m$-cells and $s$ and $t$ are the source and
target maps.  We thus obtain a forgetful functor $U$ from the category of
strict $\omega$-categories and strict $\omega$-functors to the category
\ftrcat{\scat{G}}{\Set} of globular sets.  $U$ has a left adjoint, so there
is an induced monad $(T, \id \goby{\eta} T, T^2 \goby{\mu} T)$ on
\ftrcat{\scat{G}}{\Set}.

\paragraph{Collections}

We define a monoidal category \fcat{Coll} of collections.  Let $1$ be the
terminal globular set.  A \demph{(globular) collection} is a map $C \goby{d}
T1$ into $T1$ in \ftrcat{\scat{G}}{\Set}; a \demph{map of collections} is a
commutative triangle.  The \demph{tensor product} of collections $C \goby{d}
T1$, $C' \goby{d'} T1$ is the composite along the top row of
\begin{diagram}[size=2em]
\SEpbk C \otimes C'&\rTo&TC'	&\rTo^{Td'}	&T^2 1	&\rTo^{\mu_1}&T1\\
\dTo		&	&\dTo>{T!}&		&	&	&	\\
C		&\rTo^d	&T1,	&		&	&	&	\\
\end{diagram}
where the right-angle symbol means that the square containing it is a
pullback, and $!$ denotes the unique map to $1$.  The \demph{unit} for the
tensor is $1 \goby{\eta_1} T1$.

\paragraph{Globular Operads}

A \demph{(globular) operad} is a monoid in the monoidal category \fcat{Coll};
a \demph{map of operads} is a map of monoids.  

\paragraph{Algebras}

Any operad $C$ induces a monad $C\cdot\dashbk$ on \ftrcat{\scat{G}}{\Set}.
For an object $A$ of \ftrcat{\scat{G}}{\Set}, this is defined by pullback: 
\begin{diagram}[size=2em]
\SEpbk C\cdot A &\rTo		&TA		\\
\dTo		&		&\dTo>{T!}	\\
C		&\rTo^d		&T1.		\\
\end{diagram}
The multiplication and unit of the monad come from the multiplication and
unit of the operad.  A \emph{$C$-algebra} is an algebra for the monad
$C\cdot\dashbk$.  Note that every $C$-algebra has an underlying globular
set.

\concept{Contractions}

\paragraph{Contractions}

Let $C \goby{d} T1$ be a collection.  For $m\geq 0$ and $\nu\in (T1)(m)$,
write 
$
C(\nu) = \{ \theta \in C(m) \such d(\theta) = \nu \}.
$ 
For $m\geq 2$ and $\pi\in (T1)(m)$, define
\[
P_C(\pi) = 
\{ \pr{\theta_0}{\theta_1} \in C(s(\pi)) \times C(t(\pi)) \such
s(\theta_0) = s(\theta_1) \mbox{ and } t(\theta_0) = t(\theta_1)\},
\]
and for $\pi\in (T1)(1)$, define
$
P_C(\pi) = C(s(\pi)) \times C(t(\pi)).
$
A \demph{contraction} $\gamma$ on $C$ is a family of functions
\[
(\gamma_\pi: P_C(\pi) \go C(\pi))_{m\geq 1, \pi\in (T1)(m)} 
\]
satisfying
\[
s(\gamma_\pi \pr{\theta_0}{\theta_1} ) = \theta_0,
\diagspace
t(\gamma_\pi \pr{\theta_0}{\theta_1} ) = \theta_1
\]
for every $m\geq 1$, $\pi\in (T1)(m)$ and $\pr{\theta_0}{\theta_1} \in
P_C(\pi)$. 

\paragraph{Initial Object}

Let \fcat{OC} be the category in which an object is an operad equipped with a
contraction and a map is a map of operads preserving the specified
contraction.  Then \fcat{OC} can be shown to have an initial object, whose
underlying operad will be written $L$.

\concept{The Definitions}

\paragraph{Definition \ds{L1}} 

A \demph{weak $\omega$-category} is an $L$-algebra.  (Maps of
$L$-algebras should be regarded as \emph{strict} $\omega$-functors.)

\paragraph{Definition \ds{L2}}

A \demph{weak $\omega$-category} is a pair $(C,A)$, where $C$ is an operad on
which there exists a contraction and satisfying $C(0) \iso 1$, and $A$ is a
$C$-algebra.

\paragraph{Weak $n$-Categories}

Let $n\geq 0$.  A globular set $A$ is \emph{$n$-dimensional} if for all
$m\geq n$,
\[
s=t: A(m+1) \go A(m)
\]
and this map is an isomorphism.  A \demph{weak $n$-category} is a weak
$\omega$-category whose underlying globular set is $n$-dimensional.  This can
be interpreted according to either \ds{L1} or \ds{L2}.

\clearpage

\lowsdimsheading{L}

\concept{Definition \ds{L1}}

An alternative way of handling weak $n$-categories is to work with only $n$-
(not infinite-) dimensional structures throughout.  So we replace $\scat{G}$
by its full subcategory $\scat{G}_n$ with objects $0, \ldots, n$, and $T$ by
the free strict $n$-category monad $T_n$, to obtain definitions of
\demph{$n$-collection}, \demph{$n$-operads}, and their \demph{algebras}.
\demph{Contractions} are defined as before, except that we only speak of
contractions on $C$ if
\begin{equation}	\label{eq:n-contr-l}
\forall \nu \in (T_n 1)(n), \forall \theta_0, \theta_1 \in C(\nu), \ 
s(\theta_0) = s(\theta_1) \ \&\ t(\theta_0) = t(\theta_1) \implies 
\theta_0 = \theta_1
\end{equation}
(taking $C(-1)=1$ to understand this when $n=0$).  There is an initial
$n$-operad $L_n$ equipped with a contraction, and the category of weak
$n$-categories turns out to be equivalent to the category of $L_n$-algebras.
The latter is easier to analyse.

\paragraph{$n=0$} 

We have $\ftrcat{\scat{G}_0}{\Set} \iso \Set$, $T_0 = \id$, and
$0\hyph\fcat{Coll} \iso \Set$; a 0-operad $C$ is a monoid, and a $C$-algebra
is a set with a $C$-action.  By~\bref{eq:n-contr-l}, the only
0-operad with a contraction is the one-element monoid, so a weak $0$-category
is just a set.

\paragraph{$n=1$}

$\ftrcat{\scat{G}_1}{\Set}$ is the category of directed graphs and $T_1$ is
the free category monad.  $L_1$ is the terminal 1-operad, by arguments
similar to those under `$n=2$' below.  It follows that the induced monad $L_1
\cdot \dashbk$ is just $T_1$, and so a weak $1$-category is just a
$T_1$-algebra, that is, a category.

\paragraph{$n=2$} 

A functor $A: \scat{G}_2 \go \Set$ consists of a set of $0$-cells (drawn
$\gzeros{a}$), a set of $1$-cells ($\gfsts{a}\gones{f}\glsts{b}$), and a set
of $2$-cells ($\gfsts{a}\gtwos{f}{g}{\alpha}\glsts{b}$).  A 2-collection $C$
consists of a set $C(0)$, a set $C(\nu_k)$ for each $k\geq 0$ (where $\nu_k$
indicates the `1-pasting diagram' $\gfstsu\gonesu\ \ldots\ \gonesu\glstsu$
with $k$ arrows), and a set $C(\pi)$ for each `2-pasting diagram' $\pi$ such
as the $\pi_i$ in Fig.~\ref{fig:op-comp-l}, together with source and target
functions.

A 2-operad is a 2-collection $C$ together with `composition' functions such as 
\[
\begin{array}{ccc}
C(\nu_3) \times [C(\nu_2) \times_{C(0)} C(\nu_1) \times_{C(0)} C(\nu_2)] 
&\go &
C(\nu_5), \\
C(\pi_1) \times [C(\pi_2) \times_{C(\nu_2)} C(\pi_3)] &\go & C(\pi_4).
\end{array}
\]
In the first, the point is that there are $3$ terms $2,1,2$ and their sum is
$5$.  This makes sense if an element of $C(\nu_k)$ is regarded as an
operation which takes a string of $k$ 1-cells and turns it into a single
$1$-cell.  (The $\times_{C(0)}$'s denote pullbacks.)  Similarly for the
second; see Fig.~\ref{fig:op-comp-l}.  There are also identities for the
compositions.
\begin{figure}
\piccy{compoppic.ps}
\caption{Composition of operations in a globular operad}
\label{fig:op-comp-l}
\end{figure}
A $C$-algebra is a functor $A: \scat{G}_2 \go \Set$ together with functions 
\[
\begin{array}{l}
\ovln{\psi}:	A(0) \go A(0) 
\textrm{ for each } \psi \in C(0),
\\
\ovln{\phi}: 	\{ \textrm{diagrams } \gfsts{a_0}\gones{f_1}\ \cdots \
\gones{f_k} \glsts{a_k} \textrm{ in } A \} \go A(1)			
\textrm{ for each } \phi \in C(\nu_k),
\\
\ovln{\theta}:	\{ \textrm{diagrams }
\gfsts{a}%
\gthrees{f}{g}{h}{\alpha}{\beta}%
\gfbws{\ b}\gtwos{l}{m}{\!\gamma}\glsts{c} 
\textrm{ in } A \} \go A(2)
\textrm{ for each } \theta \in C(\gfstsu\gthreesu\gzersu\gtwosu\glstsu)
\end{array}
\]
(etc), all compatible with the source, target, composition and identities in
$C$.

$L_2$ is generated from the empty collection by adding in the minimal amount
to obtain a 2-operad-with-contraction.  We have the identity $1\in L_2(0)$.
Then, contraction gives an element of $L_2(\nu_k)$ for each $k$, so that
composition in $L_2$ gives an element of $L_2(\nu_k)$ for each $k$-leafed
tree.  Contraction at the next level (with~\bref{eq:n-contr-l}) says that if
$\pi$ is a 2-pasting diagram of width $k$ then $L_2(\pi) = L_2(\nu_k) \times
L_2(\nu_k)$.  So $L_2(0) = \{1\}$, $L_2(\nu_k) = \{ k \textrm{-leafed trees}
\}$, $L_2(\pi) = \{ k \textrm{-leafed trees} \}^2$.  An $L_2$-algebra is,
then, an `unbiased bicategory': that is, just like a bicategory except that
there is specified $k$-fold composition for every $k\geq 0$ rather than just
$k=2$ (binary composition) and $k=0$ (identities).  Since these are
essentially the same as ordinary bicategories, so too are weak
$2$-categories.

\concept{Definition \ds{L2}}

Definition \ds{L2} of weak $n$-category refers to infinite-dimensional
globular operads.  So in order to do a concrete analysis of $n\leq 2$, we
redefine a \demph{weak $n$-category} as a pair $(C,A)$ where $C$ is a
$n$-operad admitting a contraction and with $C(0) \iso 1$, and $A$ is a
$C$-algebra.  (Temporarily, call such a $C$ \demph{good}.)  I do not know to
what extent this is equivalent to \ds{L2}, but the spirit, at least, is the
same.

\paragraph{$n=0$} 

From `$n=0$' above we see that a weak $0$-category is just a set.

\paragraph{$n=1$}

The only good $1$-operad is the terminal $1$-operad, so by `$n=1$' above, a
weak $1$-category is just a category.

\paragraph{$n=2$}

A \demph{(non-symmetric) classical operad} $D$ is a sequence $(D(k))_{k\geq
0}$ of sets together with an element (the \demph{identity}) of $D(1)$ and for
each $k, r_1, \ldots, r_k \geq 0$ a map $ D(k) \times D(r_1) \times \cdots
\times D(r_k) \go D(r_1 + \cdots + r_k) $ (\demph{composition}), obeying unit
and associativity laws.  It turns out that good $2$-operads $C$ correspond
one-to-one with classical operads $D$ such that $D(k) \neq \emptyset$ for
each $k$, via $D(k) = C(\nu_k)$.  A $C$-algebra is then something like a
2-category or bicategory, with one way of composing a string of $k$ $1$-cells
for each element of $D(k)$, and all the appropriate coherence $2$-cells.
E.g.\ if $D=1$ then a $C$-algebra is a 2-category; if $D(k)$ is the set of
$k$-leafed trees in which each vertex has either $0$ or $2$ edges coming up
out of it then a $C$-algebra is a bicategory.  $C$ can, therefore, be
regarded as a theory of (more or less weak) 2-categories, and $A$ as a model
for such a theory.

\defnheading{\lp}	\label{p:lprime}


\concept{Globular Multicategories}

\paragraph{Globular Sets}
 
Let \scat{G} be the category whose objects are the natural numbers
$0,1,\ldots$, and whose arrows are generated by
$
\sigma_m, \tau_m: m \go m-1
$
for each $m\geq 1$, subject to equations
\[
\sigma_{m-1} \of \sigma_m = \sigma_{m-1} \of \tau_m,
\diagspace
\tau_{m-1} \of \sigma_m = \tau_{m-1} \of \tau_m
\]
($m\geq 2$).  A functor $A: \scat{G} \go \Set$ is called a \demph{globular
set}; I will write $s$ for $A(\sigma_m)$, and $t$ for $A(\tau_m)$. 

\paragraph{The Free Strict $\omega$-Category Monad}

Any (small) strict $\omega$-category has an underlying globular set $A$, in
which $A(m)$ is the set of $m$-cells and $s$ and $t$ are the source and
target maps.  We thus obtain a forgetful functor $U$ from the category of
strict $\omega$-categories and strict $\omega$-functors to the category
\ftrcat{\scat{G}}{\Set} of globular sets.  $U$ has a left adjoint, so there
is an induced monad $(T, \id \goby{\eta} T, T^2 \goby{\mu} T)$ on
\ftrcat{\scat{G}}{\Set}.

\paragraph{Globular Graphs}

For each globular set $A$, we define a monoidal category $\Graph_A$.  An
object of $\Graph_A$ is a \demph{(globular) graph on $A$}: that is, a globular
set $R$ together with maps of globular sets
\begin{diagram}[size=2em]
	&		&R	&		&	\\
	&\ldTo<\dom	&	&\rdTo>\cod	&	\\
TA	&		&	&		&A.	\\
\end{diagram}
A \demph{map} $(R,\dom,\cod) \go (R',\dom',\cod')$ of graphs on $A$ is a map $R
\go R'$ making the evident triangles commute.  The \demph{tensor product} of
graphs $(R,\dom,\cod)$, $(R',\dom',\cod')$ is given by composing along the
upper edges of the following diagram, in which the right-angle symbol means that the square containing it is a
pullback:
\begin{diagram}[size=2em]
&&	&	&	&	&R\otimes R'\Spbk&&	&	&	\\
&&	&	&	&\ldTo	&	&\rdTo	&	&	&	\\
&&	&	&TR'	&	&	&	&R	&	&	\\
&&	&\ldTo<{T\dom'}&&\rdTo<{T\cod'}&&\ldTo>\dom&	&\rdTo>\cod&	\\
&&T^2A	&	&	&	&TA	&	&	&	&A.	\\
&\ldTo<{\mu_A}&&&	&	&	&	&	&	&	\\
TA&&	&	&	&	&	&	&	&	&	\\
\end{diagram}
The \demph{unit} for the tensor is the graph
\begin{diagram}[size=2em]
	&		&A	&		&	\\
	&\ldTo<{\eta_A}	&	&\rdTo>1	&	\\
TA	&		&	&		&A.	\\
\end{diagram}

\paragraph{Globular Multicategories}  

A \demph{globular multicategory} is a globular set $A$ together with a
monoid in $\Graph_A$.  A globular multicategory $\cat{A}$ therefore consists of a
globular set $A$, a graph $(R,\dom,\cod)$ on $A$, and maps $\comp: R \otimes
R \go R$ and $\ids: A \go R$ compatible with \dom\ and \cod\ and obeying associativity and identity laws.

\concept{Contractible Maps}

\begin{figure}
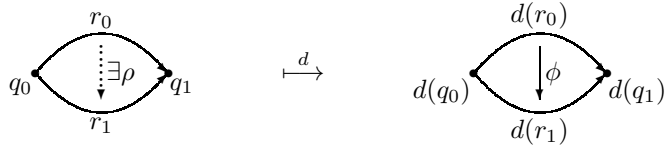

\begin{center}
$
\gfst{q_0}\gtwodotty{r_0}{r_1}{\exists\rho}\glst{q_1}
\mbox{\hspace{2.5em}}
\stackrel{d}{\goesto}
\mbox{\hspace{2.5em}}
\gfst{d(q_0)}\gtwo{d(r_0)}{d(r_1)}{\phi}\glst{d(q_1)}
$
\end{center}
\vspace*{-1em}
\caption{Part~\bref{part:major} of the definition of contractibility, shown
for $m=1$}
\label{fig:contr}
\end{figure}

A map $d: R \go S$ of globular sets is \demph{contractible}
(Figure~\ref{fig:contr}) if
\begin{enumerate}
\item the function $d_0: R(0) \go S(0)$ is bijective, and
\item 	\label{part:major}
for every 
\begin{description}
\item $m\geq 0$,
\item $r_0,r_1 \in R(m)$ with $s(r_0)=s(r_1)$ and $t(r_0)=t(r_1)$,
\item $\phi \in S(m+1)$ with $s(\phi)=d_m(r_0)$ and $t(\phi)=d_m(r_1)$,
\end{description}
there exists $\rho \in R(m+1)$ with $s(\rho)=r_0$, $t(\rho)=r_1$, and
$d_{m+1}(\rho)=\phi$.  In the case $m=0$ we drop the (nonsensical) conditions
that $s(r_0)=s(r_1)$ and $t(r_0)=t(r_1)$.
\end{enumerate}

\concept{The Definition}

\paragraph{Weak $\omega$-Categories}

A \demph{weak $\omega$-category} is a globular multicategory 
$\cat{A} = (A,R,\dom,\cod,\comp,\ids)$ such that $\dom: R \go TA$ is
contractible. 

\paragraph{Weak $n$-Categories}

Let $n\geq 0$.  A globular set $A$ is \demph{$n$-dimensional} if for all
$m\geq n$,
\[
s=t: A(m+1) \go A(m)
\]
and this map is an isomorphism.  A \demph{weak $n$-category} is a weak
$\omega$-category $\cat{A}$ such that the globular sets $A$ and $R$ are
$n$-dimensional.

\clearpage

\lowdimsheading{\lp}

An alternative way of handling weak $n$-categories is to work with only
$n$-dimensional (not infinite-dimensional) structures throughout.  Thus we
replace $\scat{G}$ by its full subcategory $\scat{G}_n$ with objects $0,
\ldots, n$, replace $T$ by the free strict $n$-category monad $T_n$, and so
obtain a definition of \demph{$n$-globular multicategory}.  We also modify
part~\bref{part:major} of the definition of contractibility by changing
`$m\geq 0$' to `$n-1\geq m\geq 0$', and `there exists $\rho$' to `there
exists a unique $\rho$' in the case $m=n-1$.  From these ingredients we get a
new definition of weak $n$-category.

The new and old definitions give two different, but equivalent, categories of
weak $n$-categories (with maps of multicategories as the morphisms); the
analysis of $n\leq 2$ is more convenient with the new definition.

\concept{$n=0$}

We have $\ftrcat{\scat{G}_0}{\Set} \iso \Set$ and $T_0=\id$, and the
contractible maps are the bijections.  So a weak $0$-category is a category
whose domain map is a bijection; that is, a discrete category; that is,
a set.

\concept{$n=1$}

$\ftrcat{\scat{G}_1}{\Set}$ is the category of directed graphs, $T_1$ is the
free category monad on it, and a map of graphs is contractible if and only if
it is an isomorphism.  So a weak $1$-category is essentially a 1-globular
multicategory whose underlying 1-globular graph looks like $T_1 A \ogby{1}
T_1 A \goby{\cod} A$.  Such a graph has at most one multicategory structure,
and it has one if and only if $\cod$ is a $T_1$-algebra structure on $A$.  So
a weak $1$-category is just a $T_1$-algebra, i.e., a category.

\concept{$n=2$}

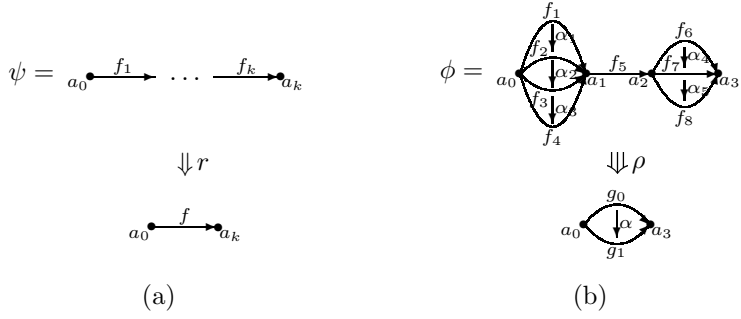
\begin{figure}
\centering
\setlength{\unitlength}{1mm}
\begin{picture}(40,40)
\cell{0}{30}{l}{\psi = 
	\gfsts{a_0} \gones{f_1} \ \ldots \ \gones{f_k} \glsts{a_k}}
\cell{23.5}{18}{c}{\Downarrow}
\cell{25}{18}{l}{r}
\cell{23.5}{10}{c}{\gfsts{a_0} \gones{f} \glsts{a_k}}
\cell{20}{2}{t}{\textrm{(a)}}
\end{picture}
\hspace{15mm}
\begin{picture}(40,40)
\cell{0}{30}{l}{\phi = 
	\gfsts{a_0} \gfours{f_1}{f_2}{f_3}{f_4}{\alpha_1}{\alpha_2}{\alpha_3}
	\grgts{a_1} \gones{f_5}
	\glfts{a_2} \gthrees{f_6}{f_7}{f_8}{\alpha_4}{\alpha_5}
	\glsts{a_3}}
\cell{24.5}{20.8}{c}{\Ddownarrow}
\cell{25.5}{18}{l}{\rho}
\cell{23.5}{10}{c}{\gfsts{a_0} \gtwos{g_0}{g_1}{\alpha} \glsts{a_3}}
\cell{20}{2}{t}{\textrm{(b)}}
\end{picture}
\caption{(a) A 1-cell, and (b) a typical 2-cell, of $R$.  Here $a_i, f_i, f,
g_i, \alpha_i$ and $\alpha$ are all cells of $A$}
\label{fig:reasons}
\end{figure}

The free 2-category $T_2 A$ on a 2-globular set $A \in
\ftrcat{\scat{G}_2}{\Set}$ has the same 0-cells as $A$; 1-cells of $T_2 A$
are formal paths $\psi$ in $A$ as in Fig.~\ref{fig:reasons}(a); and a typical
2-cell of $T_2 A$ is the diagram $\phi$ in Fig.~\ref{fig:reasons}(b).

Next, what is a 2-globular multicategory
$\cat{A}=(A,R,\dom,\cod,\comp,\ids)$?  Since we ultimately want to consider
just those \cat{A} in which \dom\ is contractible, let us assume immediately
that $R(0)=A(0)$.  Then \cat{A} consists of:
\begin{itemize}
\item a 2-globular set $A \in \ftrcat{\scat{G}_2}{\Set}$
\item for each $\psi$ and $f$ as in Fig.~\ref{fig:reasons}(a), a set of cells
$r: \psi \reason f$; such an $r$ is a 1-cell of $R$, and can be
regarded as a `reason why $f$ is a composite of $\psi$'
\item for each $\phi$ and $g_0,g_1,\alpha$ as in Fig.~\ref{fig:reasons}(b), a
set of cells $\rho: \phi\Reason\alpha$; such a $\rho$ is a
2-cell of $R$, and can be regarded as a `reason why $\alpha$ is a
composite of $\phi$'
\item source and target functions $R(2) \parpairu R(1)$, which, for instance,
assign to $\rho$ a reason $s(\rho)$ why $g_0$ is a composite of
$\gfsts{a_0}\gones{f_1}\gblws{a_1}\gones{f_5}
\gblws{a_2}\gones{f_6}\glsts{a_3}$
\item composition and identities: given $r$ as in Fig.~\ref{fig:reasons}(a)
and similarly $r_i: (f_i^1, \ldots, f_i^{p_i}) \reason f_i$ for each $i=1,
\ldots, k$, there is a composite $r\,\of\,(r_1, \ldots, r_k): (f_1^1, \ldots,
f_k^{p_k}) \reason f$; and similarly for 2-cells and for identities,
\end{itemize}
such that the composition and identities satisfy associativity and identity
axioms and are compatible with source and target.

Contractibility says that for each $\psi$ as in
Fig.~\ref{fig:reasons}(a) there is at least one pair $(r,f)$ as in
Fig.~\ref{fig:reasons}(a), and that for each $\phi$ as in
Fig.~\ref{fig:reasons}(b) and each $r_0: (f_1,f_5,f_6) \reason g_0$ and $r_1:
(f_4,f_5,f_8) \reason g_1$, there is exactly one pair $(\rho,\alpha)$ as in
Fig.~\ref{fig:reasons}(b) satisfying $s(\rho)=r_0$ and $t(\rho)=r_1$.  That
is: every diagram of 1-cells has at least one composite, and every diagram of
2-cells has exactly one composite once a way of composing the 1-cells
along its boundary has been chosen.

When $\phi = \gfsts{a_0} \gthrees{f_0}{f_1}{f_2}{\alpha_1}{\alpha_2}
\glsts{a_1}$, the identity reasons for $f_0$ and $f_2$ give via
contractibility a composite $\alpha_2 \of \alpha_1: f_0 \go f_2$, and in this
way the 1- and 2-cells between $a_0$ and $a_1$ form a category
$\cat{A}(a_0,a_1)$.  Now suppose that $\psi$ is as in
Fig.~\ref{fig:reasons}(a) and $r:\psi\reason f$, $r':\psi\reason f'$.
Applying contractibility to the degenerate 2-cell diagram $\phi$ which looks
exactly like $\psi$, we obtain a 2-cell
$\gfsts{a_0}\gtwos{f}{f'}{}\glsts{a_k}$; and similarly the other way round;
so by the uniqueness property of the $\rho$'s, $f\iso f'$ in
$\cat{A}(a_0,a_k)$.  Thus any two composites of a string of 1-cells are
canonically isomorphic.

A weak 2-category is essentially what is known as an `anabicategory'.  To see
how one of these gives rise to a bicategory, choose for each
$\gfsts{a_0}\gones{f}\gblws{a_1}\gones{g}\glsts{a_2}$ in $A$ a reason
$r_{f,g}: (f,g) \reason h$ and write $h=(g\of f)$; and similarly for
identities.  Then, for instance, the horizontal composite of 2-cells $
\gfsts{a_0} \gtwos{f_0}{f_1}{} \gfbws{a_1} \gtwos{g_0}{g_1}{} \glsts{a_2} $
comes via contractibility from $r_{f_0,g_0}$ and $r_{f_1,g_1}$, the
associativity cells arise from the degenerate 2-cell diagram $ \phi =
\gfsts{a_0} \gones{f_1}\gblws{a_1} \gones{f_2}\gblws{a_2}
\gones{f_3}\glsts{a_3} $, and the coherence axioms come from the uniqueness
of the $\rho$'s.

\defnheading{Si}		\label{p:si}

\concept{Simplicial Objects}

\paragraph{The Simplicial Category}

Let \Del\ be a skeleton of the category of nonempty finite totally ordered
sets: that is, \Del\ has objects $[k]=\{0,\ldots,k\}$ for $k\geq 0$, and maps
are order-preserving functions (with respect to the usual ordering $\leq$).

\paragraph{Some Maps in \Del}

Let $\sigma, \tau: [0] \go [1]$ be the maps in \Del\ with respective values
$0$ and $1$.  Given $k\geq 0$, let $\iota_1, \ldots, \iota_k: [1] \go [k]$
denote the `embeddings' of $[1]$ into $[k]$, defined by
$\iota_j(0) = j-1$ and $\iota_j(1) = j$.

\paragraph{The Segal Maps}

Let $k\geq 0$.  Then the following diagram in \Del\ commutes:
\begin{diagram}[width=2em,height=2em]
	&	&	&	&	&	&	&[k]	&	&
	&	&	&	&	&	\\
	&	&	&	&	&	&\ruTo(7,2)<{\iota_1}	
						 \ruTo(3,2)<{\iota_2}
							&	&\ldots\luTo(7,2)>{\iota_k}	&
	&	&	&	&	&	\\
[1]	&	&	&	&[1]	&	&	&	&	&
\ldots	&	&	&	&	&[1]	\\
	&\luTo<\tau&	&\ruTo>\sigma&	&\luTo<\tau&	&\ruTo>\sigma&	&
	&	&\luTo<\tau&	&\ruTo>\sigma&	\\
	&	&[0]	&	&	&	&[0]	&	&	&
\cdots	&	&	&[0].	&	&	\\
\end{diagram}
Let $X: \Delop \go \cat{E}$ be a functor into a category \cat{E} possessing
finite limits, and write $ X[1] \times_{X[0]} \cdots \times_{X[0]} X[1] $
(with $k$ occurrences of $X[1]$) for the limit of the diagram
\begin{diagram}[width=2em,height=2em]
X[1]	&	&	&	&X[1]	&	&	&	&	&
\ldots	&	&	&	&	&X[1]	\\
	&\rdTo<{X\tau}&	&\ldTo>{X\sigma}&	&\rdTo<{X\tau}&	&\ldTo>{X\sigma}&	&
	&	&\rdTo<{X\tau}&	&\ldTo>{X\sigma}&	\\
	&	&X[0]	&	&	&	&X[0]	&	&	&
\cdots	&	&	&X[0]	&	&	\\
\end{diagram}
(with, again, $k$ occurrences of $X[1]$) in \cat{E}.  Then by commutativity
of the first diagram, there is an induced map in \cat{E}---a \demph{Segal
map}---
\begin{equation}	\label{eq:Segal-Si}
X[k] \go X[1] \times_{X[0]} \cdots \times_{X[0]} X[1].
\end{equation}

\concept{Contractibility}

\paragraph{Sources and Targets}

If $0\leq p\leq r$, write $I_p$ for the object $ (\underbrace{[1], \ldots,
[1]}_p , \underbrace{[0], \ldots, [0]}_{r-p}) $ of $\Deln{r}$.  Let $X:
\Delnop{r} \go \Set$, $0 \leq p \leq r$, and $x, x' \in X(I_p)$.  Then
$x, x'$ are \demph{parallel} if $p=0$ or if $p\geq 1$ and $s(x) =
s(x')$ and $t(x) = t(x')$; here $s$ and $t$ are the maps
\begin{diagram}
X(I_p)	&
\pile{\rTo^{X(\id, \ldots, \id, X\sigma, \id, \ldots, \id)}\\
\rTo_{X(\id, \ldots, \id, X\tau, \id, \ldots, \id)}}	&
X(I_{p-1}).\\
\end{diagram}

\paragraph{Contractible Maps}

Let $r\geq 1$ and let $X,Y: \Delnop{r} \go \Set$.  A natural transformation
$\phi: X \go Y$ is \demph{contractible} if
\begin{itemize}
\item 	
the function $\phi_{I_0}: X(I_0) \go Y(I_0)$ is surjective
\item  	
given $p\in \{0, \ldots, r-1\}$, parallel $x,x'\in X(I_p)$, and $h\in
Y(I_{p+1})$ satisfying
\[
s(h)=\phi_{I_p}(x), \diagspace
t(h)= \phi_{I_p}(x'), \diagspace
\]
\marginpar{\centering
\fbox{%
\begin{diagram}[width=2em,height=1.3em]
x	&\rGet^g	&x'	\\
	&		&	\\
	&\dGoesto>\phi	&	\\
	&		&	\\
\phi(x)	&\rTo^h		&\phi(x')\\
\end{diagram}}}
there exists $g\in X(I_{p+1})$ satisfying 
\[
s(g)=x, \diagspace
t(g)=x',\diagspace
\phi_{I_{p+1}}(g)=h
\]
\item  	
given parallel $x,x'\in X(I_r)$ satisfying $\phi_{I_r}(x) = \phi_{I_r}(x')$,
then $x=x'$. 
\end{itemize}

If $r=0$ then $X$ and $Y$ are just sets and $\phi$ is just a function $X\go
Y$; call $\phi$ \demph{contractible} if it is bijective.

\concept{The Definition}

Let $n\geq 0$.  A \demph{weak $n$-category} is a functor $A: \Delnop{n} \go
\Set$ such that for each $m\in \{0, \ldots, n-1\}$ and $K=([k_1], \ldots,
[k_m]) \in \Deln{m}$,
\begin{enumerate}
\item 	\label{part:defn:degen-Si}
the functor
$
A(K, [0], \dashbk): \Delnop{n-m-1} \go \Set
$ 
is constant, and
\item 	\label{part:defn:main-Si} 
for each $[k]\in\Del$, the Segal map
\[
A(K, [k], \dashbk) \go 
A(K, [1],\dashbk) \times_{A(K, [0],\dashbk)} \cdots 
\times_{A(K, [0],\dashbk)} A(K, [1],\dashbk)
\]
is contractible.  (We are taking
$\cat{E}=\ftrcat{\Delnop{n-m-1}}{\Set}$ and $X[j]=A(K,[j],\dashbk)$ in the
definition of Segal map.)
\end{enumerate}

\clearpage

\lowdimsheading{Si}

\concept{$n=0$}

Parts~\bref{part:defn:degen-Si} and~\bref{part:defn:main-Si} of the
definition are vacuous, so a weak $0$-category is just a functor
$\Delnop{0}\go\Set$, that is, a set.

\concept{$n=1$}

A weak $1$-category is a functor $A: \Delop \go \Set$ (that is, a simplicial
set) satisfying~\bref{part:defn:degen-Si} and~\bref{part:defn:main-Si}.
Part~\bref{part:defn:degen-Si} is always true, and~\bref{part:defn:main-Si}
says that for each $k\geq 0$ the Segal map~\bref{eq:Segal-Si} (with $X=A$) is
a bijection---in other words, that $A$ is a \demph{nerve}.  The category of
nerves and natural transformations between them is equivalent to \Cat, where
a nerve $A$ corresponds to a certain category with object-set $A[0]$ and
morphism-set $A[1]$.  So a weak $1$-category is essentially just a category.

\concept{$n=2$}

A weak $2$-category is a
functor $A: \Delnop{2} \go \Set$ such that
\begin{enumerate}
\item the functor $A([0],\dashbk): \Delop \go \Set$ is constant
\item for each $k\geq 0$, the Segal map 
\[
A([k], \dashbk) \go 
A([1],\dashbk) \times_{A([0],\dashbk)} \cdots 
\times_{A([0],\dashbk)} A([1],\dashbk)
\]
is contractible, and for each $k_1, k\geq 0$, the Segal map
\[
A([k_1],[k]) \go 
A([k_1],[1]) \times_{A([k_1],[0])} \cdots 
\times_{A([k_1],[0])} A([k_1],[1])
\]
is a bijection.
\end{enumerate}
The second half of~\bref{part:defn:main-Si} says that $A([k_1],\dashbk)$ is a
nerve for each $k_1$, so we can regard $A$ as a functor $\Delop \go \Cat$.
Note that if $X$ and $Y$ are nerves then a natural transformation $\phi: X
\go Y$ is the same thing as a functor between the corresponding categories,
and that $\phi$ is contractible if and only if this functor is full, faithful
and surjective on objects.  So a weak $2$-category is a functor $A: \Delop
\go \Cat$ such that
\begin{enumerate}
\item $A[0]$ is a discrete category (i.e.\ the only morphisms are the
identities) 
\item for each $k\geq 0$, the Segal functor 
\[
A[k] \go A[1] \times_{A[0]} \cdots \times_{A[0]} A[1]
\]
is full, faithful and surjective on objects.  
\end{enumerate}

I will now argue that a weak $2$-category is essentially the same thing as a
bicategory.

First take a weak $2$-category $A: \Delop \go \Cat$, and let us construct a
bicategory $B$.  The object-set of $B$ is $A[0]$.  The two functors
$s,t: A[1] \go A[0]$ express the category $A[1]$ as a disjoint union
$\coprod_{a,b\in A[0]}B(a,b)$ of categories; the $1$-cells from $a$ to $b$
are the objects of $B(a,b)$, and the $2$-cells are the morphisms.

Vertical composition of $2$-cells in $B$ is composition in each
$B(a,b)$.  To define horizontal composition of $1$-~and $2$-cells,
first choose for each $k$ a pseudo-inverse
\[
A[1] \times_{A[0]} \cdots \times_{A[0]} A[1] \goby{\psi_k} A[k]
\]
to the Segal functor $\phi_k$ (an equivalence of categories), and
natural isomorphisms $\eta_k: 1 \go \psi_k\,\of\,\phi_k$, $\epsln_k:
\phi_k\,\of\,\psi_k \go 1$.  Horizontal composition is given by
\[
A[1] \times_{A[0]} A[1] \goby{\psi_2} A[2] \goby{A\delta} A[1],
\]
where $\delta:[1] \go [2]$ is the injection whose image omits $1\in [2]$.
The associativity isomorphisms are built up from $\eta_k$'s and $\epsln_k$'s,
and the pentagon commutes just as long as the equivalence
$(\phi_k,\psi_k,\eta_k,\epsln_k)$ was chosen to be an adjunction too (which
is always possible).  Identities work similarly.

Conversely, take a bicategory $B$ and construct a weak 2-category $A:
\Delnop{2}$\linebreak$\go \Set$ (its `2-nerve') as follows.  An element of
$A([j],[k])$ is a quadruple
\renewcommand{\arraystretch}{0.8}	
\[
((a_u)_{0\leq u\leq j},
(f_{uv}^z)_{\begin{array}{cc} \scriptstyle
	0\leq u < v\leq j\\ \scriptstyle
	0\leq z \leq k
	\end{array}},
(\alpha_{uv}^z)_{\begin{array}{cc} \scriptstyle
	0\leq u < v \leq j\\ \scriptstyle
	1\leq z \leq k
	\end{array}},
(\iota_{uvw}^z)_{\begin{array}{cc} \scriptstyle
	0\leq u < v < w \leq j\\ \scriptstyle
	0 \leq z \leq k
	\end{array}})
\]
\renewcommand{\arraystretch}{1}		
where
\begin{itemize}
\item $a_u$ is an object of $B$
\item $f_{uv}^z: a_u \go a_v$ is a 1-cell of $B$
\item $\alpha_{uv}^z: f_{uv}^{z-1} \go f_{uv}^z$ is a
2-cell of $B$
\item $\iota_{uvw}^z: f_{vw}^z \of f_{uv}^z 
\goiso f_{uw}^z$ is an invertible 2-cell of $B$ 
\end{itemize}
such that
\begin{itemize}
\item $\iota_{uvw}^z \,\of\, (\alpha_{vw}^z * \alpha_{uv}^z) =
\alpha_{uw}^z \,\of\, \iota_{uvw}^{z-1}$ 
whenever $0\leq u < v < w \leq j$, $1\leq z \leq k$
\item $\iota_{uwx}^z \,\of\, (1_{f_{wx}^z} * \iota_{uvw}^z)
\,\of\, (\textrm{associativity isomorphism})
=
\iota_{uvx}^z \,\of\, (\iota_{vwx}^z * 1_{f_{uv}^z})$
whenever $0\leq u < v < w < x \leq j$, $0\leq z \leq k$.
\end{itemize}
This defines the functor $A$ on objects of $\Deln{2}$; it is defined on maps
by a combination of inserting identities and forgetting data.

To get a rough picture of $A$, consider the analogous construction for
strict 2-categories, in which we insist that the isomorphisms
$\iota_{uvw}^z$ are actually equalities.  Then an element of
\marginpar{\centering\fbox{%
$
\begin{array}{c}
\gfstsu\gfoursu\gzersu\gfoursu\glstsu\\
j=2,\, k=3
\end{array}
$
}}
$A([j],[k])$ is a grid of $jk$ $2$-cells, of width $j$ and height $k$.  (When
$j=0$ this is just a single object of $B$, regardless of $k$.)  The
bicategorical version is a suitable weakening of this construction.

Finally, passing from a bicategory to a weak 2-category and back again gives
a bicategory isomorphic (in the category of bicategories and weak functors)
to the original one.  Passing from a weak 2-category to a bicategory and back
again gives a weak 2-category which is `equivalent' to the original one in a
sense which we do not have quite enough vocabulary to make precise here.

\clearpage

\defnheading{Ta}		\label{p:ta}

\concept{Simplicial Objects}

\paragraph{The Simplicial Category}

Let \Del\ be a skeleton of the category of nonempty finite totally ordered
sets: that is, \Del\ has objects $[k]=\{0,\ldots,k\}$ for $k\geq 0$, and maps
are order-preserving functions (with respect to the usual ordering $\leq$).

\paragraph{Some Objects and Morphisms}

Given $k\geq 0$, let $\iota_1, \ldots, \iota_k: [1] \go [k]$ denote the
`embeddings' of $[1]$ into $[k]$, defined by $\iota_j(0) = j-1$ and
$\iota_j(1) = j$.  Let $\sigma, \tau: [0] \go [1]$ be the maps in \Del\ with
respective values $0$ and $1$.  Given $p\geq 0$, write $0^p = ([0], \ldots,
[0]) \in \Deln{p}$ and $1^p = ([1], \ldots, [1]) \in \Deln{p}$.
Let $X: \Delnop{r} \go \Set$, $0 \leq p \leq r$, and $x, x' \in X(1^p,
0^{r-p})$.  Then $x, x'$ are \demph{parallel} if $p=0$ or if $p\geq 1$ and
$s(x) = s(x')$ and $t(x) = t(x')$; here $s$ and $t$ are the maps
\begin{diagram}
X(1^p,0^{r-p})	&
\pile{\rTo^{X(\id, \ldots, \id, X\sigma, \id, \ldots, \id)}\\
\rTo_{X(\id, \ldots, \id, X\tau, \id, \ldots, \id)}}	&
X(1^{p-1},0^{r-p+1}).\\
\end{diagram}

\paragraph{The Segal Maps}

Let $k\geq 0$.  Then the following diagram in \Del\ commutes:
\begin{diagram}[width=2em,height=2em]
	&	&	&	&	&	&	&[k]	&	&
	&	&	&	&	&	\\
	&	&	&	&	&	&\ruTo(7,2)<{\iota_1}	
						 \ruTo(3,2)<{\iota_2}
							&	&\ldots\luTo(7,2)>{\iota_k}	&
	&	&	&	&	&	\\
[1]	&	&	&	&[1]	&	&	&	&	&
\ldots	&	&	&	&	&[1]	\\
	&\luTo<\tau&	&\ruTo>\sigma&	&\luTo<\tau&	&\ruTo>\sigma&	&
	&	&\luTo<\tau&	&\ruTo>\sigma&	\\
	&	&[0]	&	&	&	&[0]	&	&	&
\cdots	&	&	&[0].	&	&	\\
\end{diagram}
Let $X: \Delop \go \cat{E}$ be a functor into a category
\cat{E} possessing finite limits, and write
$
X[1] \times_{X[0]} \cdots \times_{X[0]} X[1]
$
(with $k$ occurrences of $X[1]$) for the limit of the diagram  
\begin{diagram}[width=2em,height=2em]
X[1]	&	&	&	&X[1]	&	&	&	&	&
\ldots	&	&	&	&	&X[1]	\\
	&\rdTo<{X\tau}&	&\ldTo>{X\sigma}&	&\rdTo<{X\tau}&	&\ldTo>{X\sigma}&	&
	&	&\rdTo<{X\tau}&	&\ldTo>{X\sigma}&	\\
	&	&X[0]	&	&	&	&X[0]	&	&	&
\cdots	&	&	&X[0]	&	&	\\
\end{diagram}
(with, again, $k$ occurrences of $X[1]$) in \cat{E}.  Then by commutativity
of the first diagram, there is an induced map in \cat{E}---a \demph{Segal
map}---
\begin{equation}	\label{eq:Segal}
X[k] \go X[1] \times_{X[0]} \cdots \times_{X[0]} X[1].
\end{equation}

\paragraph{Nerves}

Call $X: \Delop \go \Set$ a \demph{nerve} if for each $k\geq
0$, the Segal map~\bref{eq:Segal} is a bijection.  The category of nerves and
natural transformations is equivalent to \Cat, where a nerve $X$ corresponds
to a category with object-set $X[0]$ and morphism-set $X[1]$.   
Let $QX$ be the set of isomorphism classes of objects of the category
corresponding to $X$, and let $\pi_X: X[0] \go QX$ be the quotient map.

\concept{Truncatability}

For each $r\geq 0$, we define what it means for $X: \Delnop{r} \go
\Set$ to be \demph{truncatable}, writing $\Trunc{r}$ for the
category of truncatable functors $\Delnop{r} \go \Set$ and natural
transformations between them.  We also define functors $\mr{ob}^{(r)}$,
$Q^{(r)}: \Trunc{r} \go \Set$ and a natural transformation
$\pi^{(r)}: \mr{ob}^{(r)} \go Q^{(r)}$.

The functor $\mr{ob}^{(r)}$ is given by $\mr{ob}^{(r)}X = X(0^r)$.  All
functors $\Delnop{0} \go \Set$ are truncatable, and $Q^{(0)}$ and
$\pi^{(0)}$ are identities.  Inductively, when $r\geq 1$, a functor $X:
\Delnop{r} \go \Set$ is truncatable if 
\begin{itemize}
\item for each $k\geq 0$, the functor $X([k],\dashbk): \Delnop{r-1} \go
\Set$ is truncatable 
\item the functor $\widehat{X}: \Delop \go \Set$ defined by $[k] \goesto
Q^{(r-1)}(X([k],\dashbk))$ is a nerve.
\end{itemize}
If $X$ is truncatable then we define $Q^{(r)}(X)=Q(\widehat{X})$ and
$\pi^{(r)}_X = \pi_{\widehat{X}} \,\of\, \pi^{(r-1)}_{X([0],\dashbk)}$.

\concept{Equivalence}

\paragraph{Internal Equivalence}

Let $0\leq p\leq r$, let $X: \Delnop{r} \go \Set$ be truncatable, and let
$x_1,x_2 \in X(1^p, 0^{r-p})$.  We call $x_1$ and $x_2$ \demph{equivalent},
and write $x_1\inteqv x_2$, if $x_1$ and $x_2$ are parallel and
$\pi^{(r-p)}_{X(1^p,\dashbk)}(x_1) = \pi^{(r-p)}_{X(1^p,\dashbk)}(x_2)$.

\paragraph{External Equivalence}

Let $r\geq 0$.  A natural transformation $\phi: X \go Y$ of truncatable
functors $X, Y: \Delnop{r} \go \Set$ is called an \demph{equivalence} if
\begin{itemize}
\item 
for each $y\in Y(0^r)$ there exists $x\in X(0^r)$ with $\phi_{0^r}(x)
\inteqv y$, and this $x$ is unique up to equivalence
\item
for all $0 \leq p\leq r-1$, parallel $x, x' \in X(1^p,0^{r-p})$, and $h\in
Y(1^{p+1},0^{r-p-1})$ satisfying
\[
s(h)=\phi_{(1^p,0^{r-p})}(x),	\diagspace 
t(h)=\phi_{(1^p,0^{r-p})}(x'),
\]
there is an element $g\in X(1^{p+1},0^{r-p-1})$, unique up to equivalence,
satisfying 
\marginpar{\centering
\fbox{%
\begin{diagram}[width=2em,height=1.3em]
x	&\rGet^g	&x'	\\
	&		&	\\
	&\dGoesto>\phi	&	\\
	&		&	\\
\phi(x)	&\rTo^h		&\phi(x')\\
\end{diagram}}}
\[
s(g)=x, 		\diagspace
t(g)=x', 		\diagspace
\phi_{(1^{p+1},0^{r-p-1})}(g) \inteqv h.
\]
\end{itemize}

\concept{The Definition}

Let $n\geq 0$.  A \demph{weak $n$-category} is a truncatable functor $A:
\Delnop{n} \go \Set$ such that for each $m\in \{0, \ldots, n-1\}$ and
$K=([k_1], \ldots, [k_m]) \in \Deln{m}$,
\begin{enumerate}
\item 	\label{part:defn:degen}
the functor
$
A(K, [0], \dashbk): \Delnop{n-m-1} \go \Set
$ 
is constant, and
\item 	\label{part:defn:main} 
for each $[k]\in\Del$, the Segal map
\begin{equation}	\label{eq:defining}
A(K, [k], \dashbk) \go 
A(K, [1],\dashbk) \times_{A(K, [0],\dashbk)} \cdots 
\times_{A(K, [0],\dashbk)} A(K, [1],\dashbk)
\end{equation}
is an equivalence.  (We are taking $\cat{E}=\ftrcat{\Delnop{n-m-1}}{\Set}$
and $X[j]=A(K,[j],\dashbk)$ in the definition of Segal map, and we can check
that both the domain and the codomain of~\bref{eq:defining} are truncatable.)
\end{enumerate}

\clearpage

\lowdimsheading{Ta}

\concept{$n=0$}

Parts~\bref{part:defn:degen} and~\bref{part:defn:main} of the definition are
vacuous, and truncatability is automatic, so a weak $0$-category is just a
functor $\Delnop{0}\go\Set$, that is, a set.

\concept{$n=1$}

Note that a functor $A: \Delop \go \Set$ is truncatable exactly when it is a
nerve; that a functor $X: \Delnop{0} \go \Set$ is merely a set, and two
elements of $X$ are equivalent just when they are equal; and that a map
$\phi: X \go Y$ of functors $X, Y: \Delnop{0} \go \Set$ is an equivalence
just when it is a bijection.  A weak $1$-category is a truncatable functor
$A: \Delop \go \Set$ satisfying~\bref{part:defn:degen}
and~\bref{part:defn:main}.  Part~\bref{part:defn:degen} is trivially true,
and both truncatability and~\bref{part:defn:main} say that $A$ is a nerve.
So a weak $1$-category is just a nerve---that is, essentially just a
category.

\concept{$n=2$}

First note that if $X: \Delop \go \Set$ is a nerve then two elements of
$X[0]$ are equivalent just when they are isomorphic (as objects of
the category corresponding to $X$), and two elements of $X[1]$ are equivalent
just when they are equal.  Note also that a map $\phi: X \go Y$ of nerves is
an equivalence if and only if (regarded as a functor between the
corresponding categories) it is full, faithful and essentially surjective on
objects---that is, an equivalence of categories.

A weak $2$-category is a truncatable functor $A: \Delnop{2} \go \Set$ such
that
\begin{enumerate}
\item the functor $A([0],\dashbk): \Delop \go \Set$ is constant
\item for each $k\geq 0$, the Segal map 
\[
A([k], \dashbk) \go 
A([1],\dashbk) \times_{A([0],\dashbk)} \cdots 
\times_{A([0],\dashbk)} A([1],\dashbk)
\]
is an equivalence, and for each $k_1, k\geq 0$, the Segal map
\[
A([k_1],[k]) \go 
A([k_1],[1]) \times_{A([k_1],[0])} \cdots 
\times_{A([k_1],[0])} A([k_1],[1])
\]
is a bijection.
\end{enumerate}
The second half of~\bref{part:defn:main} says that $A([k_1],\dashbk)$ is a
nerve for each $k_1$, so we can regard $A$ as a functor $A: \Delop \go \Cat$;
then the first half of~\bref{part:defn:main} says that the Segal
map~\bref{eq:Segal} (with $X=A$)
is an equivalence of categories.  Truncatability of $A$ says that the functor
$\Delop \go \Set$ given by
$
[k] \goesto \{
$%
isomorphism classes of objects of
$
A[k] \}
$
is a nerve, which follows anyway from the other conditions.  So a weak
$2$-category is a functor $A: \Delop \go \Cat$ such that
\begin{enumerate}
\item $A[0]$ is a discrete category (i.e.\ the only morphisms are the
identities) 
\item for each $k\geq 0$, the Segal functor 
$
A[k] \go A[1] \times_{A[0]} \cdots \times_{A[0]} A[1]
$
is an equivalence of categories. 
\end{enumerate}

It seems that a weak $2$-category is essentially just a bicategory.

First take a weak $2$-category $A: \Delop \go \Cat$, and let us construct a
bicategory $B$.  The object-set of $B$ is $A[0]$.  The two functors
$s,t: A[1] \go A[0]$ express the category $A[1]$ as a disjoint union
$\coprod_{a,b\in A[0]}B(a,b)$ of categories; the $1$-cells from $a$ to $b$
are the objects of $B(a,b)$, and the $2$-cells are the morphisms.

Vertical composition of $2$-cells in $B$ is composition in each
$B(a,b)$.  To define horizontal composition of $1$-~and $2$-cells,
first choose for each $k$ a pseudo-inverse
\[
A[1] \times_{A[0]} \cdots \times_{A[0]} A[1] \goby{\psi_k} A[k]
\]
to the Segal functor $\phi_k$, and choose natural isomorphisms $\eta_k: 1 \go
\psi_k\,\of\,\phi_k$, $\epsln_k: \phi_k\,\of\,\psi_k \go 1$.  Horizontal
composition is then given as
\[
A[1] \times_{A[0]} A[1] \goby{\psi_2} A[2] \goby{A\delta} A[1],
\]
where $\delta:[1] \go [2]$ is the injection whose image omits $1\in [2]$.
The associativity isomorphisms are built up from $\eta_k$'s and $\epsln_k$'s,
and the pentagon commutes just as long as the equivalence
$(\phi_k,\psi_k,\eta_k,\epsln_k)$ was chosen to be an adjunction too (which
is always possible).  Identities work similarly.

Conversely, take a bicategory $B$ and construct a weak 2-category $A:
\Delnop{2}$\linebreak$\go \Set$ (its `2-nerve') as follows.  An element of
$A([j],[k])$ is a quadruple
\renewcommand{\arraystretch}{0.8}	
\[
((a_u)_{0\leq u\leq j},
(f_{uv}^z)_{\begin{array}{cc} \scriptstyle
	0\leq u < v\leq j\\ \scriptstyle
	0\leq z \leq k
	\end{array}},
(\alpha_{uv}^z)_{\begin{array}{cc} \scriptstyle
	0\leq u < v \leq j\\ \scriptstyle
	1\leq z \leq k
	\end{array}},
(\iota_{uvw}^z)_{\begin{array}{cc} \scriptstyle
	0\leq u < v < w \leq j\\ \scriptstyle
	0 \leq z \leq k
	\end{array}})
\]
\renewcommand{\arraystretch}{1}		
where
\begin{itemize}
\item $a_u$ is an object of $B$
\item $f_{uv}^z: a_u \go a_v$ is a 1-cell of $B$
\item $\alpha_{uv}^z: f_{uv}^{z-1} \go f_{uv}^z$ is a
2-cell of $B$
\item $\iota_{uvw}^z: f_{vw}^z \of f_{uv}^z 
\goiso f_{uw}^z$ is an invertible 2-cell of $B$ 
\end{itemize}
such that
\begin{itemize}
\item $\iota_{uvw}^z \,\of\, (\alpha_{vw}^z * \alpha_{uv}^z) =
\alpha_{uw}^z \,\of\, \iota_{uvw}^{z-1}$ 
whenever $0\leq u < v < w \leq j$, $1\leq z \leq k$
\item $\iota_{uwx}^z \,\of\, (1_{f_{wx}^z} * \iota_{uvw}^z)
\,\of\, (\textrm{associativity isomorphism})
=
\iota_{uvx}^z \,\of\, (\iota_{vwx}^z * 1_{f_{uv}^z})$
whenever $0\leq u < v < w < x \leq j$, $0\leq z \leq k$.
\end{itemize}
This defines the functor $A$ on objects of $\Deln{2}$; it is defined on maps
by a combination of inserting identities and forgetting data.

To get a rough picture of $A$, consider the analogous construction for
strict 2-categories, in which we insist that the isomorphisms
$\iota_{uvw}^z$ are actually equalities.  Then an element of
\marginpar{\centering\fbox{%
$
\begin{array}{c}
\gfstsu\gfoursu\gzersu\gfoursu\glstsu\\
j=2,\, k=3
\end{array}
$
}}
$A([j],[k])$ is a grid of $jk$ $2$-cells, of width $j$ and height $k$.  (When
$j=0$ this is just a single object of $B$, regardless of $k$.)  The
bicategorical version is a suitable weakening of this construction.

Finally, it appears that passing from a bicategory to a weak 2-category and
back again gives a bicategory isomorphic (by weak functors) to the original
one, and that passing from a weak 2-category to a bicategory and back again
gives a weak 2-category equivalent to the original one.

\clearpage

\defnheading{J}		\label{p:j}

\concept{Disks}

\paragraph{Disks}

A \demph{disk} $D$ is a diagram of sets and functions
\vspace*{-3ex}
\[
\cdots\diagspace
D_m \bundleint{p_m}{u_m}{v_m} D_{m-1} 
\diagspace\cdots\diagspace
\bundleint{p_2}{u_2}{v_2} D_1
\bundleint{p_1}{u_1}{v_1} D_0 = 1
\]
equipped with a total order on the fibre $p_m^{-1}(d)$ for each $m\geq 1$
and $d\in D_{m-1}$, such that for each $m\geq 1$ and $d\in D_{m-1}$,
\begin{itemize}
\item $u_m(d)$ and $v_m(d)$ are respectively the least and greatest elements
of $p_m^{-1}(d)$
\item $u_m(d)=v_m(d) \iff d \in \textrm{image}(u_{m-1}) \cup
\textrm{image}(v_{m-1})$.
\end{itemize}
When $m=1$, the second condition is to be interpreted as saying that $u_1\neq
v_1$ (or equivalently, that $D_1$ has at least two elements).

A \demph{map $D\goby{\psi} D'$ of disks} is a family of functions $(D_m
\goby{\psi_m} D'_m)_{m\geq 0}$ commuting with the $p$'s, $u$'s and $v$'s and
preserving the order in each fibre.  (The last condition means that if $d\in
D_{m-1}$ and $b,c\in p_m^{-1}(d)$ with $b\leq c$, then $\psi_m(b)\leq
\psi_m(c) \in p'^{-1}_m(\psi_{m-1}(d))$.)  Call $\psi$ a \demph{surjection}
if each $\psi_m$ is a surjection.

\paragraph{Interiors, Volume, Dimension}

Let $D$ be a disk.  For $m\geq 1$, define
\[
\iota D_m = D_m \backslash (\textrm{image}(u_m) \cup 
\textrm{image}(v_m)),
\]
(the \demph{interior} of $D_m$), and define $\iota D_0 = D_0$.  If the set
$\coprod_{m\geq 1} \iota D_m$ is finite then we call $D$ \demph{finite} and
define the \demph{volume} $|D|$ of $D$ to be its cardinality.  In this case
we may also define the \demph{dimension} of $D$ to be the largest $m\geq 0$
for which $\iota D_m \neq \emptyset$.

\paragraph{Finite Disks}

Write \scat{D} for a skeleton of the category of finite disks and maps
between them.  In other words, take the category of all finite disks and
choose one object in each isomorphism class; the objects of $\scat{D}$ are
all these chosen objects, and the morphisms in $\scat{D}$ are all disk maps
between them.  Thus \scat{D} is equivalent to the category of finite disks
and no two distinct objects of \scat{D} are isomorphic.

\concept{Faces and Horns, Cofaces and Cohorns}

\paragraph{Cofaces}

Let $D \in \scat{D}$.  A \demph{(covolume $1$) coface} of $D$ is a
surjection $D \goby{\phi} E$ in \scat{D} where $|E| = |D|-1$.  We call $\phi$
an \demph{inner coface} of $D$ if $\phi_m(\iota D_m) \sub \iota E_m$ for all
$m\geq 0$.

\paragraph{Cohorns}

For each $D \in \scat{D}$ and coface $D \goby{\phi} E$ of \scat{D}, define
the \demph{cohorn} 
\[
\Lambda^D_\phi: \scat{D} \go \Set
\]
by
\[
\Lambda^D_\phi(C) = \{ \psi \in \scat{D}(D,C) \such
\psi \textrm{ factors through some coface of }D
\textrm{ other than }\phi \}.
\]
That is, a map $\psi: D \go C$ is a member of $\Lambda^D_\phi(C)$ if and
only if there is a coface $(D \goby{\phi'}E') \neq (D \goby{\phi}E)$ of $D$
and a map $\chi: E' \go C$ such that
\begin{diagram}[height=2em]
D	&\rTo^{\phi'}	&E'			\\
	&\rdTo<\psi	&\dTo>{\chi}		\\
	&		&C			\\
\end{diagram}
commutes.  There is an inclusion $\Lambda^D_\phi(C) \rIncl \scat{D}(D,C)$ for
each $C$, and $\Lambda^D_\phi$ is thus a subfunctor of $\scat{D}(D,\dashbk)$.
Write $i^D_\phi: \Lambda^D_\phi \rIncl \scat{D}(D, \dashbk)$ for the
inclusion.

\paragraph{Fillers}

Let $A: \scat{D} \go \Set$, let $D\in\scat{D}$, and let $\phi$ be a coface of
$D$.  A \demph{$(D,\phi)$-cohorn in $A$} is a natural transformation
$h:\Lambda^D_\phi \go A$; if $\phi$ is an inner coface then $h$ is an
\demph{inner cohorn in $A$}.

A \demph{filler} for a $(D,\phi)$-cohorn $h$ in $A$ is a natural
transformation $\ovln{h}: \scat{D}(D,\dashbk)$ $\go A$ making the following
diagram commute:
\begin{diagram}[height=2em]
\Lambda^D_\phi	&\rIncl^{i^D_\phi}	&\scat{D}(D,\dashbk)		\\
		&\rdTo<h		&\dTo>{\ovln{h}}		\\
		&			&A.				\\
\end{diagram}

\concept{The Definition}

\paragraph{Weak $\omega$-Categories}

A \demph{weak $\omega$-category} is a functor $A: \scat{D} \go \Set$ such
that there exists a filler for every inner cohorn in $A$.

\paragraph{Weak $n$-Categories}

Let $n\geq 0$.  A functor $A: \scat{D} \go \Set$ is \demph{$n$-dimensional}
if, whenever $\psi: D \go E$ is a map in \scat{D} such that
\begin{itemize}
\item $D$ has dimension $n$ 
\item $\psi_m$ is a bijection for every $m\leq n$,
\end{itemize}
then $A(\psi)$ is a bijection.  A \demph{weak $n$-category} is an
$n$-dimensional weak $\omega$-category.

\clearpage

\lowdimsheading{J}

Let $n\geq 0$.  An \demph{$n$-disk} is defined in the same way as a disk,
except that $D_m$, $p_m$, $u_m$ and $v_m$ are now only defined for $m\leq n$:
so an $n$-disk is essentially the same thing as a disk of dimension $\leq n$.
Write $\scat{D}_n$ for a skeleton of the category of finite $n$-disks.  An
$n$-dimensional functor $\scat{D} \go \Set$ is determined by its effect on
disks of dimension $\leq n$, and conversely any functor $\scat{D}_n \go \Set$
extends uniquely to become an $n$-dimensional functor $\scat{D} \go \Set$.
So the category of $n$-dimensional functors $\scat{D}\go \Set$ is equivalent
to $\ftrcat{\scat{D}_n}{\Set}$.

Take an $n$-dimensional functor $A: \scat{D} \go \Set$ and its restriction
$\twid{A}: \scat{D}_n$ $\go \Set$.  Then there is automatically a unique
filler for every cohorn of dimension $\geq n+2$ in $A$ (that is, cohorn
$\Lambda^D_\phi \go A$ where $D$ has dimension $\geq n+2$).  Moreover, there
exists a filler for every inner cohorn of dimension $n+1$ in $A$ if and only
if there is at most one filler for every inner cohorn of dimension $n$ in
$\twid{A}$.  (We do not prove this, but the idea of the method is in the last
sentence of `$n=2$'.)  So: a weak $n$-category is a functor $\scat{D}_n \go
\Set$ such that every inner cohorn has a filler, unique when the cohorn is of
dimension $n$.

\concept{$n=0$}

$\scat{D}_0$ is the terminal category \One.  The unique $0$-disk has no
cofaces, so a weak $0$-category is merely a functor $\One \go \Set$, that is,
a set.

\concept{$n=1$}

An \demph{interval} is a totally ordered set with a least and a greatest
element, and is called \demph{strict} if these elements are distinct.
$\scat{D}_1$ is (a skeleton of) the category of finite strict intervals, so
we can take its objects to be the intervals $\intvl{k} = \{0, \ldots, k+1\} $
for $k\geq 0$ and its morphisms to be the interval maps.

The cofaces of \intvl{k} are the surjections $\intvl{k}\go\intvl{k-1}$
(assuming $k\geq 1$; if $k=0$ then there are none).  They are $\phi_0,
\ldots, \phi_k$, where $\phi_i$ identifies $i$ and $i+1$; of these, $\phi_1,
\ldots, \phi_{k-1}$ are inner.  The cohorn $\Lambda^{\intvl{k}}_{\phi_i}:
\scat{D}_1 \go \Set$ sends \intvl{l} to
$
\{ \psi: \intvl{k}\go\intvl{l} \such
\psi \textrm{ factors through }\phi_{i'} \textrm{ for some } 
i'\in\{0, \ldots, i-1, i+1, \ldots, k\} \}.
$

Now, let $\Del$ be a skeleton of the category of nonempty finite totally
ordered sets, with objects $[k] = \{0, \ldots, k\}$ ($k\geq 0$).  Then
$\scat{D}_1 \iso \Del^\op$, with $\intvl{k}$ corresponding to $[k]$, the
cofaces $\phi_i: \intvl{k}\go\intvl{k-1}$ to the usual face maps $[k-1] \go
[k]$, and the inner cofaces to the inner faces (i.e.\ all but the first and
last).  Trivially, cohorns $\Lambda^{\intvl{k}}_{\phi_i}$ correspond to horns
in the standard sense, and fillers to fillers.  Hence a weak $1$-category is
a functor $A: \Delop \go \Set$ in which every inner horn has a unique
filler---exactly the condition that $A$ is the nerve of a category.  So a
weak $1$-category is just a category.

\concept{$n=2$}

\begin{figure}
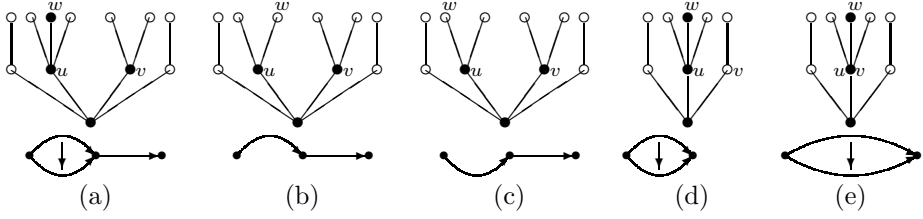

\[
\begin{array}{ccccc}
\begin{tree}
\enode		&\enode		&\lanode{w}	&\enode		&	&
\enode		&		&\enode		&\enode		\\
\dn		&		&\rt{1}\dn\lt{1}&		&	&
		&\rt{1}\lt{1}	&		&\dn		\\
\enode		&		&\lnode{u}	&		&	&
		&\lnode{v}	&		&\enode		\\
		&\rt{4}		&		&\rt{2}		&	&
\lt{2}		&		&\lt{4}		&		\\
		&		&		&		&\node	&
		&		&		&		\\
\end{tree}
\ &
\begin{tree}
\enode		&\enode		&		&\laenode{w}	&	&
\enode		&		&\enode		&\enode		\\
\dn		&		&\rt{1}\lt{1}	&		&	&
		&\rt{1}\lt{1}	&		&\dn		\\
\enode		&		&\lnode{u}	&		&	&
		&\lnode{v}	&		&\enode		\\
		&\rt{4}		&		&\rt{2}		&	&
\lt{2}		&		&\lt{4}		&		\\
		&		&		&		&\node	&
		&		&		&		\\
\end{tree}
\ &
\begin{tree}
\enode		&\laenode{w}	&		&\enode		&	&
\enode		&		&\enode		&\enode		\\
\dn		&		&\rt{1}\lt{1}	&		&	&
		&\rt{1}\lt{1}	&		&\dn		\\
\enode		&		&\lnode{u}	&		&	&
		&\lnode{v}	&		&\enode		\\
		&\rt{4}		&		&\rt{2}		&	&
\lt{2}		&		&\lt{4}		&		\\
		&		&		&		&\node	&
		&		&		&		\\
\end{tree}
\ &
\begin{tree}
\enode		&\enode		&\lanode{w}	&\enode		&\enode	\\
\dn		&		&\rt{1}\dn\lt{1}&		&\dn	\\
\enode		&		&\lnode{u}	&		&\lenode{v}\\
		&\rt{2}		&\dn		&\lt{2}		&	\\
		&		&\node		&		&	\\
\end{tree}
\ &
\begin{tree}
\enode		&\enode		&\lanode{w}	&\enode		&\enode	\\
\dn		&		&\rt{1}\dn\lt{1}&		&\dn	\\
\enode		&		&\lnode{\!\!\!\!\!\!\!u\ v}&	&\enode	\\
		&\rt{2}		&\dn		&\lt{2}		&	\\
		&		&\node		&		&	\\
\end{tree}
\\
\gzersu\gtwosu\gzersu\gonesu\gzersu	
&	
\gzersu\gunhappysu\gzersu\gonesu\gzersu
&	
\gzersu\ghappysu\gzersu\gonesu\gzersu
&
\gzersu\gtwosu\gzersu\ghole
&
\gzersu\gtwowidesu\gzersu
\\
\textrm{(a)}	&\textrm{(b)}	&\textrm{(c)}	&
\textrm{(d)}	&\textrm{(e)}	
\end{array}
\]
\caption{The duality.  In the upper row, $\node$ denotes an interior element
and $\enode$ an endpoint of a fibre, and the labels $u, v, w$ show what the
coface maps `$\phi$' do}
\label{fig:disks}
\end{figure}

Again we use a duality.  Given natural numbers $l_1, \ldots, l_k$, let
$T_{l_1, \ldots, l_k}$ be the strict $2$-category generated by objects $x_0,
\ldots, x_k$, $1$-cells $p_i^j: x_{i-1} \go x_i$ ($1\leq i\leq k$, $0\leq j
\leq l_i$), and $2$-cells $\xi_i^j: p_i^{j-1} \go p_i^j$ ($1\leq i\leq k$,
$1\leq j \leq l_i$).  (E.g.\ the lower half of Fig.~\ref{fig:disks}(a) shows
$T_{1,0}$.)  Let $\Delta_2$ be the category whose objects are sequences
$(l_1, \ldots, l_k)$ with $k,l_i\geq 0$, and whose maps $(l_1, \ldots, l_k)
\go (l'_1, \ldots, l'_{k'})$ are the strict $2$-functors $T_{l_1, \ldots,
l_k} \go T_{l'_1, \ldots, l'_{k'}}$.  Then $\scat{D}_2 \iso \Delta_2^\op$.
On objects, this says that a finite 2-disk is just a finite sequence of
numbers, e.g.\ $(1,0)$ in Fig.~\ref{fig:disks}(a).

Any bicategory $B$ has a `nerve' $A: \Delta_2^\op \go \Set$, where $A(l_1,
\ldots, l_k) = \{$weak functors $T_{l_1, \ldots, l_k} \go B$ strictly
preserving identities$\}$.  We can recover $B$ from $A$, so weak
$2$-categories are the same as bicategories just as long as the definition
gives the right conditions on functors $\scat{D}_2 \iso \Delta_2^\op \go
\Set$.  I do not have a full proof that this is so, hence there are gaps in
what follows.

Defining \demph{faces} of an object of $\Delta_2$ as cofaces of the
corresponding object of $\scat{D}_2$, and similarly \demph{horns}, a weak
$2$-category is a functor $\Delta_2^\op \go \Set$ in which every 1-
(respectively, 2-) dimensional horn has a filler (respectively, unique
filler).  Faces are certain subcategories: e.g.\ Fig.~\ref{fig:disks}(b)--(e)
shows the 4 cofaces of a disk and correspondingly the 4 faces of $T_{1,0}$,
of which only (e) is inner.

For the converse of the nerve construction, we take a weak $2$-category $A$
and form a bicategory $B$.  Its graph $B_2 \parpairu B_1 \parpairu B_0$ is
the image under $A$ of the diagram $\gfstsu\gtwosu\glstsu \pile{\lIncl \\
\lIncl} \gfstsu\gonesu\glstsu \pile{\lIncl \\ \lIncl} \gzeros{}$ in
$\Delta_2$.  A diagram $\gfsts{a}\gones{f}\gblws{b}\gones{g}\glsts{c}$ in $B$
is a horn in $A$ for the unique inner face of $T_{0,0} =
\gfstsu\gonesu\gzersu\gonesu\glstsu$; choose a filler $K_{f,g}$ and write
$g\of f$ for its third face.  This gives 1-cell composition; vertical 2-cell
composition works similarly but without choice.  Next, a diagram
$\gfsts{a}\gtwos{f}{f'}{\alpha}\gfbws{b}\gones{g}\glsts{c}$, with $K_{f,g}$
and $K_{f',g}$, forms a horn for the unique inner face of $T_{1,0}$
(Fig.~\ref{fig:disks}), so has a unique filler $K_{\alpha,g}$; write $g\of
\alpha : g\of f \go g\of f'$ for face~(e) of $K_{\alpha,g}$.  Horizontal
2-cell composition is defined via this construction, its dual, and vertical
composition.  Next, $\gfsts{a} \gones{f} \gblws{b} \gones{g} \gblws{c}
\gones{h} \glsts{d}$, together with $K_{f,g}, K_{g,h}, K_{g\sof f,h}$, gives
an inner horn for $T_{0,0,0}$.  There's a (unique?)  filler, whose final face
$L_{f,g,h}$ is itself a filler of $\gfsts{a} \gones{f} \gblws{b} \gones{h\sof
g} \glsts{d}$ with third face $h\of (g\of f)$.  Considering
$\gfsts{a}\gtwos{f}{f}{1} \gfbws{b} \gones{h\sof g}\glsts{d}$ with
$K_{f,h\sof g}$ and $L_{f,g,h}$ gives an invertible 2-cell $(h\of g)\of f \go
h\of (g\of f)$.

\defnheading{St}		\label{p:st}

\concept{Simplicial Sets}

\paragraph{The Simplicial Category}

Let \Del\ be a skeleton of the category of nonempty finite totally ordered
sets: that is, \Del\ has objects $[m]=\{0,\ldots,m\}$ for $m\geq 0$, and maps
are order-preserving functions (with respect to the usual ordering $\leq$).
A \demph{simplicial set} is a functor $\Delop \go \Set$.

\paragraph{Maps in \Del}

Let $m\geq 1$: then there are injections $\delta_0, \ldots, \delta_m: [m-1]
\go [m]$ in \Del, determined by saying that the image of $\delta_i$ is $[m]
\without \{ i \}$.  

Let $A: \Delop \go \Set$ and $m\geq 0$.  An element $a\in A[m]$ is called
\demph{degenerate} if there exist a natural number $m'<m$, a surjection
$\sigma: [m] \go [m']$, and an element $a' \in A[m']$ such that $a =
(A\sigma)a'$.

\paragraph{Horns}

Given $0\leq k\leq m$, we define the \demph{horn} $\Lambda^k_m: \Delop \go
\Set$ by
\[
\Lambda^k_m[m'] = \{ \psi \in \Del([m'],[m]) \such \textrm{image}(\psi)
\not\supseteq [m]\without\{k\} \}.
\]
That is, $\Lambda^k_m[m']$ is the set of all maps $\psi: [m'] \go [m]$ in
\Del\ except for the surjections and the maps with image $ \{ 0, \ldots, k-1,
k+1, \ldots, m \}$.  So for each $m'$ we have an inclusion $\Lambda^k_m[m']
\rIncl \Del([m'],[m])$, and $\Lambda^k_m$ is thus a subfunctor of
$\Del(\dashbk,[m])$.  Write $i^k_m: \Lambda^k_m \rIncl \Del(\dashbk,[m])$ for
the inclusion.

Let $A$ be a simplicial set.  A \demph{horn in $A$} is a natural
transformation $h: \Lambda^k_m \go A$, for some $0\leq k\leq m$.  A
\demph{filler} for the horn $h$ is a natural transformation $\ovln{h}:
\Del(\dashbk,[m]) \go A$ making the following diagram commute:
\begin{diagram}[height=2em]
\Lambda^k_m	&\rIncl^{i^k_m}	&\Del(\dashbk,[m])	\\
		&\rdTo<h	&\dTo>{\ovln{h}}	\\
		&		&A.			\\
\end{diagram}

\concept{Orientation}

\paragraph{Alternating Sets} 

A set of natural numbers is \demph{alternating} if its elements,
when written in ascending order, alternate in parity.

Let $0\leq k\leq m$, and write $k^\pm = \{k-1,k,k+1\} \cap [m]$.  A subset $S
\sub [m]$ is \demph{$k$-alternating} if
\begin{itemize}
\item $k^\pm \sub S$
\item the set $k^\pm \cup ([m] \without S)$ is alternating.
\end{itemize}

\paragraph{Admissible Horns}

A \demph{simplicial set with hollowness} is a simplicial set $A$ together
with a subset $H_m\sub A[m]$ for each $m\geq 1$, whose elements are called
the \demph{hollow} elements of $A[m]$ (and may also be thought of as `thin'
or `universal').

Let $(A,H)$ be a simplicial set with hollowness, and $0\leq k\leq m$.  A horn
$h: \Lambda^k_m \go A$ is \demph{admissible} if for every $m'\geq 1$ and
$\psi\in\Lambda^k_m[m']$,
\[
\textrm{image}(\psi) \textrm{ is a } k\textrm{-alternating subset of }[m]
\ \implies\ 
h_{[m']}(\psi) \textrm{ is hollow.}	
\]

\concept{The Definition}

\paragraph{Weak $\omega$-Categories}

A \demph{weak $\omega$-category} is a simplicial set with hollowness $(A,H)$
such that
\begin{enumerate}
\item  	\label{part:degen} 
for $m\geq 1$, $H_m \supseteq \{ \textrm{degenerate elements of }
A[m] \}$
\item 	\label{part:filler}
for $m\geq 1$ and $0\leq k\leq m$, every admissible horn $h:
\Lambda^k_m \go A$ has a filler $\ovln{h}$ satisfying $\ovln{h}_{[m]}(1_{[m]})
\in H_m$ (`every admissible horn has a hollow filler')
\item  	\label{part:comp}
for $m\geq 2$ and $0\leq k\leq m$, if $a\in H_m$ has the property that
$(A\delta_i)a \in H_{m-1}$ for each $i \in [m]\without \{k\}$ then also
$(A\delta_k)a \in H_{m-1}$. 
\end{enumerate}

\paragraph{Weak $n$-Categories}

Let $n\geq 0$.  A \demph{weak $n$-category} is a weak $\omega$-category
$(A,H)$ such that 
\renewcommand{\theenumi}{\roman{enumi}$'$}
\begin{enumerate}
\item 	\label{part:high-dims}
for $m > n$, $H_m = A[m]$
\item 	\label{part:unique}
in condition~\bref{part:filler} above, when $m > n$ there is a
\emph{unique} filler $\ovln{h}$ for $h$ (which necessarily satisfies
$\ovln{h}_{[m]}(1_{[m]}) \in H_m$). 
\end{enumerate}
\renewcommand{\theenumi}{\roman{enumi}}

\clearpage

\lowdimsheading{St}

Let $m\geq 0$ and let $S$ be a nonempty subset of $[m]$; then in \Del\ there
is a unique injection $\phi$ into $[m]$ with image $S$.  Given a simplicial
set $A$ and an element $a\in A[m]$, the \demph{$S$-face} of $a$ is the
element $(A\phi)a$ of $A[l]$, where $l+1$ is the cardinality of $S$.
Similarly, the \demph{$S$-face} of a horn $h: \Lambda^k_m \go A$ is
$h_{[l]}(\phi) \in A[l]$ (which makes sense as long as $S \not\supseteq
[m]\without \{k\}$).

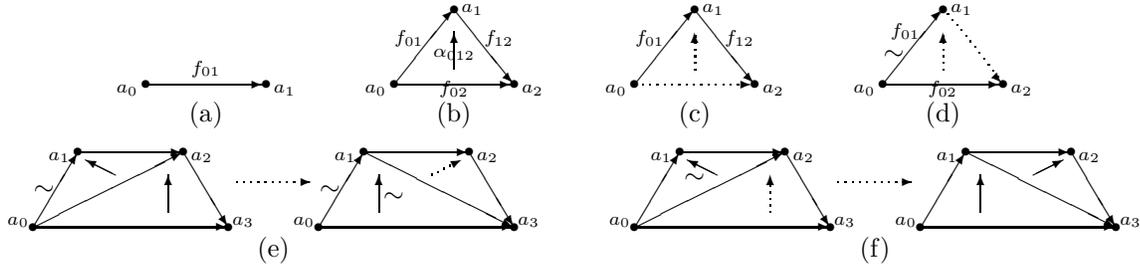
\begin{figure}[t]
\setlength{\unitlength}{1mm}
\begin{picture}(122,35) 
\cell{12}{20}{c}{\textrm{(a)}}
\put(0,23){%
\begin{picture}(24,12)
\cell{4}{1}{c}{\zmark}
\cell{20}{1}{c}{\zmark}
\put(4,1){\vector(1,0){15.5}}
\cell{3}{1}{tr}{\scriptstyle a_0}
\cell{21}{1}{tl}{\scriptstyle a_1}
\cell{12}{2}{b}{\scriptstyle f_{01}}
\end{picture}
}
\cell{45}{20}{c}{\textrm{(b)}}
\put(33,23){%
\begin{picture}(24,12)
\cell{4}{1}{c}{\zmark}
\cell{20}{1}{c}{\zmark}
\cell{12}{11}{c}{\zmark}
\put(4,1){\vector(1,0){15.5}}
\put(4,1){\line(4,5){7.8}}
\put(11.8,10.8){\vector(1,1){0}}
\put(12,11){\line(4,-5){7.8}}
\put(19.8,1.2){\vector(1,-1){0}}
\put(12,3){\vector(0,1){5}}
\cell{3}{1}{tr}{\scriptstyle a_0}
\cell{13}{11}{l}{\scriptstyle a_1}
\cell{21}{1}{tl}{\scriptstyle a_2}
\cell{8}{6}{br}{\scriptstyle f_{01}}
\cell{16}{6}{bl}{\scriptstyle f_{12}}
\cell{12}{1.5}{t}{\scriptstyle f_{02}}
\cell{12}{5}{c}{\scriptstyle \alpha_{012}}
\end{picture}
}
\cell{77}{20}{c}{\textrm{(c)}}
\put(65,23){%
\begin{picture}(24,12)
\cell{4}{1}{c}{\zmark}
\cell{20}{1}{c}{\zmark}
\cell{12}{11}{c}{\zmark}
\qbezier[16](4,1)(11.8,1)(19.6,1)
\put(19.5,1){\vector(1,0){0}}
\put(4,1){\line(4,5){7.8}}
\put(11.8,10.8){\vector(1,1){0}}
\put(12,11){\line(4,-5){7.8}}
\put(19.8,1.2){\vector(1,-1){0}}
\qbezier[5](12,3)(12,5.2)(12,7.4)
\put(12,8){\vector(0,1){0}}
\cell{3}{1}{tr}{\scriptstyle a_0}
\cell{13}{11}{l}{\scriptstyle a_1}
\cell{21}{1}{tl}{\scriptstyle a_2}
\cell{8}{6}{br}{\scriptstyle f_{01}}
\cell{16}{6}{bl}{\scriptstyle f_{12}}
\end{picture}
}
\cell{110}{20}{c}{\textrm{(d)}}
\put(98,23){%
\begin{picture}(24,12)
\cell{4}{1}{c}{\zmark}
\cell{20}{1}{c}{\zmark}
\cell{12}{11}{c}{\zmark}
\put(4,1){\vector(1,0){15.5}}
\put(4,1){\line(4,5){7.8}}
\put(11.8,10.8){\vector(1,1){0}}
\qbezier[13](12,11)(15.9,6.6)(19.8,1.2)
\put(19.8,1.2){\vector(1,-1){0}}
\qbezier[5](12,3)(12,5.2)(12,7.4)
\put(12,8){\vector(0,1){0}}
\cell{3}{1}{tr}{\scriptstyle a_0}
\cell{13}{11}{l}{\scriptstyle a_1}
\cell{21}{1}{tl}{\scriptstyle a_2}
\cell{9}{7}{br}{\scriptstyle f_{01}}
\cell{12}{1.5}{t}{\scriptstyle f_{02}}
\cell{5.5}{5}{c}{\sim}
\end{picture}
}
\cell{21}{2}{c}{\textrm{(e)}}
\put(-11,5){%
\begin{picture}(64,10)
\cell{0}{0}{c}{\zmark}
\cell{26}{0}{c}{\zmark}
\cell{6}{10}{c}{\zmark}
\cell{20}{10}{c}{\zmark}
\put(0,0){\vector(1,0){25.8}}
\put(0,0){\line(3,5){5.8}}
\put(5.9,9.8){\vector(2,3){0}}
\put(6,10){\vector(1,0){13.8}}
\put(20,10){\line(3,-5){5.8}}
\put(25.9,0.2){\vector(2,-3){0}}
\put(0,0){\vector(2,1){19.8}}
\put(18,2){\vector(0,1){5}}
\put(11,7){\vector(-2,1){4}}
\cell{1.5}{5}{c}{\sim}
\cell{-0.5}{0}{br}{\scriptstyle a_0}
\cell{5}{10}{tr}{\scriptstyle a_1}
\cell{21}{10}{tl}{\scriptstyle a_2}
\cell{26.5}{0}{bl}{\scriptstyle a_3}
\qbezier[10](27,6)(31.8,6)(36.6,6)
\put(37,6){\vector(1,0){0}}
\cell{38}{0}{c}{\zmark}
\cell{64}{0}{c}{\zmark}
\cell{44}{10}{c}{\zmark}
\cell{58}{10}{c}{\zmark}
\put(38,0){\vector(1,0){25.8}}
\put(38,0){\line(3,5){5.8}}
\put(43.9,9.8){\vector(2,3){0}}
\put(44,10){\vector(1,0){13.8}}
\put(58,10){\line(3,-5){5.8}}
\put(63.9,0.2){\vector(2,-3){0}}
\put(44,10){\vector(2,-1){19.8}}
\put(46,2){\vector(0,1){5}}
\put(57,9){\vector(2,1){0}}
\qbezier[4](53,7)(54.6,7.8)(56.2,8.6)
\cell{39.5}{5}{c}{\sim}
\cell{48}{4}{c}{\sim}
\cell{37.5}{0}{br}{\scriptstyle a_0}
\cell{43}{10}{tr}{\scriptstyle a_1}
\cell{59}{10}{tl}{\scriptstyle a_2}
\cell{64.5}{0}{bl}{\scriptstyle a_3}
\end{picture}
}
\cell{101}{2}{c}{\textrm{(f)}}
\put(69,5){%
\begin{picture}(64,12)
\cell{0}{0}{c}{\zmark}
\cell{26}{0}{c}{\zmark}
\cell{6}{10}{c}{\zmark}
\cell{20}{10}{c}{\zmark}
\put(0,0){\vector(1,0){25.8}}
\put(0,0){\line(3,5){5.8}}
\put(5.9,9.8){\vector(2,3){0}}
\put(6,10){\vector(1,0){13.8}}
\put(20,10){\line(3,-5){5.8}}
\put(25.9,0.2){\vector(2,-3){0}}
\put(0,0){\vector(2,1){19.8}}
\put(18,7){\vector(0,1){0}}
\qbezier[5](18,2)(18,4.2)(18,6.4)
\put(11,7){\vector(-2,1){4}}
\cell{8}{6.5}{c}{\sim}
\cell{-0.5}{0}{br}{\scriptstyle a_0}
\cell{5}{10}{tr}{\scriptstyle a_1}
\cell{21}{10}{tl}{\scriptstyle a_2}
\cell{26.5}{0}{bl}{\scriptstyle a_3}
\qbezier[10](27,6)(31.8,6)(36.6,6)
\put(37,6){\vector(1,0){0}}
\cell{38}{0}{c}{\zmark}
\cell{64}{0}{c}{\zmark}
\cell{44}{10}{c}{\zmark}
\cell{58}{10}{c}{\zmark}
\put(38,0){\vector(1,0){25.8}}
\put(38,0){\line(3,5){5.8}}
\put(43.9,9.8){\vector(2,3){0}}
\put(44,10){\vector(1,0){13.8}}
\put(58,10){\line(3,-5){5.8}}
\put(63.9,0.2){\vector(2,-3){0}}
\put(44,10){\vector(2,-1){19.8}}
\put(46,2){\vector(0,1){5}}
\put(53,7){\vector(2,1){4}}
\cell{37.5}{0}{br}{\scriptstyle a_0}
\cell{43}{10}{tr}{\scriptstyle a_1}
\cell{59}{10}{tl}{\scriptstyle a_2}
\cell{64.5}{0}{bl}{\scriptstyle a_3}
\end{picture}
}
\end{picture}
\caption{(a)~Element of $A[1]$, with, for instance, the $\{0\}$-face
labelled $a_0$; 
(b)~element of $A[2]$; 
(c)~(admissible) horn $\Lambda^1_2 \protect\go A$;
(d)~admissible horn $\Lambda^0_2 \protect\go A$, with $\diso$ indicating a
hollow face; 
(e)~admissible horn $\Lambda^0_3 \protect\go A$, with labels $f_{ij},
\alpha_{ijk}$ omitted; 
(f)~as (e), but for $\Lambda^1_3$.}
\label{fig:simplices}
\end{figure}

To compare weak $1$-($2$-)categories with (bi)categories, we
need to interpret elements $a\in A[m]$ as arrows pointing in some direction.
Our convention is: if $S$ is an $m$-element subset of the
$(m+1)$-element set $[m]$ and the missing element is odd, then we regard the
$S$-face of $a$ as a source; if even, a target.
See Fig.~\ref{fig:simplices}.

Suppose $(A,H)$ is a simplicial set with hollowness
satisfying~\bref{part:degen}, and let $h: \Lambda^k_m \go A$ be a horn
satisfying the defining condition for admissibility for just the
\emph{injective} $\psi\in\Lambda^k_m[m']$.  Then, in fact, $h$ is admissible.
So $h$ is admissible if and only if: \emph{for every $k$-alternating subset
$S$ of $[m]$, the $S$-face of $h$ is hollow}.  Table~\ref{table:alt} shows the
$k$-alternating subsets of $[m]$ in the cases we need.
\begin{table}[b]
\centering
\begin{tabular}{c|c}
$k$		&$k$-alternating subsets of $[m]$ of cardinality $\leq 3$\\
\hline
$0$		&$\{0,1\}, \{0,1,m\}$					\\
$1, \ldots, m-1$&$\{k-1,k,k+1\}$					\\
$m$		&$\{m-1,m\}, \{0,m-1,m\}$				
\end{tabular}
\caption{$k$-alternating subsets of $[m]$ of cardinality $\leq 3$, for $m\geq
1$} 
\label{table:alt}
\end{table}

\concept{$n=0$}

A weak $0$-category is a simplicial set $A$ in which every horn has a unique
filler---including those of shape $\Lambda^k_1$.  It follows that the functor
$A: \Delop \go \Set$ is constant, so a weak $0$-category is just a set.

\concept{$n=1$}

A horn of shape $\Lambda_m^k$ is called \emph{inner} if $0<k<m$; a simplicial
set is the nerve of a category if and only if every inner horn has a unique
filler.  If $(A,H)$ is a simplicial set with hollowness
satisfying~\bref{part:high-dims} for $n=1$ then every inner horn is
admissible, hence, if~\bref{part:filler} and~\bref{part:unique} also hold,
has a unique filler: so $A$ is (the nerve of) a category.  Working out the
other conditions, we find that a weak $1$-category is a category equipped
with a set $H_1$ of isomorphisms containing all the identity maps and closed
under composition and inverses.  So given a weak $1$-category we obtain a
category by forgetting $H$; conversely, given a category we can take $H_1 =
\{ \textrm{all isomorphisms} \}$ (or $\{ \textrm{all identities} \}$) to
obtain a weak $1$-category.

\concept{$n=2$}

A weak $2$-category is a simplicial set $A$ equipped with subsets $H_1\sub
A[1]$ and $H_2 \sub A[2]$, satisfying certain axioms.  It appears that this
is the same as a bicategory equipped with a set $H_1$ of $1$-cells which are
equivalences and a set $H_2$ of 2-cells which are isomorphisms, satisfying
closure conditions similar to those under `$n=1$' above.

So, let $(A,H)$ be a weak $2$-category.  We construct a bicategory whose $0$-
and $1$-cells are the elements of $A[0]$ and $A[1]$; a 2-cell $f\go g$ is an
element of $A[2]$ of the form
\raisebox{-5mm}{%
\setlength{\unitlength}{1mm}
\begin{picture}(24,10)
\cell{4}{1}{c}{\zmark}
\cell{20}{1}{c}{\zmark}
\cell{12}{9}{c}{\zmark}
\put(4,1){\vector(1,0){15.5}}
\put(4,1){\vector(1,1){7.8}}
\put(12,9){\vector(1,-1){7.8}}
\put(12,3){\vector(0,1){4}}
\cell{3}{2}{tr}{\scriptstyle a}
\cell{13}{9}{l}{\scriptstyle a}
\cell{21}{2}{tl}{\scriptstyle b}
\cell{8}{5}{br}{\scriptstyle 1_a}
\cell{16}{5}{bl}{\scriptstyle g}
\cell{11}{1.5}{c}{\scriptstyle f}
\end{picture}},
where $1_a$ indicates a degenerate $1$-cell.  Composition of $1$-cells is
defined by making a random choice of hollow filler for each horn of shape
$\Lambda^1_2$; composition of $2$-cells is defined by filling in
$3$-dimensional horns $\Lambda^k_3$; identities are got from degeneracies.
See Fig.~\ref{fig:simplices}.  The associativity and unit isomorphisms are
certain hollow cells, and the coherence axioms hold because of the uniqueness
of certain fillers.

Conversely, let $B$ be a bicategory, and construct a weak $2$-category
$(A,H)$ as follows.  $A[0]$ and $A[1]$ are the sets of $0$- and $1$-cells of
$B$; an element of $A[2]$ as in Fig~\ref{fig:simplices}(b) is a $2$-cell
$\alpha_{012}: f_{02} \go f_{12}\of f_{01}$ in $B$.  (In general, an element
of $A[m]$ is a `strictly unitary colax morphism' $[m] \go B$, where $[m]$ is
regarded as a $2$-category whose only $2$-cells are identities.)  $H_1 \sub
A[1]$ is the set of $1$-cells which are equivalences, and $H_2$ is the set of
$2$-cells which are isomorphisms.  Then all the axioms for a weak 2-category
are satisfied.

\concept{Variant}

We could add to conditions~\bref{part:degen}--\bref{part:comp} on $(A,H)$ the
further condition that $H$ is maximal with respect
to~\bref{part:degen}--\bref{part:comp}: that is, if $(A,H')$ also
satisfies~\bref{part:degen}--\bref{part:comp} and $H'_m \supseteq H_m$ for
all $m$ then $H'=H$.  (Compare the issue of maximal atlases in the definition
of smooth manifold.)  With this addition, a weak $1$-category is essentially
just a category, and a weak $2$-category a bicategory.  This is in contrast
to the original \ds{St} (as analysed above), where the flexibility in the
choice of $H$ means that the correspondence between weak $1$-($2$-)categories
and (bi)categories is inexact.

\clearpage

\defnheading{X}		\label{p:x}

This definition is not intended to be rigorous, although it can be made so.

\concept{Opetopic Sets}

An \demph{opetopic set} $A$ is a commutative diagram of sets and functions
\begin{diagram}[width=1.7em,tight,alignlabels]
\cdots		&\ 	&\rTo^s_{\raisebox{-.5ex}{$t$}}	&&		&
		A'_2	&\rTo^s_{\raisebox{-.5ex}{$t$}}	&&		&
		A'_1	&\rTo^s_{\raisebox{-.5ex}{$t$}}	&&		&
		A'_0=A_0\\
\cdots 		&\ 	&\rdTo(4,2)			&&\ruTo(4,2)	&
			&\rdTo(4,2)			&&\ruTo(4,2)	&
			&\rdTo(4,2)			&&\ruTo(4,2)	&
			\\
\cdots		&\ 	&\rTo_t^{\raisebox{.7ex}{$s$}}	&&		&
		A_2	&\rTo_t^{\raisebox{.7ex}{$s$}}	&&		&
		A_1	&\rTo_t^{\raisebox{.7ex}{$s$}}	&&		&
		A_0	\\
\end{diagram}
where for each $m\geq 1$, the set $A'_m$ and the functions $s: A'_m \go
A'_{m-1}$ and $t: A'_m \go A_{m-1}$ are defined from the sets $A_m, A'_{m-1},
A_{m-1}, \ldots, A'_1, A_1, A_0$ and the functions $s,t$ between them in the
following way.

An element $a \in A_0$ is regarded as a $0$-cell, and drawn \gzeros{a}.  An
element $f \in A_1$ is regarded as a $1$-cell $\gfsts{a} \gones{f}
\glsts{b}$, where $a=s(f)$ and $b=t(f)$.  $A'_1$ is the set of `1-pasting
diagrams' in $A$, that is, diagrams of 1-cells pasted together, that is,
paths $\gfsts{a_0} \gones{f_1} \gblws{a_1} \gones{f_2} \ldots \gones{f_k}
\glsts{a_k}$ ($k\geq 0$) in $A$.  An element $\alpha \in A_2$ has a source
$s(\alpha)$ of this form and a target $t(\alpha)$ of the form $\gfsts{a_0}
\gones{g} \glsts{a_k}$, and is drawn as
\begin{equation}	\label{eq:two-ope}
\raisebox{-6.5mm}{%
\setlength{\unitlength}{1mm}
\begin{picture}(36,13)(0,-2)
\cell{0}{0}{c}{\zmark}
\cell{6}{8}{c}{\zmark}
\cell{36}{0}{c}{\zmark}
\put(0,0){\vector(3,4){6}}
\put(6,8){\vector(3,1){9}}
\put(30,8){\vector(3,-4){6}}
\put(0,0){\vector(1,0){36}}
\put(18,7){\vector(0,-1){5}}
\cell{22}{9.5}{c}{\cdots}
\cell{-2.5}{0}{c}{a_0}
\cell{4}{9}{c}{a_1}
\cell{39}{0}{c}{a_k}
\cell{1}{5}{c}{f_1}
\cell{10}{11.5}{c}{f_2}
\cell{35.5}{5}{c}{f_k}
\cell{18}{-1.5}{c}{g}
\cell{20}{4.5}{c}{\alpha}
\end{picture}}
\end{equation}
Next, $A'_2$ is the set of `2-pasting diagrams', that is, diagrams of cells
of the form~\bref{eq:two-ope} pasted together, such as
\begin{equation}	\label{eq:two-pd}
\raisebox{-10.5mm}{%
\setlength{\unitlength}{1mm}
\begin{picture}(50,21)(-5,-3)
\cell{0}{0}{c}{\zmark}
\cell{-5}{7.5}{c}{\zmark}
\cell{0}{12.5}{c}{\zmark}
\cell{7.5}{10}{c}{\zmark}
\cell{22.5}{10}{c}{\zmark}
\cell{26.25}{15}{c}{\zmark}
\cell{36.25}{12.5}{c}{\zmark}
\cell{43.75}{7.5}{c}{\zmark}
\cell{36.25}{2.5}{c}{\zmark}
\cell{30}{0}{c}{\zmark}
\put(0,0){\vector(-2,3){5}}
\put(-5,7.5){\vector(1,1){5}}
\put(0,12.5){\vector(3,-1){7.5}}
\put(7.5,10){\vector(1,0){15}}
\put(22.5,10){\vector(3,4){3.75}}
\put(26.25,15){\vector(4,-1){10}}
\put(36.25,12.5){\vector(3,-2){7.5}}
\put(43.75,7.5){\vector(-3,-2){7.5}}
\put(36.25,2.5){\line(-5,-2){6.25}}
\put(30,0){\vector(-3,-1){0}}
\put(36.25,12.5){\vector(0,-1){10}}
\put(0,0){\vector(3,4){7.5}}
\put(22.5,10){\vector(3,-4){7.5}}
\put(0,0){\vector(1,0){30}}
\put(15,7.5){\vector(0,-1){5}}
\put(-2,9){\vector(1,-1){4}}
\put(32,10.5){\vector(-1,-1){4}}
\put(41.25,7){\vector(-1,0){4}}
\cell{-1}{-1}{c}{\scriptstyle a_0}
\cell{-7}{7.5}{c}{\scriptstyle a_1}
\cell{0}{14}{c}{\scriptstyle a_2}
\cell{8.5}{11.5}{c}{\scriptstyle a_3}
\cell{21}{11}{c}{\scriptstyle a_4}
\cell{27}{16.5}{c}{\scriptstyle a_5}
\cell{37}{14}{c}{\scriptstyle a_6}
\cell{46}{7.5}{c}{\scriptstyle a_7}
\cell{38}{1.5}{c}{\scriptstyle a_8}
\cell{30}{-1.5}{c}{\scriptstyle a_9}
\cell{-3.5}{3}{c}{\scriptstyle f_1}
\cell{-3.5}{11}{c}{\scriptstyle f_2}
\cell{4}{12.5}{c}{\scriptstyle f_3}
\cell{15}{11.5}{c}{\scriptstyle f_4}
\cell{23.5}{13.5}{c}{\scriptstyle f_5}
\cell{31.5}{15}{c}{\scriptstyle f_6}
\cell{41.5}{11}{c}{\scriptstyle f_7}
\cell{41.5}{4}{c}{\scriptstyle f_8}
\cell{34}{0.25}{c}{\scriptstyle f_9}
\cell{34.5}{6.5}{c}{\scriptstyle f_{10}}
\cell{25}{4}{c}{\scriptstyle f_{11}}
\cell{5}{4}{c}{\scriptstyle f_{12}}
\cell{15}{-1.5}{c}{\scriptstyle f_{13}}
\cell{17}{5}{c}{\scriptstyle \alpha_1}
\cell{1}{8}{c}{\scriptstyle \alpha_2}
\cell{28.5}{9.5}{c}{\scriptstyle \alpha_3}
\cell{40}{8}{c}{\scriptstyle \alpha_4}
\end{picture}}
\end{equation}
Note that the arrows go in compatible directions: e.g.\ the target or
`output' edge $f_{11}$ of $\alpha_3$ is a source or `input' edge of
$\alpha_1$.  The source of this element of $A'_2$ is $\gfsts{a_0} \gones{f_1}
\cdots \gones{f_9} \glsts{a_9} \in A'_1$, and the target is $f_{13} \in A_1$.
Next, if $\gamma \in A_3$ and $s(\gamma)$ is~\bref{eq:two-pd} then
$t(\gamma)$ is of the form~\bref{eq:two-ope} with $k=9$ and $g = f_{13}$, and
we picture $\gamma$ as a $3$-dimensional cell with a flat bottom face
(labelled $\alpha$) and four curved faces on top (labelled $\alpha_1,
\alpha_2, \alpha_3, \alpha_4$).  Carrying on, $A'_3$ is the set of
$3$-pasting diagrams, $A_4$ is the set of $4$-cells, etc.

We need some terminology concerning cells.  Let $\Phi \goby{\alpha} g$ be an
$m$-cell: that is, let $\alpha \in A_m$ with $s(\alpha) = \Phi \in A'_{m-1}$
and $t(\alpha) = g \in A_{m-1}$.

For any $p$-cell $e$, there is a $p$-pasting diagram $\langle e \rangle$
consisting of $e$ alone.  If $\langle g \rangle \goby{\beta} h$ then
$\alpha$ and $\beta$ can be pasted to obtain $(\Phi \goby{\beta_*(\alpha)} h)
\in A'_m$.

The \demph{faces} of $\Phi$ are the $(m-1)$-cells which have been pasted
together to form it, e.g.\ if $\alpha$ is as in~\bref{eq:two-ope} then $\Phi$
has faces $f_1, \ldots, f_k$, and the $\Phi$ of~\bref{eq:two-pd} has faces
$\alpha_1, \alpha_2, \alpha_3, \alpha_4$.  If $f$ is a face of $\Phi$ and $e$
is a cell \demph{parallel} to $f$ (i.e.\ $e\in A_{m-1}$ with $s(e)=s(f)$,
$t(e)=t(f)$) then we obtain a new pasting diagram $\Phi(e/f) \in A_{m-1}$ by
replacing $f$ with $e$ in $\Phi$.  (Read $\Phi(e/f)$ as `$\Phi$ with $e$
replacing $f$'.)  If also $\langle e \rangle \goby{\beta} f$ then $\alpha$
and $\beta$ can be pasted to obtain $(\Phi(e/f) \goby{\beta^*(\alpha)} g) \in
A'_m$.

\concept{Universal Cells}

Let $A$ be an opetopic set and fix $n\geq 0$.  We define what it means for a
cell $\Phi \goby{\epsln} g$ of $A$ to be `universal', and, when $f$ is a face
of $\Phi$, what it means for $f$ to be `liminal' in the cell.  The
definitions depend on $n$, so I should really say `$n$-universal' rather than
just `universal', and similarly `$n$-liminal'; but I will drop the `$n$'
since it is regarded as fixed.  

The two definitions are given inductively in an interdependent way.

\paragraph{Universality} 

Let $m\geq n+1$.  A cell $(\Phi \goby{\epsln} g) \in A_m$ is
\demph{universal} if whenever $(\Phi \goby{\epsln'} g')\in A_m$, then
$\epsln'=\epsln$.

Let $1\leq m\leq n$.  A cell $(\Phi \goby{\epsln} g) \in A_m$ is
\demph{universal} if
\begin{enumerate}
\item 	\label{part:univ-univ}
for every $\alpha: \Phi \go h$, there exist $\ovln{\alpha}: \langle
g\rangle \go h$ and a universal cell $U: \ovln{\alpha}_*(\epsln) \go \alpha$,
and 
\item 	\label{part:univ-liminal}
for every $\alpha: \Phi \go h$, $\ovln{\alpha}: \langle g\rangle \go h$
and universal $U: \ovln{\alpha}_*(\epsln) \go \alpha$, $\ovln{\alpha}$ is
liminal in $U$.
\end{enumerate}

\paragraph{Liminality}

Let $m\geq 1$, let $(\Phi \goby{\epsln} g) \in A_m$, and let $f$ be a face of
$\Phi$.  Then $f$ is \demph{liminal in $\epsln$} if $m\geq n+2$ or
\begin{enumerate}
\item 	\label{part:liminal-univ}
for every cell $e$ parallel to $f$ and $\beta: \Phi(e/f) \go g$, there
exist $\ovln{\beta}: \langle e\rangle \go f$ and a universal cell $U:
\ovln{\beta}^*(\epsln) \go \beta$, and
\item 	\label{part:liminal-liminal}
for every $e$ parallel to $f$, $\beta: \Phi(e/f) \go g$, $\ovln{\beta}:
\langle e\rangle \go f$  and universal $U: \ovln{\beta}^*(\epsln) \go \beta$,
$\ovln{\beta}$ is liminal in $U$.
\end{enumerate}

\concept{The Definition}

Let $n\geq 0$.  A \demph{weak $n$-category} is an opetopic set $A$ satisfying
\begin{description}
\item[existence of universal fillers:] for every $m\geq 0$ and $\Phi \in
A'_m$, there exists a universal cell of the form $\Phi \goby{\epsln} g$
\item[closure of universals under composition:] if $m\geq 2$, $(\Phi
\goby{\epsln} g) \in A_m$, each face of $\Phi$ is universal, and $\epsln$ is
universal, then $g$ is universal.
\end{description}

\clearpage

\lowdimsheading{X}

First consider high-dimensional cells in an $n$-category $A$.  For every
$m\geq n+1$ and $\Phi \in A'_{m-1}$, there is a unique $\epsln \in A_m$ whose
source is $\Phi$: in other words, $s: A_m \go A'_{m-1}$ is a bijection.  This
means that the entire opetopic set is determined by the part of dimension
$\leq n$ and the map $t: A_{n+1} \go A_n$.  Identifying $A_{n+1}$ with
$A'_n$, this map $t$ assigns to each $n$-pasting diagram $\Phi$ the target
$g$ of the unique $(n+1)$-cell with source $\Phi$, and we regard $g$ as the
composite of $\Phi$.  (In general, we may regard the target of a universal
cell as a composite of its source; note that all cells of dimension $>n$ are
universal.)  This composition of $n$-cells is strictly associative and
unital.  A weak $n$-category therefore consists of a commutative diagram 
\begin{equation}	\label{eq:trunc-ope}
\begin{diagram}[width=1.7em,tight,alignlabels]
A'_n	&\rTo^s_{\raisebox{-.5ex}{$t$}}	&&		&
\ &\cdots&\	
	&\rTo^s_{\raisebox{-.5ex}{$t$}}	&&		&
A'_1	&\rTo^s_{\raisebox{-.5ex}{$t$}}	&&		&
A'_0=A_0\\
\dTo<\of &\rdTo(4,2)			&&\ruTo(4,2)	&
&&	
	&\rdTo(4,2)			&&\ruTo(4,2)	&
	&\rdTo(4,2)			&&\ruTo(4,2)	&
	\\
A_n 	&\rTo_t^{\raisebox{.7ex}{$s$}}	&&		&
\ &\cdots&\	
	&\rTo_t^{\raisebox{.7ex}{$s$}}	&&		&
A_1	&\rTo_t^{\raisebox{.7ex}{$s$}}	&&		&
A_0	\\
\end{diagram}
\end{equation}
of sets and functions such that $n$-dimensional composition $\of$ obeys
strict laws and the defining conditions on existence and closure of
universals hold in lower dimensions.

\concept{$n=0$}

A weak $0$-category is a set $A_0$ together with a function $\of: A_0 \go
A_0$ obeying strict laws---which say that $\of$ is the
identity.  So a weak $0$-category is just a set.

\concept{$n=1$}

When $n=1$, diagram~\bref{eq:trunc-ope} is a directed graph $A_1
\parpair{s}{t} A_0$ together with a map $\of: A'_1 \go A_1$ compatible with
the source and target maps: in other words, assigning to each string of edges
\[
a_0 \goby{f_1} a_1 \goby{f_2} \ \cdots \ \goby{f_k} a_k
\]
in $A$ a new edge $a_0 \goby{(f_k \sof\cdots\sof f_1)} a_k$.  The axioms on
$\of$ say that
\[
((f_k^{r_k} \of\cdots\of f_k^1) \of\cdots\of (f_1^{r_1} \of\cdots\of f_1^1))
= (f_k^{r_k} \of\cdots\of f_1^1),
\diagspace
(f) = f
\]
---in other words, that $A$ forms a category.  A weak $1$-category is
therefore a category satisfying the extra conditions that every object is the
domain of some universal morphism and that the composite of universal
morphisms is universal.  I claim that the universal morphisms are the
isomorphisms, so that these conditions hold automatically and a weak
$1$-category is just a category. 

So: a morphism $\epsln: a \go b$ is universal if~\bref{part:univ-univ} every
morphism $\alpha: a \go c$ factors as $\alpha = \ovln{\alpha} \of \epsln$,
and~\bref{part:univ-liminal} such an $\ovln{\alpha}$ is always liminal in the
unique 2-cell
\raisebox{-6mm}{%
\setlength{\unitlength}{1mm}
\begin{picture}(24,12)(0,-2)
\cell{4}{1}{c}{\zmark}
\cell{20}{1}{c}{\zmark}
\cell{12}{9}{c}{\zmark}
\put(4,1){\vector(1,0){15.5}}
\put(4,1){\vector(1,1){7.8}}
\put(12,9){\vector(1,-1){7.8}}
\put(12,6){\vector(0,-1){4}}
\cell{3}{2}{tr}{\scriptstyle a}
\cell{13}{9}{l}{\scriptstyle b}
\cell{21}{2}{tl}{\scriptstyle c}
\cell{8}{5}{br}{\scriptstyle \epsln}
\cell{16}{5}{bl}{\scriptstyle \ovln{\alpha}}
\cell{12}{-0.5}{c}{\scriptstyle \alpha}
\cell{13}{4}{l}{\scriptstyle U}
\end{picture}}.
Liminality of $\ovln{\alpha}$ in $U$ means, in turn, that if $\twid{\alpha}:
b \go c$ satisfies $\twid{\alpha} \of \epsln = \alpha$ then $\twid{\alpha} =
\ovln{\alpha}$.  (This is part~\bref{part:liminal-univ} of the definition of
liminality; part~\bref{part:liminal-liminal} holds trivially.)  So $\epsln: a
\go b$ is universal if and only if every morphism out of $a$ factors uniquely
through $\epsln$, which holds if and only if $\epsln$ is an isomorphism.

\concept{$n=2$}

A weak $2$-category is essentially the same thing as a bicategory.  More
precisely, the category of bicategories and weak functors is equivalent to
the category whose objects are weak $2$-categories and whose morphisms are
those maps of opetopic sets which send universal cells to universal cells.  
The equivalence works as follows.

Given a bicategory $B$, define a weak $2$-category $A$ by taking $A_0$ and
$A_1$ to be the sets of $0$- and $1$-cells in $B$, and a
$2$-cell~\bref{eq:two-ope} in $A$ to be a $2$-cell $(f_k \of \cdots \of f_1)
\go g$ in $B$.  Here $(f_k \of \cdots \of f_1)$ is defined inductively as
$f_k \of (f_{k-1} \of\cdots\of f_1)$ if $k\geq 1$, or as $1$ if $k=0$; any
other iterated composite would do just as well.  Composition $\of: A'_2 \go
A_2$ is pasting of $2$-cells in $B$.  Then a 1-cell (respectively, 2-cell) in
$A$ turns out to be universal if and only if the corresponding 1-cell
(2-cell) in $B$ is an equivalence (isomorphism), and it follows that $A$ is a
weak $2$-category.

Conversely, take a weak $2$-category $A$ and construct a bicategory $B$ as
follows.  The 0- and 1-cells of $B$ are just those of $A$, and a 2-cell of
$B$ is an element of $A_2$ of the form $\topeavar{\scriptstyle
a}{\scriptstyle b}{\scriptstyle f}{\scriptstyle g}{\scriptstyle \alpha}$
(i.e.\ $\alpha: \langle f \rangle \go g$).  For each diagram $\gfsts{a}
\gones{f} \gblws{b} \gones{g} \glsts{c}$ of 1-cells, choose at random a
universal filler
\raisebox{-6mm}{%
\setlength{\unitlength}{1mm}
\begin{picture}(24,12)(0,-2)
\cell{4}{1}{c}{\zmark}
\cell{20}{1}{c}{\zmark}
\cell{12}{9}{c}{\zmark}
\put(4,1){\vector(1,0){15.5}}
\put(4,1){\vector(1,1){7.8}}
\put(12,9){\vector(1,-1){7.8}}
\put(11,6){\vector(0,-1){4}}
\cell{3}{2}{tr}{\scriptstyle a}
\cell{13}{9}{l}{\scriptstyle b}
\cell{21}{2}{tl}{\scriptstyle c}
\cell{8}{5}{br}{\scriptstyle f}
\cell{16}{5}{bl}{\scriptstyle g}
\cell{12}{-0.5}{c}{\scriptstyle g\sof f}
\cell{11.5}{4}{l}{\scriptstyle \epsln_{f,g}}
\end{picture}},
where by definition $g\of f = t(\epsln_{f,g})$; this defines composition of
1-cells.  Vertical composition of 2-cells comes from $\of: A'_2 \go A_2$.  To
define the horizontal composite of $\gfsts{a} \gtwos{f}{g}{\alpha} \gfbws{a'}
\gtwos{f'}{g'}{\alpha'} \glsts{a''}$, consider pasting $\epsln_{g,g'}$ to
$\alpha$ and $\alpha'$, and then use the universality of $\epsln_{f,f'}$.
Next observe that given two universal fillers $\Phi \goby{\epsln} g$, $\Phi
\goby{\epsln'} g'$ for a 1-pasting diagram $\Phi = (\gfsts{a_0} \gones{f_1}
\ldots \gones{f_k} \glsts{a_k})$, there is a unique 2-cell $\langle g \rangle
\goby{\delta} g'$ such that the composite $\of (\delta_* (\epsln))$ is
$\epsln'$.  Applying this observation to a certain pair of universal fillers
for $(\gfsts{} \gones{f} \gblws{} \gones{g} \gblws{} \gones{h} \glsts{})$
gives the associativity isomorphism, and the word `unique' in the observation
gives the pentagon axiom.  Identities work similarly, where this time we
choose a random universal filler for each degenerate 1-pasting diagram
\gzeros{a}.

\section*{Further Reading}	\label{p:biblio}

This section contains the references and historical notes missing in the
main text.  It is not meant to be a survey of the literature.  Where I have
omitted relevant references it is almost certainly a result of my own
ignorance, and I hope that the authors will forgive me.

First are some references to introductory and general material, and a very
brief account of the history of higher-dimensional category theory.  Then
there are references for each of the sections in turn: `Background',
followed by the ten definitions.  Finally there are references to some
proposed definitions of $n$-category which I didn't include, and a very few
references to areas of mathematics related to $n$-categories. 

Citations such as \url{math\dt CT\slsh 9810058} and \url{alg-geom\slsh
9708010} refer to the electronic mathematics archive at
\url{http:\dblslsh arXiv\dt org}.
Readers unfamiliar with the archive may find it easiest to go straight to the
address of the form \url{http:\dblslsh arXiv\dt org\slsh abs\slsh math\dt
CT\slsh 9810058} .

\subsection*{Introductory Texts}

Introductions to $n$-categories come slanted towards various different
audiences.  One for theoretical computer scientists and logicians is 
\bibentry{Pow}{%
A. J. Power, 
Why tricategories?, 
\jnl{Information and Computation}{120}{1995}{no.~2, 251--262}; 
also 
LFCS report ECS-LFCS-94-289, April 1994, 
\web{http:\dblslsh www\dt lfcs\dt informatics\dt ed\dt ac\dt uk}}%
and another with a logical slant, but this time with foundational concerns,
is 
\bibentry{MakTCF}{%
M. Makkai,
Towards a categorical foundation of mathematics,
\contrib{Logic Colloquium '95 (Haifa)}
{Lecture Notes in Logic~\textbf{11}}
{Springer}{1998}{153--190}.}%

Moving to introductions for those more interested in topology, geometry and
physics, one which starts at a very basic level (the definition of category)
is 
\bibentry{BaezTNC}{%
John C. Baez, 
A tale of $n$-categories,
\web{http:\dblslsh math\dt ucr\dt edu\slsh home\slsh baez\slsh week73\dt
html}, 1996--97.}%
With similar themes but at a more advanced level, there are
\bibentry{BaezINC}{%
John C. Baez,
An introduction to $n$-categories,
\contrib{Category theory and computer science (Santa Margherita Ligure,
1997)}
{Lecture Notes in Computer Science}{1290}
{Springer}{1997}{1--33};
also
\epr{q-alg\slsh 9705009}{1997}}%
(especially sections~1--3) and
\bibentry{THDCTRI}{%
Tom Leinster, 
Topology and higher-dimensional category theory: the rough idea,
\eprint{math\dt CT\slsh 0106240}{2001}{15}.}%
The ambitious might, if they can find a copy, like to look at the highly
discursive 600-page letter of Grothendieck to Quillen,
\bibentry{Gro}{%
A. Grothendieck, 
Pursuing stacks,
manuscript, 1983,}%
in which (amongst many other things) the idea is put that tame topology is
really the study of weak $\omega$-groupoids.  A more accessible discussion of
what higher-dimensional algebra might `do', especially in the context of
topology, is
\bibentry{Bro}{%
Ronald Brown,
Higher dimensional group theory,
\web{http:\dblslsh www\dt bangor\dt ac\dt uk\slsh $\sim$mas010}.}%

\subsection*{General Comments and History}

The easiest way to begin a history of $n$-categories is as follows.

\begin{sloppypar}
$0$-categories---sets or classes---came into the mathematical consciousness
around the end of the 19th century.  $1$-categories---categories---arrived in
the middle of the 20th century.  Strict $2$-categories and, implicitly,
strict $n$-categories, made their presence felt around the late 1950s and
early 1960s, with the work of Ehresmann \cyte{Ehr}.  Weak $2$-categories were
first introduced by B\'enabou \cyte{Ben} in 1967, under the name of
bicategories, and thereafter the question was in the air: `what might a weak
$n$-category be?'  The first precise proposal for a definition was given by
Street \cyte{StrAOS} in 1987.  This was followed by three more proposals
around 1995: Baez and Dolan's \cyte{BDHDA3}, Batanin's \cyte{BatMGC}, and
Tamsamani's \cyte{TamNNC}.  A constant stream of further proposed definitions
has issued forth since then, and will doubtless continue for a while.
Work on low values of $n$ was also going on at the same time: an
axiomatic definition of weak $3$-category was proposed in
\bibentry{GPS}{%
R. Gordon, A. J. Power, Ross Street, 
\emph{Coherence for Tricategories}, 
Memoirs of the American Mathematical Society~\textbf{117}, no.~558,
1995,}%
and a proposal in similar vein for $n=4$ was made in
\bibentry{TriDT}{%
Todd Trimble, 
The definition of tetracategory, 
manuscript, 1995.}%
Crucially, it was shown in~\cyte{GPS} that not every tricategory is
equivalent to a strict $3$-category (in contrast to the situation for $n=2$),
from which it follows that the theory of weak $n$-categories is genuinely
different from that of strict ones. 
\end{sloppypar}

But this is far too simplistic.  A realistic history must take account of
categorical structures other than $n$-categories \emph{per se}: for instance,
the various kinds of monoidal category (plain, symmetric, braided,
tortile/ribbon, \ldots), and of monoidal 2-categories and monoidal
bicategories.  The direct importance of these is that a monoidal category is
a bicategory with only one $0$-cell, and similarly a braided monoidal
category is a tricategory with only one $0$-cell and one $1$-cell.  The basic
reference for braided monoidal categories is
\bibentry{JS}{%
Andr\'e Joyal, Ross Street,
Braided tensor categories, 
\jnl{Advances in Mathematics}{102}{1993}{no. 1, 20--78},}%
and they can also be found in the new edition of Mac Lane's book
\cyte{MacCWM}.

Moreover, around the same time as the theory of $n$-categories was starting
to develop, another theory was emerging with which it was later to converge:
the theory of multicategories and operads.  Multicategories first appeared in
\bibentry{Lam}{%
Joachim Lambek, 
Deductive systems and categories II: standard constructions
and closed categories,
\contribau{Category Theory, Homology Theory and their Applications, I
(Battelle Institute Conference, Seattle, 1968, Vol. One)}
{ed.\ Hilton}{Lecture Notes in Mathematics~\textbf{86}}
{Springer}{1969}{76--122}.}%
(The definition is on page~103.)  A multicategory is like a category, but
each arrow has as its source or input a \emph{sequence} of objects (and, as
usual, as its target or output a single object).  An operad is basically just
a multicategory with only one object.  For this reason, multicategories are
sometimes called `coloured operads', and the objects are then named after
colours (black, white, etc.).  The development of operads is generally
attributed to Boardman, Vogt and May:
\bibentry{BV}{%
J. M. Boardman, R. M. Vogt, 
\emph{Homotopy Invariant Algebraic Structures on
Topological Spaces},
Lecture Notes in Mathematics~\textbf{347},
Springer, 1973,}%
\bibentry{MayGIL}{%
J. P. May, 
\emph{The Geometry of Iterated Loop Spaces},
Lectures Notes in Mathematics~\textbf{271}, Springer, 1972,}%
although I am told that essentially the same idea was the subject of 
\bibentry{Laz}{%
Michel Lazard,
Lois de groupes et analyseurs,
\jnl{Annales Scientifiques de l'\'Ecole Normale Sup\'erieure (3)}
{72}{1955}{299--400}}%
(where operads go by the name of `analyseurs').

It seems to have taken a long time before it was realized that operads and
multicategories were so closely related; I do not know of any pre-1995 text
which mentions both Lambek and Boardman-Vogt or May in its bibliography.
This can perhaps be explained by the different fields in which they were
being studied: multicategories were introduced in the context of logic and
found application in linguistics, whereas operads were used for the theory of
loop spaces.  Moreover, if one uses the terms in their original senses then
it is not strictly true that an operad is the same thing as a one-object
multicategory; operads are also equipped with a symmetric structure, and the
`hom-sets' (sets of operations) are topological spaces rather than just sets.
(It is also very natural to consider multicategories with both these pieces
of extra structure, but historically this is beside the point.)

Many short introductions to operads have appeared as section~1 of papers by
topologists and quantum algebraists.  The interested reader may also find
useful the following texts dedicated to the subject:
\bibentry{MayDOA}{%
J. P. May, 
Definitions: operads, algebras and modules,
\contrib{Operads: Proceedings of Renaissance Conferences (Hartford,
CT\slsh Luminy, 1995)} 
{Contemporary Mathematics~\textbf{202}}
{AMS}{1997}{1--7};
also
\web{http:\dblslsh www\dt math\dt uchicago\dt edu\slsh $\sim$may},}%
\bibentry{MayOAM}{%
J. P. May, 
Operads, algebras and modules,
\contrib{Operads: Proceedings of Renaissance Conferences (Hartford,
CT\slsh Luminy, 1995)} 
{Contemporary Mathematics~\textbf{202}}
{AMS}{1997}{15--31};
also
\web{http:\dblslsh www\dt math\dt uchicago\dt edu\slsh $\sim$may},}%
\bibentry{MSS}{%
Martin Markl, Steve Shnider, Jim Stasheff, 
\emph{Operads in Algebra, Topology and Physics}, 
book in preparation.}%

A glimpse of the role of operads and multicategories in higher-dimensional
category theory can be seen in the definitions of weak $n$-category above.
Often the `operads' and `multicategories' used are not the original kinds,
but more general kinds adapted for the different shapes and dimensions which
occur in the subject; for more references, see `Definitions \ds{B} and
\ds{L}' below.

\subsection*{Background}

\paragraph*{Category Theory}

Almost any book on the subject will provide the necessary 
background.  The second edition of the classic book by Mac Lane,
\bibentry{MacCWM}{%
Saunders Mac Lane,
\emph{Categories for the Working Mathematician},
second edition, Graduate Texts in Mathematics~\textbf{5}, Springer, 1998,}%
is especially useful, containing as it does two new chapters on such
topics as bicategories and nerves of categories.  

\paragraph*{Strict $n$-Categories}

I do not know of any good text introducing strict $n$-categories.
Ehresmann's original book
\bibentry{Ehr}{%
Charles Ehresmann, 
\emph{Cat\'egories et Structures},
Dunod, Paris, 1965}%
could be consulted, but is generally regarded as a very demanding read.  
Probably more useful is 
\bibentry{KS}{%
G. M. Kelly, Ross Street, 
Review of the elements of $2$-categories, 
\contrib{Category Seminar (Sydney, 1972\slsh 1973)}
{Lecture Notes in Mathematics~\textbf{420}}
{Springer}{1974}{75--103},}%
which only covers strict $2$-categories (traditionally just called
`2-categories') but should give a good idea of strict $n$-categories for
general $n$.  This could usefully be supplemented by
\bibentry{EK}{%
Samuel Eilenberg, G. Max Kelly, 
Closed categories,
\contb{Proceedings of Conference on Categorical Algebra (La Jolla,
California, 1965)}
{Springer}{1966}{421--562}}%
(see e.g.\ page~552), which also covers enrichment.  For another reference on
enriched categories, see chapter~6 of
\bibentry{Borx2}{%
Francis Borceux, 
\emph{Handbook of Categorical Algebra 2: Categories and
Structures},
Encyclopedia of Mathematics and its Applications~\textbf{51},
Cambridge University Press, 1994.}%

\paragraph*{Bicategories}

Bicategories were first explained by B\'enabou:
\bibentry{Ben}{%
Jean B\'enabou, 
Introduction to bicategories,
\contribau{Reports of the Midwest Category Seminar}
{ed.\ B\'enabou et al}
{Lecture Notes in Mathematics~\textbf{47}}{Springer}{1967}{1--77},}%
and further important work on them is in
\bibentry{Gray}{%
John W. Gray,
\emph{Formal Category Theory: Adjointness for 2-Categories}, 
Lecture Notes in Mathematics~\textbf{391}, Springer, 1974.}%
(At least) two texts contain summaries of the `basic theory' of bicategories:
that is, the definitions of bicategory and of weak functor (homomorphism),
transformation and modification between bicategories, together with the
result that any bicategory is in a suitable sense equivalent to a strict
2-category.  These are section 9 of
\bibentry{StrCS}{%
Ross Street,
Categorical structures, 
\contbau{Handbook of Algebra~\textbf{1}}
{ed.\ M. Hazewinkel}
{North-Holland}{1996}{529--577}}%
and the whole of
\bibentry{BB}{%
Tom Leinster,
Basic bicategories,
\eprint{math\dt CT\slsh 9810017}{1998}{11}.}%

\subsection*{Definition \ds{Tr}}

The definition was given in a talk,
\bibentry{TriWFN}{%
Todd Trimble, 
What are `fundamental $n$-groupoids'?,
seminar at \mbox{DPMMS}, Cambridge, 
24 August 1999,}%
and has not been written up previously.  Trimble used the term `flabby
$n$-category' rather than `weak $n$-category'.  

As the title of the talk suggests, the idea was not to develop the weakest
possible notion of $n$-category, but to provide (in his words) `a sensible
niche for discussing fundamental $n$-groupoids'.  In a world where all the
definitions have been settled, it may be that fundamental $n$-groupoids of
topological spaces have certain special features (other than the
invertibility of their cells) not shared by all weak $n$-categories.  Thus it
may be that the word `weak' is less appropriate for definition \ds{Tr} than
the other definitions.

Evidence that this is the case comes from two directions.  Firstly, the maps
$\gamma_{a_0, \ldots, a_k}$ describing composition of hom-$(n-1)$-categories
in an $n$-category are \emph{strict} $(n-1)$-functors.  This corresponds to
having strict interchange laws.  It therefore seems likely that a precise
analysis of $n=3$ would show that every weak $3$-category gives rise to a
tricategory (in a similar manner to $n=2$) but that not every tricategory is
triequivalent to one arising from a weak $3$-category.  Secondly, forget
\ds{Tr} for the moment and consider in naive terms what the fundamental
$n$-groupoid of a space $X$ might look like when $n\geq 3$.  $0$-cells would
be points of $X$, $1$-cells could very reasonably be maps $[0,1] \go X$, and
similarly $2$-cells could be maps $[0,1]^2 \go X$ satisfying suitable
boundary conditions.  Composition of $1$-cells could be defined by travelling
each path at double speed, in the fashion customary to homotopy theorists,
and similarly for vertical and horizontal composition of $2$-cells.  The
point now is that although none of these compositions is strictly associative
or unital, the interchange law between horizontal and vertical $2$-cell
composition \emph{is} obeyed strictly.  This provides the kind of `special
feature' of fundamental $n$-groupoids referred to above.

Prospects for comparing \ds{Tr} with \ds{B} and \ds{L} look bright: it seems
very likely that weak $n$-categories according to \ds{Tr} are just the
algebras for a certain globular $n$-operad (in the sense of \ds{B} or
\ds{L}).

Other ideas on fundamental $n$-groupoids, $n$-categories, and how they tie
together can be found in Grothendieck's letter \cyte{Gro}.  More practical
material on fundamental 1- and 2-groupoids is in
\bibentry{KP}{%
K. H. Kamps, T. Porter, 
\emph{Abstract Homotopy and Simple Homotopy Theory},
World Scientific Publishing Co., 1997.}%

For more on operads, see the references under `General Comments and History'
above.  The name of the operad $E$ was not only chosen to stand for
`endpoint-preserving', but also because it comes after $D$ for `disk'---the
idea being that $E$ is something like the little disks operad $D$ (crucial
in the theory of loop spaces).  A touch more precisely, $E$ seems to play the
same kind of role for paths as $D$ does for closed loops.

More about how bicategories comes from the operad of trees can be found in
Appendix~A (and Chapter~1) of my thesis, \cyte{OHDCT}, and in my \cyte{WMC}.

\subsection*{Definition \ds{P}}

Definition \ds{P} of weak $\omega$-category is in
\bibentry{Pen}{%
Jacques Penon, 
Approche polygraphique des $\infty$-categories non strictes, 
\jnl{Cahiers de Topologie et G\'eom\'etrie Diff\'erentielle}
{40}{1999}{no.~1, 31--80}.}%
His chosen term for weak $\omega$-category is `prolixe', whose closest
English translation is perhaps `waffle'.  As far as I can see there is no 
actual definition of weak $n$-category or $n$-dimensional prolixe in the
paper, although he clearly has one in mind on page~48:
\begin{quote}
Les prolixes de dimension $\leq 2$ s'identifient exactement aux
bicat\'egories [\ldots] la preuve de ce r\'esultat sera montr\'e dans un
article ult\'erieur
\end{quote}
(`waffles of dimension $\leq 2$ correspond exactly to bicategories [\ldots]
the proof of this result will be given in a forthcoming paper').

Other translations: my category $\ftrcat{\scat{R}}{\Set}$ of reflexive
globular sets is his category $\infty\hyph\mathbb{G}rr$ of reflexive
$\infty$-graphs; my $s$ and $t$ are his $s$ and $b$; my strict
$\omega$-categories are his $\infty$-categories; my category \cat{Q} is
called by him $\mathbb{E}tC$, the category of \'etirements cat\'egoriques
(`categorical stretchings'); my contractions $\gamma$ are written
$[\dashbk,\dashbk]$ (with the arguments reversed: $\gamma_m(f_0,f_1) =
[f_1,f_0]$); and my adjunction $F\ladj U$ is called $\hat{\mathcal{E}} \ladj
\hat{V}$.

The word `magma' is borrowed from Bourbaki, who used it to mean a set
equipped with a binary operation.  It is a slightly inaccurate borrowing, in
that $\omega$-magmas are equipped with (nominal) identities as well as binary
compositions; put another way, it would have been more suitable if Bourbaki
had used the word to mean a set equipped with a binary operation and a
distinguished basepoint.

It seems plausible that Penon's construction can be generalized to provide
weak versions of structures other than $\omega$-categories (e.g.\
up-to-homotopy topological monoids).  Batanin has done something
precise along the lines of generalizing Penon's definition and comparing it
to his own:
\bibentry{BatPMW}{%
M. A. Batanin,  
On the Penon method of weakening algebraic structures, 
to appear in \emph{Journal of Pure and Applied Algebra};
also
\webprint{http:\dblslsh www\dt math\dt mq\dt edu\dt au\slsh
$\sim$mbatanin\slsh papers\dt html}{2001}{25}.}%

\subsection*{Definitions \ds{B} and \ds{L}}

Batanin gave his definition, together with an examination of $n=2$, in 
\bibentry{BatMGC}{%
M. A. Batanin,  
Monoidal globular categories as a natural environment for the theory of weak
$n$-categories, 
\jnl{Advances in Mathematics}{136}{1998}{no.~1, 39--103};
also 
\web{http:\dblslsh www\dt math\dt mq\dt edu\dt au\slsh $\sim$mbatanin\slsh
papers\dt html}.}%
Another account of it is
\bibentry{StrRMB}{%
Ross Street,
The role of Michael Batanin's monoidal globular categories,
\contrib{Higher category theory (Evanston, IL, 1997)}
{Contemporary Mathematics~\textbf{230}}
{AMS}{1998}{99--116};
also
\web{http:\dblslsh www\dt math\dt mq\dt edu\dt au\slsh $\sim$street}.}%
The definition of weak $n$-category which appears as~8.7 in \cyte{BatMGC} is
(I believe) what is here called definition \ds{B1}.  More precisely, let
\cat{O} be the category whose objects are (globular) operads on which
\emph{there exist} a contraction and a system of compositions, and whose maps
are just maps of operads.  What Batanin does is to construct an operad $K$
which is weakly initial in \cat{O}.  `Weakly initial' means that there is at
least one map from $K$ to any other object of \cat{O}, so this does not
determine $K$ up to isomorphism; one needs some further information.  But in
Remark~2 just before Definition~8.6, Batanin suggests that, once given the
appropriate extra structure, $K$ is initial in the category \fcat{OCS} of
operads \emph{equipped with} a contraction and a system of compositions,
which does determine $K$.  This is the approach taken in \ds{B1}.

A weak $n$-category according to \ds{B2} is (I believe) almost exactly what
Batanin calls a `weak $n$-categorical object in \emph{Span}' in his
Definition~8.6.  The only difference is my extra condition that the operad
$C$ is (in his terminology) \demph{normalized}: $C(0) \iso 1$.  Now $C(0)$ is
the set of operations in the operad which take a $0$-cell of an algebra and
turn it into another $0$-cell, so normality means that there are no such
operations except, trivially, the identity.  This seems reasonable in the
context of $n$-categories, since one expects to have operations for composing
$m$-cells only when $m\geq 1$.  The lack of normality in Batanin's version
ought to be harmless, since the contraction means that all the operations on
$0$-cells are in some sense equivalent to the identity operation, but it does
make the analysis of $n\leq 2$ a good deal messier.  (Note that the operad
$K$ is normalized, so any weak $\omega$-/$n$-category in the sense of \ds{B1}
is also one in the sense of \ds{B2}; the same goes for \ds{L1} and \ds{L2}.) 

My modification \ds{L1} of Batanin's definition first appeared in
\bibentry{SHDCT}{%
Tom Leinster,
Structures in higher-dimensional category theory,
\webprint{http:\dblslsh www\dt dpmms\dt cam\dt ac\dt uk\slsh
$\sim$leinster}{1998}{80},}%
but a more comprehensive and, I think, comprehensible account is in
\bibentry{OHDCT}{%
Tom Leinster,
Operads in higher-dimensional category theory,
Ph.D. thesis, University of Cambridge, 2000;
also
\eprint{math\dt CT\slsh 0011106}{2000}{viii + 127}.}%
(There, $C\otimes C'$ is called $C\of C'$, $C\cdot \dashbk$ is $T_C$, and
$P_C(\pi)$ is $P_\pi(C)$.)  \cyte{OHDCT} also contains a precise analysis of
\ds{L1} for $n\leq 2$, including proofs of~(a) the equivalence of the two
different categories of weak $n$-categories (for \emph{finite} $n$)
mentioned at the start of the analysis of $n\leq 2$ above, and~(b) the
equivalence of the category of unbiased bicategories and weak functors with
that of (classical) bicategories and weak functors.  Definition \ds{L2} has not
appeared before, and has just been added here for symmetry. 

The globular operads in \ds{B} and \ds{L} are called `$\omega$-operads in
\emph{Span}' by Batanin in \cyte{BatMGC}.  They are a special case of
\demph{generalized operads}, a family of higher-dimensional categorical
structures which are perhaps as interesting and applicable as $n$-categories
themselves.  Briefly, the theory goes as follows.  Given a monad $T$ on a
category $\cat{E}$, satisfying some natural conditions, one can define a
category of \demph{$T$-multicategories}.  For example, when $T$ is the
identity monad on the category \cat{E} of sets, a $T$-multicategory is just a
category, and when $T$ is the free-monoid monad on \Set, a $T$-multicategory
is just an ordinary multicategory (see `General Comments and History' above).
A \demph{$T$-operad} is a one-object $T$-multicategory, so in the first of
these examples it is a monoid and in the second it is an operad in the
original sense (but without symmetric or topological structure).  Now take
$T$ to be the free strict $\omega$-category monad on the category \cat{E} of
globular sets, as in \ds{B} and \ds{L}: a $T$-operad is then exactly a
globular operad.  Algebras for $T$-multicategories can be defined in the
general context, and again this notion specializes to the one in \ds{B} and
\ds{L}.

Generalized (operads and) multicategories were first put forward in
\bibentry{Bur}{%
Albert Burroni,  
$T$-cat\'egories (cat\'egories dans un triple),
\jnl{Cahiers de Topologie et G\'eom\'etrie
Diff\'erentielle}{12}{1971}{215--321}}%
and were twice rediscovered independently:
\bibentry{HerRM}{%
Claudio Hermida,  
Representable multicategories,
\jnl{Advances in Mathematics}{151}{2000}{no.~2, 164--225};
also 
\web{http:\dblslsh www\dt cs\dt math\dt ist\dt utl\dt pt\slsh s84\dt www\slsh
cs\slsh claudio\dt html},}%
\bibentry{GOM}{%
Tom Leinster,
General operads and multicategories,
\eprint{math\dt CT\slsh 9810053}{1997}{35}.}%
As far as I know, the notion of algebra for a $T$-multicategory only appears
in the third of these.  (\cyte{GOM} also appears, more or less, as Chapter~I
of \cyte{SHDCT} and Chapter~2 of \cyte{OHDCT}.) 

The difference between definitions \ds{L1} and \ds{B1} can be summarized by
saying that \ds{L1} takes \ds{B1}, dispenses with the notions of system of
compositions and \ds{B}-style contraction, and merges them into a single more
powerful notion of contraction.  A few more words on the difference are in
section~4.5 of \cyte{OHDCT}.  The operad $L$ canonically carries a
\ds{B}-style contraction and a system of compositions, so there is a
canonical map $K \go L$ of operads, and this induces a functor in the
opposite direction on the categories of algebras.  Hence every weak
$\omega$-/$n$-category in the sense of \ds{L1} gives rise canonically to one
in the sense of \ds{B1}.

\subsection*{Definition \lp}

This is the first time in print for definition \lp.  Once we have the
language of generalized multicategories (described in the previous section)
and the theory of free strict $\omega$-categories, it is very quickly stated.
My papers \cyte{OHDCT} and \cyte{SHDCT} (and to some extent \cyte{GOM}) cover
generalized multicategories and globular operads, but not specifically
globular multicategories.  The $1$-dimensional case, $1$-globular
multicategories, are the `$\mathbf{fc}$-multicategories' described briefly in
\bibentry{FCM}{%
Tom Leinster,
$\mathbf{fc}$-multicategories,
\eprint{math\dt CT\slsh 9903004}{1999}{8},}%
at a little more length in
\bibentry{GEC}{%
Tom Leinster,
Generalized enrichment of categories,
to appear in \emph{Journal of Pure and Applied Algebra},}%
and in detail in
\bibentry{GECM}{%
Tom Leinster,
Generalized enrichment for categories and multicategories,
\eprint{math\dt CT\slsh 9901139}{1999}{79}.}%

Logicians might like to view \lp\ through proof-theoretic spectacles,
substituting the word `proof' for `reason'.  They (and others) might also be
interested to read
\bibentry{MakAAC}{%
M. Makkai,
Avoiding the axiom of choice in general category theory,
\jnl{Journal of Pure and Applied Algebra}{108}{1996}{no.~2, 109--173};
also 
\web{http:\dblslsh www\dt math\dt mcgill\dt ca\slsh makkai}}%
in which Makkai defines anafunctors and anabicategories and discusses the
philosophical viewpoint which led him to them.  In the same vein, see also
Makkai's \cyte{MakTCF} and the remarks on `a composite' \emph{vs.}\ `the
composite' towards the end of the Introduction to the present paper. 

A weak $\omega$-/$n$-category in the sense of \ds{L2} (and so \ds{L1} too)
gives rise to one in the sense of \lp.  For just as `algebras' for a category
$C$ (functors $C \go \Set$) correspond one-to-one with discrete opfibrations
over $C$, via the so-called Gro\-th\-en\-dieck construction, so the same is
true in a suitable sense for globular multicategories.  This generalization
is explained in section~4.2 of \cyte{GOM}, section~I.3 of \cyte{SHDCT}, and
section~3.4 of \cyte{OHDCT} (any one of which will do, but they are listed in
increasing order of clarity).  What this means is that an algebra for a
globular operad gives rise to a globular multicategory (the domain of the
opfibration), and if the operad admits a contraction in the sense of \ds{L}
then the resulting multicategory is a weak $\omega$-category in the sense of
\lp.

Midway between \lp\ and \ds{J} is another possible definition of weak
$\omega$-category, which for various reasons I have not included here.  It
was presented in a talk,
\bibentry{LeiNQJ}{%
Tom Leinster,
Not quite Joyal's definition of $n$-category (a.k.a.~`algebraic nerves'),
seminar at \mbox{DPMMS}, Cambridge,
22 February 2001,}%
notes from which, in the $(2+2)$-page format of this paper, are available on
request.  The idea behind it can be traced back to Segal's formalization of
the notion of up-to-homotopy topological commutative monoid, \demph{special
$\Gamma$-spaces} (and their non-commutative counterparts, \demph{special
$\Delta$-spaces}).  The analogy is that just as Segal took the theory of
honest topological commutative monoids and did something to it to obtain an
up-to-homotopy version, so we take the theory of strict $n$-categories and do
something similar to obtain a weak version.  Segal's original paper is
\bibentry{SegCCT}{%
Graeme Segal, 
Categories and cohomology theories,
\jnl{Topology}{13}{1974}{293--312}.}%
A different generalization of his idea defines up-to-homotopy algebras for
any (classical) operad.  This is done at length in my
paper~\cyte{HAO}; or a much briefer explanation of the idea is
\bibentry{UTHM}{%
Tom Leinster,
Up-to-homotopy monoids,
\eprint{math\dt QA\slsh 9912084}{1999}{8}.}%

\subsection*{Definitions \ds{Si} and \ds{Ta}}

Tamsamani's original definition appeared in
\bibentry{TamNNC}{%
Zouhair Tamsamani,
Sur des notions de $n$-cat\'egorie et $n$-groupo\"{\i}de non strictes via des
ensembles multi-simpliciaux,  
\jnl{$K$-Theory}{16}{1999}{no.~1, 51--99};  
also
\epr{alg-geom\slsh 9512006}{1995}.}%
What I have called truncatability of a functor $\Delnop{r} \go \Set$ is
called `$r$-troncabilit\'e' by Tamsamani.  It is not immediately obvious that
the two conditions are equivalent, but a thoroughly mundane induction shows
that they are.  Other translations: my $1^p$ is his $I_p$, my $s$ and $t$ are
his $s$ and $b$, my $Q^{(m)}$ is his $T^m$, my $\pi^{(m)}$ is his $T^m$, my
internal equivalence of cells $x_1, x_2$ (as in the text of \ds{Ta}) is his
$(r-p)$-\'equivalence int\'erieure, and my external equivalence of functors
$\Delnop{r} \go \Set$ is his $r$-\'equivalence ext\'erieure.  His term for a
weak $n$-category is `$n$-nerf' or `$n$-cat\'egorie large'.  (`Large' has
nothing to do with large and small categories: it means broad or generous,
and can perhaps be translated here as `lax'; compare the English word
`largesse'.)

Tamsamani also offers a proof that his weak $2$-categories are essentially
the same as bicategories, but I believe that it is slightly flawed, in that
he has omitted a necessary axiom for the $2$-nerve of a bicategory (the last
bulleted item in `Definition \ds{Ta} for $n\leq 2$', starting
`$\iota_{uwx}^z$').  Without this, the constructed functor $\Delnop{2} \go
\Set$ will not necessarily be a weak $2$-category in the sense of \ds{Ta}.
(In this context, my $a$'s are his $x$'s, my $\alpha$'s are his $\lambda$'s,
and my $\iota$'s are his $\epsln$'s.)

Working with him in Toulouse, Simpson produced a simplified version of
Tamsamani's definition, which first appeared in
\bibentry{SimCMS}{%
Carlos Simpson,
A closed model structure for $n$-categories, internal
\textit{\underline{Hom}}, $n$-stacks and generalized Seifert-Van Kampen,  
\eprint{alg-geom\slsh 9704006}{1997}{69}.}%
He used the term `easy $n$-category' for his weak $n$-categories, and `easy
equivalence' for what is called a contractible map in \ds{Si}.

The simplification lies in the treatment of equivalences.  Weak
$1$-categories according to either \ds{Ta} or \ds{Si} are just categories,
but whereas a \ds{Ta}-style equivalence of weak $1$-categories is a functor
which is full, faithful and essentially surjective on objects (that is, an
ordinary equivalence of categories), an easy equivalence is a functor which
is full, faithful and \emph{genuinely} surjective on objects.  The latter
property of functors is expressible at a significantly more primitive
conceptual level than the former, since it is purely in terms of the
underlying directed graphs and has nothing to do with the actual category
structure.  For this reason, \ds{Si} is much shorter than \ds{Ta}.  (But to
develop the theory of weak $n$-categories we still need Tamsamani's more
general notion of equivalence; this is, for instance, the missing piece of
vocabulary referred to at the very end of `Definition \ds{Si} for $n\leq
2$'.)

As one would expect from this description, any easy equivalence (contractible
map) is an equivalence in the sense of \ds{Ta}.  So as long as it is true
that any weak $n$-category $\Delnop{n} \go \Set$ in the sense of \ds{Si} is
truncatable (which I cannot claim to have proved), it follows that any weak
$n$-category in the sense of \ds{Si} is also one in the sense of \ds{Ta}.

Following on from his definition, Tamsamani investigated homotopy
$n$-gr\-ou\-poids of spaces:
\bibentry{TamETH}{%
Zouhair Tamsamani,
Equivalence de la th\'eorie homotopique des $n$-groupo\"{\i}des
et celle des espaces topologiques $n$-tronqu\'es,
\eprint{alg-geom\slsh 9607010}{1996}{24}.}%
Numerous papers by Simpson, using a mixture of his definition and Tamsamani's
and largely in the language of Quillen model categories, push the theory of
weak $n$-categories further along:
\bibentry{SimLNC}{%
Carlos Simpson,
Limits in $n$-categories,
\eprint{alg-geom\slsh 9708010}{1997}{92},}%
\bibentry{SimHTS}{%
Carlos Simpson,
Homotopy types of strict 3-groupoids,
\eprint{math\dt CT\slsh 9810059}{1998}{29},}%
\bibentry{SimBBD}{%
Carlos Simpson,
On the Breen-Baez-Dolan stabilization hypothesis for Tamsamani's weak
$n$-categories,  
\eprint{math\dt CT\slsh 9810058}{1998}{36},}%
\bibentry{SimCMB}{%
Carlos Simpson,
Calculating maps between $n$-categories,
\eprint{math\dt CT\slsh 0009107}{2000}{13}.}%
Toen has also applied Tamsamani's definition, as in
\bibentry{ToenDTS1}{%
B. Toen,
Dualit\'e de Tannaka sup\'erieure I: structures monoidales,
Max-Planck-Institut preprint MPI-2000-57, 
\web{http:\dblslsh www\dt mpim-bonn\dt mpg\dt de},
2000, 71 pages}%
and
\bibentry{ToenNHC}{%
B. Toen,
Notes on higher categorical structures in topological quantum field theory,
\webprint{http:\dblslsh guests\dt mpim-bonn\dt mpg\dt de\slsh rosellen\slsh
etqft00\dt html}{2000}{14},}%
and the theory finds its way into some very grown-up mathematics in
\bibentry{SimAAH}{%
Carlos Simpson,
Algebraic aspects of higher nonabelian Hodge theory,
\eprint{math\dt AG\slsh 9902067}{1999}{186}.}%

The connection between categories and their nerves is covered briefly in one
of the new chapters of Mac Lane's book \cyte{MacCWM}; the more-or-less
original source is
\bibentry{SegCSS}{%
Graeme Segal, 
Classifying spaces and spectral sequences, 
\jnl{Institut des Hautes \'Etudes Scientifiques Publications Math\'ematiques}
{34}{1968}{105--112}.}%
Presumably the `Segal maps' are so named because of the prominent role they
play in Segal's paper \cyte{SegCCT} on loop spaces and homotopy-algebraic
structures.

The basic method by which a Simpson or Tamsamani weak 2-category gives rise
to a bicategory seems implicit in Segal's \cyte{SegCCT}, is made explicit in
section~3 of my \cyte{UTHM}, and is done in even more detail in section~3.3
of
\bibentry{HAO}{%
Tom Leinster,
Homotopy algebras for operads,
\eprint{math\dt QA\slsh 0002180}{2000}{101}.}%
(Actually, these last two papers only describe the method for monoidal
categories rather than bicategories in general, but there is no substantial
difference.)  There is also a discussion of the converse process in
section~4.4 of \cyte{HAO}, and the idea behind this is once more implicit in
the work of Segal.

\subsection*{Definition \ds{J}}

Joyal gave his definition in an unpublished note,
\bibentry{Joy}{%
A. Joyal,
Disks, duality and $\Theta$-categories,
preprint, \emph{c}.~1997, 6 pages.}%
There he defined a notion of weak $\omega$-category, which he called
`$\theta$-category'.  He also wrote a few informal words about structures
called $\theta^n$-categories, and how one could derive from them a definition
of weak $n$-category; but I was unable to interpret his meaning, and
consequently definition \ds{J} of weak $n$-category might not be what
he envisaged.

The term `disk' comes from the case where, in the notation of \ds{J}, $D_m$
is the closed $m$-dimensional unit disk ($=$ ball) in $\mathbb{R}^m$, $p_m$
is projection onto the first $(m-1)$ coordinates, and the order on the fibres
is given by the usual order on the real numbers.  The second bulleted
condition in the paragraph headed `Disks' holds at a point $d$ of $D_m$ if
and only if $d$ is on the boundary of $D_m$.  From another point of view,
this condition can be regarded as a form of exactness.

The handling of faces in \ds{J} is not necessarily equivalent to that in
\cyte{Joy}; again, I had trouble understanding the intended meaning and made
my own path.  In fact, Joyal works the duality discussed under $n\leq 2$ into
the definition itself, putting $\Theta = \scat{D}^\op$ and calling $\Theta$
the category of `Batanin cells' (for reasons suggested by
Figures~\ref{fig:op-comp-b} and~\ref{fig:disks}).  So he does not speak of
cofaces and cohorns in \scat{D}, but rather of faces and horns in $\Theta$.

Of the analyses of $n\leq 2$ for the ten definitions, that for \ds{J} is
probably the furthest from complete.  It appears to be the case that in a
weak $n$-category $A: \scat{D}_n \go \Set$, any cohorn $\Lambda^D_\phi \go A$
where $D$ has volume $>n$ has a unique filler.  (We know that this is true
when the dimension of $D$ is $n$.)  If this conjecture holds then we can
complete the proof (sketched in `$n\leq 2$') that any weak $2$-category gives
rise to a bicategory; for instance, applied to $T_{0,0,0}$ it tells us that
there is a canonical choice of associativity isomorphism, and applied to
$T_{0,0,0,0}$ it gives us the pentagon axiom.  However, I have not been able
to find a proof (or counterexample).

Introductory material on simplicial sets and horns can be found in, for
instance, Kamps and Porter's book \cyte{KP}.

The duality between the skeletal category $\Del$ of nonempty finite totally
ordered sets and the skeletal category \scat{I} of finite strict intervals
has been well-known for a long time.  Nevertheless, I have been unable to
trace the original reference, or even a text where it is explained
directly---except for Joyal's preprint \cyte{Joy}, which the reader may have
trouble obtaining.  Put briefly, the duality comes from mapping into the
2-element ordered set: if $k$ is a natural number then the set
$\Del([k],[1])$ naturally has the structure of an interval (isomorphic to
$\langle k\rangle$) and the set $\scat{I}(\langle k\rangle, \langle
0\rangle)$ naturally has the structure of a totally ordered set (isomorphic
to $[k]$).  This provides functors $\Del(\dashbk,[1]): \Delop \go \scat{I}$
and $\scat{I}(\dashbk, \langle 0\rangle): \scat{I}^\op \go \Del$ which are
mutually inverse, so $\Delop \iso \scat{I}$.

The higher duality has been the subject of detailed investigation by Makkai
and Zawadowski:
\bibentry{MZ}{%
Mihaly Makkai, Marek Zawadowski, 
Duality for simple $\omega$-categories and disks,
\jnl{Theory and Applications of Categories}{8}{2001}{114--243},}%
\bibentry{ZawDBD}{%
Marek Zawadowski,
Duality between disks and simple categories,
talk at 70th Peripatetic Seminar on Sheaves and Logic, Cambridge, 1999,}%
\bibentry{ZawDTD}{%
Marek Zawadowski,
A duality theorem on disks and simple $\omega$-categories, with applications to
weak higher-dimensional categories, 
talk at CT2000, Como, Italy, 2000.}%
(Slides and notes from Zawadowski's talks have the virtue of containing some
pictures absent in the published version.)  More on this duality and on the
relationship between definitions \ds{J} and \ds{B} is in
\bibentry{Ber}{%
Clemens Berger,
A cellular nerve for higher categories,
Universit\'e de Nice---Sophia Antipolis Pr\'epublication 602 (2000), 50 pages;
also 
\web{http:\dblslsh math\dt unice\dt fr\slsh $\sim$cberger}}%
(where a closed model category structure on $\ftrcat{\scat{D}}{\Set}$ is
also discussed) and in
\bibentry{BS}{%
Michael Batanin, Ross Street,
The universal property of the multitude of
trees, 
\jnl{Journal of Pure and Applied Algebra}{154}{2000}{no.~1-3, 3--13};
also
\web{http:\dblslsh www\dt math\dt mq\dt edu\dt au\slsh $\sim$mbatanin\slsh
papers\dt html}.}%
As mentioned above, there is another way to define weak $n$-category which
has strong connections to both \ds{J} and \ds{L}: \cyte{LeiNQJ}.

\subsection*{Definition \ds{St}}

Street proposed his definition of weak $\omega$-category in a very tentative
manner, in the final sentence of
\bibentry{StrAOS}{%
Ross Street,
The algebra of oriented simplexes,
\jnl{Journal of Pure and Applied Algebra}{49}{1987}{no.~3, 283--335}.}%
He did not explicitly formulate a notion of weak $n$-category for finite
$n$; this small addition is mine, as is the Variant at the end of the
section on $n\leq 2$.  

There is one minor but material difference, and a small number of cosmetic
differences, between Street's definition and \ds{St}.  The material
difference is that in a weak $\omega$-category as proposed in \cyte{StrAOS},
the only hollow $1$-cells are the degenerate ones.  One terminological
difference is that a pair $(A,H)$ is called a `simplicial set with
hollowness' in \cyte{StrAOS} only when~\bref{part:degen} and the
aforementioned condition on $1$-cells hold: so the term has a narrower
meaning there than here.  (Street informs me that he and Verity have used the
term `stratified simplicial set' for the same purpose, either with the two
conditions or without.)  Another is that he uses `$\omega$-category' in a
wider sense: his potentially have infinite-dimensional cells, and the
category of strict $\omega$-categories in the sense of the present paper is
denoted $\omega$-Cath.  Further translations: I say that a subset $S \sub
[m]$ is $k$-alternating where Street says that the set $[m]\without S$ is
`$k$-divided', and he calles a map $[l] \go [m]$ `$k$-monic' if its image is
a $k$-alternating subset of $[m]$.

The focus of \cyte{StrAOS} is actually on \emph{strict} $n$- and
$\omega$-categories.  To this end he considers the condition on $1$-cells
mentioned above, and conditions \bref{part:degen}--\bref{part:comp} of
\ds{St} with `unique' inserted before the word `filler'
in~\bref{part:filler}.  Having spent much of the paper constructing the nerve
of a strict $\omega$-category (this being a simplicial set with hollowness),
he conjectures that a given simplicial set with hollowness is the nerve of
some strict $\omega$-category if and only if all the conditions just
mentioned hold.  (The conjecture was, I believe, a result of joint work with
John Roberts.)  The necessity of these conditions was proved soon afterwards
in
\bibentry{StrFN}{%
Ross Street,
Fillers for nerves,
\contrib{Categorical algebra and its applications (Louvain-La-Neuve, 1987)}
{Lecture Notes in Mathematics~\textbf{1348}}
{Springer}{1988}{337--341}.}%
A proof of their sufficiency was supplied by Dominic Verity; this has not
appeared in print, but was presented at various seminars in Berkeley, Bangor
and Sydney around 1993.

It is entirely possible that most of the detailed work for $n\leq 2$ has
already been done by Duskin.  A short account of his work on this was
presented as
\bibentry{DusComo}{%
John W. Duskin, 
A simplicial-matrix approach to higher dimensional category theory,  
talk at CT2000, Como, Italy, 2000,}%
and a full-length version is in preparation:
\bibentry{DusSMA}{%
John W. Duskin, 
A simplicial-matrix approach to higher dimensional category theory I: 
nerves of bicategories,
preprint, 2001, 82 pages.}%
What Duskin does is to construct the nerve of any bicategory (this being a
simplicial set) and to give exact conditions saying which simplicial sets
arise in this way.  He moreover shows how to recover a bicategory from its
nerve.  Duskin does not deal explicitly with Street's conditions or his
notion of hollowness (although he does mention them); indeed, the results
just mentioned suggest that for $n=2$, the hollow structure on the nerve of a
bicategory is superfluous.  

The word `thin' has been used for the same purpose as `hollow', hence the
name \emph{T-complex}, as discussed in~III.2.26 and onwards in Kamps and
Porter's book \cyte{KP} (which also contains basic information on simplicial
sets and horn-filling). The original definition of T-complex was given by
M. K. Dakin in his 1975 Ph.D.\ thesis, published as
\bibentry{DakKCM}{%
M. K. Dakin, 
Kan complexes and multiple groupoid structures,
\jnl{Mathematical sketches (Esquisses Math\'ematiques)}{32}{1983}{xi+92 pages},
University of Amiens.}%
T-complexes are simplicial sets with hollowness satisfying
conditions~\bref{part:degen}--\bref{part:comp}, but with `admissible' dropped
and `filler' changed to `unique filler' in~\bref{part:filler}.  The dropping
of `admissible' means that the delicate orientation considerations of
Street's paper are ignored and any direction is as good as any
other---everything can be run backwards.  Thus, T-complexes are meant to be
like strict $\omega$-groupoids rather than strict $\omega$-categories.

\subsection*{Definition \ds{X}}

The story of \ds{X} is complicated.  Essentially it is a combination of the
ideas of Baez, Dolan, Hermida, Makkai and Power.  Baez and Dolan proposed a
definition of weak $n$-category, drawing on that of Street, in
\bibentry{BDHDA3}{%
John C. Baez, James Dolan,
Higher-dimensional algebra III: $n$-categories and the algebra of opetopes,
\jnl{Advances in Mathematics}{135}{1998}{no.~2, 145--206};
also
\epr{q-alg\slsh 9702014}{1997}.}%
An informal account is in section~4 of Baez's \cyte{BaezINC}.  In turn,
Hermida, Makkai and Power drew on the work of Baez and Dolan, producing a
modified version of Baez and Dolan's opetopic sets, which they called
multitopic sets.  (My use in \ds{X} of the former term rather than the
latter should not be interpreted as significant.)  Their original preprint
still seems to be available somewhere on the web:
\bibentry{HMPWHDpre}{%
Claudio Hermida, Michael Makkai, John Power,
On weak higher-dimensional categories, 
\webprint{http:\dblslsh fcs\dt math\dt sci\dt hokudai\dt ac\dt jp\slsh
doc\slsh info\slsh ncat\dt html}{1997}{104}}%
and is currently enjoying a journal serialization:
\bibentry{HMPWHD1}{%
Claudio Hermida, Michael Makkai, John Power,
On weak higher dimensional categories I: Part 1, 
\jnl{Journal of Pure and Applied Algebra}{154}{2000}{no.~1-3, 221--246},}%
\bibentry{HMPWHD2}{%
Claudio Hermida, Michael Makkai, John Power,
On weak higher-dimensional categories I---2, 
\jnl{Journal of Pure and Applied Algebra}{157}{2001}{no.~2-3, 247--277},}%
\bibentry{HMPWHD3}{%
Claudio Hermida, Michael Makkai, John Power,
On weak higher-dimensional categories I: third part,
to appear in 
\emph{Journal of Pure and Applied Algebra}.}%
A related paper with a somewhat different slant and in a much more elementary
style is 
\bibentry{HMPHDM}{%
Claudio Hermida, Michael Makkai, John Power,
Higher-dimensional multigraphs,
\contb{Thirteenth Annual IEEE Symposium on Logic in Computer Science
(Indianapolis, IN, 1998)}
{IEEE Computer Society, Los Alamitos, CA}{1998}{199--206}.}%
Hermida, Makkai and Power's original work did not go as far as an alternative
definition of weak $n$-category, although see the description below of
\cyte{MakMOC}.

I learned something near to definition \ds{X} from
\bibentry{Hy}{%
Martin Hyland,
Definition of lax $n$-category,
seminar at \mbox{DPMMS}, Cambridge, based on a conversation with John Power, 
18 June 1997.}%
Whether this is closer to the approach of Baez and Dolan or of Hermida,
Makkai and Power is hard to say.  The Baez-Dolan definition falls into two
parts: the definition of opetopic set, then the definition of universality.
Certainly the universality in \ds{X} is Baez and Dolan's, but the sketch of
the definition of opetopic set is very elementary, in contrast to the highly
involved definitions of opetopic/multitopic set given by both these groups of
authors.

Opetopic sets are, it is claimed in \cyte{BDHDA3}, just presheaves on a
certain category, the category of \emph{opetopes}.  (The situation can be
compared with that of simplicial sets, which are just presheaves on the
category $\Delta$.)  Multitopic sets are shown in \cyte{HMPWHDpre} to be
presheaves on a category of multitopes.  A third notion of opetope, going
(perhaps reprehensibly) by the same name, is given briefly in section~4.1 of
my \cyte{GOM}, and is laid out in more detail in Chapter~IV of my
\cyte{SHDCT}.  Roughly speaking, it is shown that all three notions are
equivalent in
\bibentry{CheROM}{%
Eugenia Cheng,
The relationship between the opetopic and multitopic approaches to weak
$n$-categories,
\webprint{http:\dblslsh www\dt dpmms\dt cam\dt ac\dt uk\slsh
$\sim$elgc2}{2000}{36}}%
(which compares Baez-Dolan's notion with Hermida-Makkai-Power's) and
\bibentry{CheEAT}{%
Eugenia Cheng,
Equivalence between approaches to the theory of opetopes,
\webprint{http:\dblslsh www\dt dpmms\dt cam\dt ac\dt uk\slsh
$\sim$elgc2}{2000}{36}}%
(which adds in my own).  More accurately, Cheng begins \cyte{CheROM} by
modifying Baez and Dolan's notion of operad; the effect of this is that the
symmetries present in Baez and Dolan's account are now handled much more
cleanly and naturally, especially when it comes to the crucial process of
`slicing'.  So this means that the Baez-Dolan opetopes are not necessarily
the same as the three equivalent kinds of opetope involved in Cheng's
result, and it remains to be seen whether they fit in.

Let us now turn from opetopic sets to universality.  The notion of liminality
does not appear in Baez and Dolan's paper, and is in some sense a substitute
for their notion of `balanced puncture niche'.  I made this change in order
to shorten the inductive definitions; it is just a rephrasing and has no
effect on the definition of universal cell.  The price to
be paid is that in isolation, liminality is probably a less meaningful
concept than that of balanced punctured niche.  

More on the formulation of universality can be found in
\bibentry{CheNUO}{%
Eugenia Cheng,
A notion of universality in the opetopic theory of $n$-categories,
\webprint{http:\dblslsh www\dt dpmms\dt cam\dt ac\dt uk\slsh
$\sim$elgc2}{2001}{12}.}%
Makkai appears to have hit upon a notion of `$\omega$-dimensional
universal properties', and thereby developed the definition of multitopic set
into a new definition of weak $\omega$-category:
\bibentry{MakMOC}{%
M. Makkai,
The multitopic $\omega$-category of all multitopic $\omega$-categories,
\webprint{http:\dblslsh mystic\dt biomed\dt mcgill\dt ca\slsh
M\_Makkai}{1999}{67}.}%
I do not, unfortunately, know enough about this to include an account here.
Nor have I included the definition of
\bibentry{LeiBMB}{%
Tom Leinster,
Batanin meets Baez and Dolan: yet more ways to define weak $n$-category,
seminar at \mbox{DPMMS}, Cambridge,
6 February 2001,}%
which uses opetopic shapes but an algebraic approach like that of \ds{L}.
This definition can be repeated for various other shapes, such as globular
(giving exactly~\ds{L}) and computads (which are like opetopes but with many
outputs as well as many inputs), and perhaps simplicial and even cubical.

Finally, the analysis of $n\leq 2$ has been done in a very precise way, in
\bibentry{CheEOC}{%
Eugenia Cheng,
Equivalence between the opetopic and classical approaches to bicategories,
\webprint{http:\dblslsh www\dt dpmms\dt cam\dt ac\dt uk\slsh
$\sim$elgc2}{2000}{68}.}%
This uses the notion of opetopic/multitopic set given by Cheng's modification
of Baez and Dolan, or by my opetopes, or by Hermida, Makkai and Power (for by
her equivalence result, all three notions are the same), together with the
Baez-Dolan notion of universality.  I am fairly confident that this gives the
same definition of weak $2$-category as is described in \ds{X} above.

\subsection*{Other Definitions of $n$-Category}

I have already mentioned several proposed definitions of weak $n$-category
which are not presented here.  My own \cyte{LeiNQJ} and \cyte{LeiBMB} are
missing.  The opetopic definitions---those related to definition \ds{X}---are
under-represented, as I have not given \emph{any} such definition in precise
terms; in particular, there is no exact presentation of Baez-Dolan's
definition~\cyte{BDHDA3}, of Cheng's modification of Baez-Dolan's definition
(\cyte{CheROM}, \cyte{CheEAT}), or of Makkai's definition \cyte{MakMOC}.

In the final stages of writing this I received a preprint,
\bibentry{MayOCA}{%
J. P. May, 
Operadic categories, $A_\infty$-categories and $n$-categories,
notes of a talk given at Morelia, Mexico on 25 May 2001,
10 pages,}%
containing another definition of weak $n$-category.  I have not had time to
assimilate this; nor have I yet digested the approach to weak $n$-categories
in
\bibentry{MiTs}{%
Hiroyuki Miyoshi, Toru Tsujishita,
Weak $\omega$-categories as $\omega$-hypergraphs,
\eprint{math\dt CT\slsh 0003137}{2000}{26},}%
\bibentry{HMT}{%
Akira Higuchi, Hiroyuki Miyoshi, Toru Tsujishita,
Higher dimensional hypercategories,
\eprint{math\dt CT\slsh 9907150}{1999}{25}.}%

\subsection*{Comparing Definitions}

It seems that not a great deal of rigorous work has been done on comparing the
proposed definitions, although there are plenty of informal ideas floating
about.  The papers that I know of are listed above under the appropriate
definitions.  

The `$n\leq 2$' sections show that there are many reasonable notions even of
weak $2$-category.  This is not diminished by restricting to one-object weak
$2$-categories, that is, monoidal categories.  So by examining and trying to
compare various possible notions of monoidal category, one can hope to get
some idea of what things will be like for weak $n$-categories in general.  A
proof of the equivalence of various `algebraic' or `definite' notions of
monoidal category is in
\bibentry{WMC}{%
Tom Leinster,
What's a monoidal category?,
poster at CT2000, Como, Italy, 2000,}%
and a similar but less general result is in Chapter~1 of my \cyte{OHDCT}
(actually stated for bicategories).  Hermida compares the indefinite with the
definite in his paper \cyte{HerRM} on representable multicategories, and a
different definite/indefinite comparison is in section~3 of my \cyte{UTHM} or
section~3.3 of my \cyte{HAO}.  (I use the terms `definite' and `algebraic' in
the sense of the Introduction.)

No-one who has seen the definition of tricategory given by Gordon, Power and
Street in \cyte{GPS} will take lightly the prospect of analysing the case
$n=3$.  However, it is worth pointing out an aspect of this definition less
well-known than its complexity: that it is not quite algebraic.  

In precise terms, what I mean by this is that the category whose objects are
tricategories and whose maps are strict maps of tricategories is not monadic
over the category of $3$-globular sets.  ($3$-globular sets are globular sets
as in the `Strict $n$-Categories' section of `Background', but with $m$ only
running from $0$ up to $3$.  So the graph structure of a tricategory is a
$3$-globular set.)  For whereas most of the definition of tricategory
consists of some data subject to some equations, a small part does not: in
items~(TD5) and~(TD6), it is stipulated that certain transformations of
bicategories are equivalences.  This is not an algebraic axiom; to make it
into one, we would have to add in as data a pseudo-inverse for each of these
equivalences, together with two invertible modifications witnessing the fact
that it is a pseudo-inverse, and then we would want to add more coherence
axioms (saying, amongst other things, that this data forms an \emph{adjoint}
equivalence).  The impact is that there is little chance of proving that the
category of weak $3$-categories (and strict maps) according to \ds{P},
\ds{B1}, \ds{L1}, or any other algebraic definition is equivalent to the
category of tricategories and strict maps.  This is in contrast to the
situation for $n=2$.

\subsection*{Related Areas}

I will be extremely brief here; as stated above, this is not meant to be a
survey of the literature.  However, there are two areas I feel it would be
inappropriate to omit.  Most of the references that follow are meant to
function as `meta-references', and are chosen for their comprehensive
bibliographies.

The first area is the Australian school of 2-dimensional algebra, a
representative of which is
\bibentry{BKP}{%
R. Blackwell, G. M. Kelly, A. J. Power,
Two-dimensional monad theory,
\jnl{Journal of Pure and Applied Algebra}{59}{1989}{no.~1, 1--41}.}%
The issues arising there merge into questions of coherence, one starting point
for which is the paper `On braidings, syllepses and symmetries' by
Sjoerd Crans:
\bibentry{CraBSS}{%
Sjoed Crans, 
On braidings, syllapses and symmetries,
\jnl{Cahiers de Topologie et G\'eom\'etrie Diff\'erentielle}{41}
{2000}{no.~1, 2--74};
also
\web{http:\dblslsh math\dt unice\dt fr\slsh $\sim$crans},}%
\bibentry{CraBSSerr}{%
S. Crans,  
Erratum: `On braidings, syllapses and symmetries', 
\jnl{Cahiers de Topologie et G\'eom\'etrie Diff\'erentielle}{41}
{2000}{no.~2, 156}.}%
More references for work in this area are to be found in Street's \cyte{StrCS}.

The second area is from algebraic topology: where higher-dimensional category
theorists want to take strict algebraic structures and weaken them, stable
homotopy theorists like to take strict topological-algebraic structures and do
them up to homotopy (in a more sensitive way than one might at first
imagine).  The two have much in common.  Various systematic ways of doing the
latter have been proposed, and some of these are listed on the last page of
text in my \cyte{HAO}.  Missing from that list is the method of
\bibentry{BatHCC}{%
Mikhail A. Batanin,
Homotopy coherent category theory and $A_\infty$-structures in monoidal
categories,
\jnl{Journal of Pure and Applied Algebra}{123}{1998}{no.~1-3, 67--103};
also
\web{http:\dblslsh www\dt math\dt mq\dt edu\dt au\slsh $\sim$mbatanin\slsh
papers\dt html}.}%
Another connection with homotopy theory and loop spaces is in 
\bibentry{BFSV}{%
C. Balteanu, Z. Fiedorowicz, R. Schw\"anzl, R. Vogt,
Iterated monoidal categories,
\eprint{math\dt AT\slsh 9808082}{1998}{55}.}%

Further references for these two areas and more can be found in the
`Introductory Texts' listed above.

\end{document}